\documentclass[12pt,a4paper]{amsart}

\makeindex

 \def\mmat #1,#2,#3,#4,{\text{\small\arraycolsep=3pt $
\begin{pmatrix}#1&#2\\#3&#4\end{pmatrix}$}}

\newcommand {\der} {{\partial}}
\newcommand {\one} {{\bf 1}}
\newcommand {\lra} {\leftrightarrow}
\newcommand {\llra} {\longleftrightarrow}
\newcommand {\Llra} {\Longleftrightarrow}
\newcommand {\Lra} {\Longrightarrow}
\newcommand{\del}{\partial}

\usepackage{array}
\usepackage{graphicx}
\usepackage[all]{xy}
\usepackage{tabularx,mdwtab}

\usepackage{hyperref}

\usepackage{stmaryrd}

\usepackage{DLdef1}
\makeatletter
\@ifundefined{AOwidehat}{%
  \let\AOwidehat\widehat \def\widehat{\protect\AOwidehat}%
  \let\AOwidetilde\widetilde \def\widetilde{\protect\AOwidetilde}}{}
\makeatother

\makeatletter
\theoremstyle{thm}
\@ifundefined{Theorem}{}{}
\@ifundefined{Statement}{}{}
\@ifundefined{Lemma}{}{}
\@ifundefined{Proposition}{}{}
\@ifundefined{Claim}{}{}
\@ifundefined{Corollary}{}{}
\theoremstyle{definition}
\@ifundefined{Remark}{\newtheorem*{Remark}{\secno Remark}}{}
\makeatother

\makeatletter
\@ifundefined{rank}{\DeclareMathOperator{\rank}{rank}}{}
\makeatother

\newcommand{\Size}{\text{Size}}

\newcommand {\fba}{{\mathfrak{ba}}}

\renewcommand {\ssbegin}[2][*]
 {\refstepcounter{subsection}%
\if#1*
\addcontentsline{toc}{subsection}{\thesubsection.\hskip 1pc #2}%
\else
\addcontentsline{toc}{subsection}{\thesubsection.\hskip 1pc #2. #1}%
\fi
 \def \secno {\gdef \secno {}{\ssecfont
\thesubsection.\hskip 2ex}%
 }%
 \begin{#2}}

\renewcommand {\sssbegin}[2][*]
 {\refstepcounter{subsubsection}
\if#1*
\addcontentsline{toc}{subsubsection}{\thesubsubsection.\hskip 1pc #2}%
\else
\addcontentsline{toc}{subsubsection}{\thesubsubsection.\hskip 1pc #2. #1}
\fi
 \def \secno {\gdef \secno {}{\ssecfont \thesubsubsection.\hskip 2ex}%
 }%
 \begin{#2}}

\renewcommand {\parbegin}[2][*]
 {\refstepcounter{paragraph}
\if#1*
\addcontentsline{toc}{paragraph}{\theparagraph.\hskip 1pc #2}%
\else
\addcontentsline{toc}{paragraph}{\theparagraph.\hskip 1pc #2. #1}
\fi
 \def \secno {\gdef \secno {}{\ssecfont \theparagraph.\hskip 2ex}%
 }%
 \begin{#2}}

\makeatletter

\newcommand\surtab{\small\centering\extrarowheight=1\p@ \def\arraystretch{1.4}\tabcolsep3.6\p@ \doublerulesep=0.09\p@%
 \vskip 3\p@\advance\belowcaptionskip3\p@ \mathsurround\z@ \hyphenpenalty\z@ \doublehyphendemerits\z@ }
\newcommand\tcaption[2][0pt]{%
 \def\@tempa{#2}%
 \addtocounter{table}{1}\surtab\edef\@currentlabel{\thetable}%
 \ifx\@tempa\@empty \hfill \textbf{\tablename}~\thetable\kern#1 \else\textbf{\tablename}~\thetable.\kern 1ex{\boldmath\bfseries #2}\fi
}

\makeatother

\begin{document}



\title[Maximal subalgebras of simple
vectorial Lie superalgebras]{Maximal graded
subalgebras of simple vectorial\\ Lie superalgebras with polynomial coefficients}

\author{Dimitry Leites${}^{a,*}$, Irina~Shchepochkina${}^b$}

\address{${}^*$The corresponding author\\
${}^a$Department of mathematics, Stockholm University, Albanov\"agen 28, SE-114 19 
Stockholm, Sweden; dimleites@gmail.com\\
${}^b$Independent University of Moscow,
B. Vlasievsky per., d. 11, RU-119 002 Moscow, Russia;
iparam@mail.ru}


\begin{abstract}
S. Lie was the first to try to classify maximal subalgebras of
simple vectorial (i.e., consisting of vector fields) Lie algebras
with polynomial coefficients. In Appendix it is shown that, as stated, the problem is wild: 
describing all subalgebras contains the classification of pairs of matrices up to 
simultaneous conjugation, so there are too many such subalgebras. Several researchers
distinguished various classes of maximal subalgebras important in applications.

Here, in simple (or in Cartan prolongations close to simple) vectorial Lie superalgebras
with polynomial coefficients considered with gradings associated
with Weisfeiler filtrations, we describe maximal simple (and close
to simple) graded subalgebras --- a superization of what Sophus Lie tackled.

Together with I. Shchepochkina's description of the maximal subalgebras of matrix Lie
superalgebras, see arXiv:hep-th/9702122, this paper solves the super and infinite-dimensional
version of Dynkin's problem --- the description of the maximal subalgebras of simple
finite-dimensional Lie algebras --- except for the 3 exceptional ambients, which are considered
elsewhere.
\end{abstract}


\makeatletter
\@namedef{subjclassname@2020}{\textup{2020} Mathematics Subject Classification}
\makeatother

\maketitle

\markboth{\itshape Dimitry Leites, \itshape Irina~Shchepochkina} {{\itshape Maximal graded subalgebras of simple
vectorial Lie superalgebras}}

\thispagestyle{empty}
\setcounter{tocdepth}{3}
\tableofcontents

\section{Introduction}

Hereafter, the ground field is~$\Cee$. For the details of classification of simple vectorial Lie superalgebras with polynomial coefficients, see \cite{LS1} and a~finite-dimensional case \cite{Ld}. In Section~\ref{Sspisok}, we extract from \cite{LS1} everything needed to make this article self-contained. 

\ssec{Setting of the problem and history} Any sub(super)algebra~of $\fgl(W)$ is said to be \textit{linear} or \textit{matrix} when a~ basis of $W$ is fixed.\index{Lie superalgebra, linear or matrix} 
In what follows, Lie
(super)algebras of vector fields will
be briefly called \textit{vectorial}.\index{Lie superalgebra, vectorial} In the main body of this article, we confine ourselves to vectorial Lie (super)algebras with \textbf{polynomial} coefficients; in Introduction and Discussions, we mention other coefficients (smooth, Laurent, formal power series, divided powers, \dots). 

Recall that a~given Lie
superalgebra~$\fg$ without proper ideals and of dimension $>1$ is
said to be \textit{simple}.\index{Lie superalgebra, simple}
Every infinite-dimensional $\Zee$-graded simple 
vectorial Lie superalgebra~$\fg$ has a~finite number of inequivalent invariant filtrations, 
each filtration corresponding to a~maximal subalgebra~of finite codimension, see \cite{LS1}. These
filtrations are the ones B.~Weisfeiler introduced in \cite{W}; since the gradings they induce
are what this whole article is about, here is what they are.

Let $\fL$ be a~vectorial Lie superalgebra~and let $\fL_{0}\subset\fL$ be a~maximal
subalgebra~of finite codimension containing no ideal of~$\fL$. Let $\fL_{-1}$ be a~minimal
$\fL_{0}$-invariant subspace strictly containing $\fL_{0}$ and, for $i\geq 1$, set
\begin{equation}\label{Wgen}
\fL_{-i-1}=[\fL_{-1}, \fL_{-i}]+\fL_{-i}, \qquad
\fL_{i}=\{D\in \fL_{i-1}\mid [D, \fL_{-1}]\subset \fL_{i-1}\}.
\end{equation}
What we get is a~decreasing filtration
\[
\fL=\fL_{-d}\supset \fL_{-d+1}\supset\dots\supset \fL_{0}\supset \fL_{1}\supset\cdots,
\]
called the \textit{Weisfeiler filtration},\index{Weisfeiler filtration, W-filtration} or
\textit{W-filtration}, of~$\fL$; the number $d$, finite for every algebra~considered here, is
its \textit{depth}.\index{Depth of a~filtration or a~grading} A~Lie superalgebra~taken together
with such a~filtration is said to be \textit{W-filtered},\index{Lie superalgebra, W-filtered} and
the associated graded Lie superalgebra
\[
\fg=\mathop{\oplus}_{i\geq -d}\ \fg_{i}, \text{~~where $\fg_{i}:=\fL_{i}/\fL_{i+1}$},
\]
is said to be \textit{W-graded};\index{Lie superalgebra, W-graded} its grading is called
a~\textit{W-grading}\index{Weisfeiler grading, W-grading} of~$\fg$. For vectorial Lie
superalgebras, the notion invariant under automorphisms is the filtration, not the grading; graded
algebras are, however, easier to work with, and every ambient of this article is W-graded.

Three properties of a~W-grading follow from \eqref{Wgen} alone, and are used throughout:
\begin{equation}\label{Wprop}
\begin{minipage}[l]{14cm}
$1)$ The $\fg_{0}$-module $\fg_{-1}$ is irreducible;\\
$2)$ $\fg_{-i-1}=[\fg_{-1}, \fg_{-i}]$ for $1\leq i<d$; in particular, $\fg_{-}$ is generated by
$\fg_{-1}$;\\
$3)$ $\fg$ is transitive, see Subsec.~\ref{2.6}, and it suffices to test the bracket with
$\fg_{-1}$:\\
\phantom{.}\hskip 0.4cm if $x\in \fg_{i}$ with $i\geq 0$ and $[x, \fg_{-1}]=0$, then $x=0$.
\end{minipage}
\end{equation}
Indeed, $1)$ holds because $\fL_{-1}$ is minimal among the $\fL_{0}$-invariant subspaces strictly
containing $\fL_{0}$, so the $\fL_{0}$-module $\fL_{-1}/\fL_{0}$ is irreducible, whereas
$[\fL_{1}, \fL_{-1}]\subset\fL_{0}$ by \eqref{Wgen}, so the action is through
$\fg_{0}=\fL_{0}/\fL_{1}$; item $2)$ is the first of the formulas \eqref{Wgen} read modulo
$\fL_{-i}$; and item $3)$ is the second of them read with $i+1$ in place of~$i$.
 
The following examples exhaust all simple vectorial Lie
algebras, and their non-trivial central extensions:

\underline{$\fvect(n)=\fder\, \Cee[x]$}, where $x=(x_{1}, \dots ,
x_n)$.

\underline{$\fsvect(n)=\{D\in \fvect(n) \mid \Div D=0\}=\{D\in \fvect(n) \mid L_D(\vvol)=0\}$},
where $L_D$ is the Lie derivative along the vector field $D$ and $\vvol$ is the volume
element (in fixed coordinates).

\underline{$\fk(2n+1)$}, the Lie algebra~of \textit{contact vector
fields}; it is the subalgebra~in $\fvect(2n+1)$ that preserves
the distribution given by Pfaff equation $\alpha(X) =0$ for $X\in
\fvect(2n+1)$, where
\[
\alpha= dt+\sum(p_idq_i-q_idp_i).
\]
The elements of $\fk(2n+1)$ depend on just one functional
parameter; these elements are denoted
$K_{f}\in\fk(2n+1)$, where $f\in\Cee[t, p, q]$ for $p=(p_1, \dots,
p_n)$ and ${q=(q_1, \dots, q_n)}$. One directly verifies that
$L_{K_{f}}(\alpha)= 2\pderf{f}{t}\alpha$.

\underline{$\fpo(2n)\subset \fk(2n+1)$}, the \textit{Poisson
algebra},\index{Poisson
algebra} spanned by contact vector fields $K_{f}\in\fk(2n+1)$
such that $\pderf{f}{t}=0$, i.e., $L_{K_{f}}(\alpha)=0$.

\underline{$\fh(2n)=\fpo(2n)/\fc$, where $\fc$ is the center}, the
Lie algebra~of Hamiltonian vector fields. \vskip 0.2cm

Sophus Lie was the first to consider the problem ``classify all
maximal subalgebras of $\fvect(n)$", and a version: ``classify all
finite-dimensional subalgebras of $\fvect(n)$", see 
his collected works \cite{LieT}.\footnote{For the modern mathematician, unless one is a~
German-reading person, Lie's results are not readily
accessible, so we cite an English translation \cite{LieT} of Lie's works. Aside: observe that, in the 1970s,
M.~Semenov--Tian-Shansky skimmed through this volume~ 3 and saw
that Lie introduced the form
rediscovered by F.A.Berezin, later clearly described and popularized by A.A.Kirillov-p\`ere. Definitely, still more discoveries await any historian
of science capable to read German.}

In the list of finite-dimensional subalgebras of $\fvect(n)$
that Lie described, not only the abundance of Fraktur font stymied
the reader. Lie's result was so messy for $n=3$ that Lie himself refrained from
publishing it and confined to $n=1$ and 2, see an elucidation in \cite{ABFM} and references therein, see also \cite{CPT}; it was clear that to
advance as $n$ grows will be difficult, to say the least: the
problem might turn wild.  

It does, at any rate in the following closely related sense. Recall that a~classification
problem is said to be \textbf{wild} if it contains
the classification of pairs of matrices up to simultaneous conjugation, see \cite{Dr}. The
classification, up to conjugation, of the Lie subalgebras of $\fgl(n)$ for all $n$ is wild;
already the $11$-dimensional solvable ones --- a~torus plus a~$5$-dimensional abelian ideal
with zero product --- suffice. So is, therefore, the classification of the subalgebras
$\fr\subset\fvect(0|n)$ such that $\pr_{-1}(\fr)=\fvect_{-1}(0|n)$, considered up to
$\Aut\fvect(0|n)$; for the proofs, see \ref{Awild}, where they serve to show that
``describe primitive Lie superalgebras\rq\rq\ is a~wild problem. This is one more reason why the reasonable
question is not ``describe \textit{all} subalgebras\rq\rq\ but ``describe the \textbf{maximal}
ones\rq\rq.

The problem Lie considered has, however, several reasonable
subproblems. By ``reasonable" we mean either having a~neat answer,
or with obvious and interesting applications, or both. Even if
applications are not immediately seen, the beauty of the answer
justifies studying the problem: eventually, answers of this kind find their
applications. In particular, there are very few
maximal subalgebras in $\fvect(n)$, which are \textit{simple} or ``relatives" of a simple, i.e., its (usually nontrivial)
central extensions, algebra of derivations, and iterated
combinations thereof.\index{Lie (super)algebra, relative of}

As far as we know, the next, historically speaking, who, after Lie,
considered finite-dimensional subalgebras of vectorial Lie algebras
was V.~I.~Ogievetsky \cite{Og,BOg}. He did not bother to \textit{classify
all} subalgebras, but he gave nontrivial physical
interpretations to his simple (here: mathematically obvious)
examples of simple (here: no nontrivial ideals) maximal subalgebras.

Turbiner returned to Lie's problem for needs of ``quasi-exactly
solvable differential equations", see \cite{T1}. Inspired by Turbiner's
ideas, Post \cite{P} distinguished a~class of (not necessarily
simple) finite-dimensional subalgebras of $\fvect(n)$ and considered~$n\leq 4$. Post observed that though similar results can, in principle,
be obtained by other means (e.g., Vinberg's general, but implicit,
results cited in \cite{P}), his technique is rather elementary,
explicit, and very graphic. It is not Post's fault that for $n>4$,
the problem, though tame, loses attraction, since the answer, even
for $n=5$, requires several pages to be documented. 

Now, it is time to recall that a~given Lie
superalgebra~$\fg$ is \textit{semi-simple}\index{Lie superalgebra, semi-simple}
if its radical is zero. (For a description of semi-simple Lie superalgebras, see \cite{LS1}.) We say that $\fh$ is \textit{almost simple}\index{Lie superalgebra, almost simple} if it
can be sandwiched (non-strictly) between a~simple Lie superalgebra
$\fs$ and the Lie superalgebra~$\fder(\fs)$ of derivations of $\fs$,
i.e., ${\fs\subset\fh\subset\fder(\fs)}$. A~subalgebra~is said to be \textit{semi-almost simple}\index{Lie superalgebra, semi-almost simple}
 if it is a~direct sum of almost simple algebras. 

Hereafter, in this article, speaking of Lie (super)algebras,
\textbf{the adjective ``simple\rq\rq\ is usually broadened, by abuse of notation, to include also
``relatives" of simple algebras}. We usually drop the quotation marks in these cases. We consider only maximal \textbf{graded} subalgebras in ambient Lie superalgebras; this is a~vital restriction: for example, $\fk(2n+1|2n+2)$ is maximal in $\fsvect(2n+1|2n+2)$, but non-graded. To list all maximal subalgebras is more difficult than the following problem we solve here:

\be\label{PROB}
\begin{minipage}[l]{14cm}
\textbf{In ``simple\rq\rq\ 
vectorial Lie superalgebras $\fg=\oplus_{i\geq -d}\ \fg_i$, such as Cartan prolongs of $\fg_-:=\oplus_{i<0}\ \fg_i$, classify maximal simple, semi-simple, almost simple, and semi-almost 
simple graded
Lie subsuperalgebras $\fu= \mathop{\oplus}_{i\geq -d}\ (\fu\cap \fg_{i})$}.
\end{minipage}
\ee

\ssec{The six cases} 
Recall that the depth $d$ of any simple W-graded vectorial Lie superalgebra $\fg$ is 
equal to 1 or 2 or~3. \textbf{For each of these three depths we give the complete
answer}, and the depth changes nothing: as \S\ref{Sdepthred} shows, for a~\textbf{maximal}
graded subalgebra~the components $\fu_{-j}$ with $j\geq2$ are computed from
$\fu_{-1}$ and $\fu_{\geq0}$, not chosen, so the data to run over are the depth-$1$
data and the depth is not a~parameter of the answer.

What is incomplete is one \textbf{case}, not one depth: in the case $6^o$
of~\eqref{5cases}, and for the infinite-dimensional ambients only, we give partial
results and formulate conjectures and open problems. These are, in the order in which
they occur: the caveat after~\eqref{5cases}, saying what we observe and what we do not
prove; item a\textup{)} of Subsec.~\ref{ssUMNpart}, saying which ingredient is missing,
namely a~counterpart of \cite{ShM} for the vectorial series; item $2)$ of
Remark~\ref{Rdepthnot}, saying why the reduction of the depth does not reach this case;
and \cite{LOdd} for the three exceptional finite-dimensional ambients, where the case
is \textbf{not} vacuous. The reason the case $6^o$ escapes is not its depth but its
$\fu_{-1}$: being vacuous at the depth $1$, it is the one case for which the depth-$1$
machinery offers no pattern to reduce to.

Observe that, whatever the depth of $\fg$, the set of
graded subalgebras $\fu \subset \fg$, where $\fg$ is either simple or a~Cartan prolong of $\fg_{-}$ with irreducible
$\fg_0$-action on $\fg_{-1}$ splits into the following
non-intersecting subsets:
\begin{equation}\label{5cases}
\begin{minipage}[l]{14cm}
$1^o$. $\fu_{-1}=\fg_{-1}$ and $\fu_{0}=\fg_{0}$; $\fu_{1}\neq 0$;\\
$2^o$. $\fu_{-1}=\fg_{-1}$ and $\fu_{0}=\fg_{0}$; $\fu_{1}=0$;\\
$3^o$. $\fu_{-1}=\fg_{-1}$ and $\fu_{0}\neq \fg_{0}$;\\
$4^o$. $0 \subsetneq \fu_{-1}\subsetneq \fg_{-1}$;\\
$5^o$. $\fu_{-}=0$, where $\fu_-= \mathop{\oplus}_{ i<0}\ (\fu\cap
\fg_{i})$;\\
$6^o$. $\fu_{-1}=0$, but $\fu_{-}\neq 0$.
\end{minipage}
\end{equation}
The case $6^o$ is vacuous for ambients of depth $d=1$; besides, it occurs only in
finite dimension, where it is not vacuous: in each of the three exceptional
finite-dimensional ambients $\fosp_a(4|2)$, $\fag(2)$ and $\fab(3)$ taken with its contact
grading, the even part $\fg_\ev$ is a~maximal subalgebra of the case $6^o$, see
\cite{LOdd}.
A~caveat: ``Occurs only in finite dimension\rq\rq\ is what we observe: for none of the
infinite-dimensional ambients of this paper did we find a~maximal subalgebra of the case $6^o$, but we have no proof that there is none.

Let us immediately exclude obvious subcases from future discussions.

\sssbegin[Obvious subcases of the cases \eqref{5cases}]{Lemma}[Obvious subcases of the cases \eqref{5cases}] \label{Lobvious} In case \textup{$2^o$}, $\fu= \oplus
_{i\leq 0}\ \fg_{i}$ is maximal in~$\fg$, if and only if the action of $\fg_{0}$ on $\fg_{1}$ is irreducible and $\fg_{1}$ generates the whole  part with positive components.
In case \textup{$5^o$}, $\fu$ is a~maximal subsuperalgebra
of $\fg$ if and only if $\fu=\oplus_{i\geq 0}\ \fg_{i}$. \end{Lemma}

\begin{proof} The necessity is obvious. Consider sufficiency. In the case $2^o$, let $\fw\supsetneq\fu$ be a~graded
subalgebra. Since $\fg_{-1}\subseteq\fw$ and $\fg$ is transitive by \eqref{Wprop}, a~non-zero
component $\fw_k$ with $k>0$ gives $[\fg_{-1}, \fw_k]\neq0$, and bracketing
lowers the degree by one at a~time, so $\fw_1\neq0$. Being a~$\fg_0$-submodule
of the irreducible $\fg_1$, the component $\fw_1$ is the whole of $\fg_1$; and
$\fg_1$ generates $\mathop{\oplus}_{i>0}\fg_i$, whence $\fw=\fg$.

In the case $5^o$, the condition $\fu_-=0$ means $\fu\subseteq
\mathop{\oplus}_{i\geq 0}\fg_{i}$, which is a~proper subalgebra~of $\fg$, so
a~maximal $\fu$ of this case is that subalgebra. Conversely,
$\mathop{\oplus}_{i\geq 0}\fg_{i}$ is maximal among the graded
subalgebras~whose $(-1)$st component is non-zero: such a~$\fw$ has $\fw_{-1}$
a~non-zero $\fg_0$-submodule of $\fg_{-1}$, which is irreducible by item $1)$
of \eqref{Wprop}, so $\fw_{-1}=\fg_{-1}$; and the negative part of a~W-grading
is generated by its $(-1)$st component by item $2)$ of \eqref{Wprop}, whence
$\fw\supseteq\fg_-$ and $\fw=\fg$. \end{proof}

\ssec{Our results} For each depth, 1 to 3, of all possible W-gradings of simple vectorial Lie superalgebras we describe the maximal
elements in each of the cases~\eqref{5cases}. One of the problems we
encounter here is that the description of various cases is almost as
long as their proof. For the convenience of the reader, we gathered
the results in Tables \eqref{table1}--\eqref{table4}. 

We also listed, in passing, not claiming classification, several infinite-dimensional simple maximal subalgebras. Some of these examples naturally appear in various problems (e.g., see \cite{BM, BO, BTW} and references therein).

Most of the results here, due to I.~Shchepochkina, were preprinted in Report no.~${32/1988-15}$ of Department of Mathematics of Stockholm University. Here we sharpen formulations, expound the proofs, update references and correct partial results published since 1988, e.g., \cite{BM, BLM}.

\ssec{Related results and open problems}\label{ssOpen} 1) Since some simple vectorial Lie
\textbf{super}algebras are finite-dimen\-si\-o\-nal, the problem we 
consider here (a~generalization of problems considered by Lie, Ogievetsky,
Turbiner, Sokolov) coincides (for the finite-dimensional
case) with superization of a~problem considered by
Dynkin: ``classify maximal subalgebras of simple finite-dimensional Lie algebras". For the classification of maximal subalgebras
of finite-dimensional matrix Lie superalgebras, see
\cite{ShM}; we will need some of the results obtained there, and so we 
recall them in due course.

For a classification of maximal subalgebras in the exceptional simple Lie
algebras with Cartan matrix over fields of ``good" characteristic
$p>0$, see \cite{PS}.

That the unrestricted problem --- describe \textbf{all} subalgebras --- is hopeless is
now a~theorem and not just a~feeling: the classification, up to conjugation, of the Lie
subalgebras of $\fgl(n)$ for all $n$ contains the classification of pairs of matrices up to
simultaneous conjugation, and hence is \textit{wild} in the technical sense; already the
$11$-dimensional solvable subalgebras --- a~torus plus a~$5$-dimensional abelian ideal with
zero product --- suffice for that. The same construction shows that the finite-dimensional
Lie superalgebras $\fvect_{-1}(\xi)\oplus\fu$, where
$\xi=(\xi_1,\dots, \xi_n)$ and $\fu\subset\fgl(n)$ is the ``linear" part of a~
primitive Lie superalgebra, are wild to classify, see \ref{Awild}. This is why
``describe the \textbf{maximal} subalgebras" is a~ reasonable problem (solved here).

2) For subalgebras inhomogeneous with respect to \textit{parity} (called
\textit{Volichenko algebras}),\index{Volichenko algebra} see \cite{LSe, LSe2}. 

3) There were published several papers describing
maximal sub(super)algebras in certain simple Lie
(super)algebras of vector fields over fields of characteristic
$p>0$, see, e.g., \cite{BLLM, LW, Me}. Unlike \cite{Mu}, these papers do not complete the job, and
the gap is in the same place in both \cite{BLLM} and \cite{LW}: the maximal \textbf{irreducible}
graded subalgebras are not classified but only \textit{reduced} to the classification of the
maximal irreducible subalgebras of $\fgl(m|n)$, $\fsl(m|n)$ and $\fosp(m|n)$, resp. of $\fpe(n)$,
i.e., to the modular counterpart of \cite{ShM}, which is not available; and only the standard
gradings of the four series are considered, whereas over $\Cee$ it is the regradings, see table
\eqref{1.4b}, that produce a~good half of the answers: just compare with our tables.

4) For an algorithm allowing to embed every Lie superalgebra into the Lie superalgebra of differential operators with polynomial coefficients, 
see \cite{ShVf}, where an embedding of the Lie algebra~$\fg(2)$ is given as an illustration. 
In \cite{KST}, the embedding of~$\fg(2)$ is superized. 




5) 
V.~Sokolov \cite{S} considered maximal finite-dimensional
subalgebras in $\fk(1)=\fvect(1)$ and in $\fk(3)$; he associated
interesting differential equations with simple
subalgebras (and with their contractions), see also \cite{MSS}. In
\cite{SS}, the list of simple finite-dimensional and loop
subalgebras in the stringy Lie algebra $\fk^L(5)$, where the subscript $L$ indicates that
the contact Lie algebra~is realized on the space of \textit{Laurent}
polynomials, is announced; no related
differential equations were indicated.

6) For an algorithm allowing to embed every Lie superalgebra into the Lie superalgebra of differential operators with rational coefficients, see \cite{Gel}. It is unclear how this embedding was constructed, same question applies to the embedding considered in \cite{Po1}. 

7) Maximal graded solvable subalgebras of simple finite-dimensional vectorial Lie superalgebras are classified in \cite{Shch}.

\textbf{Open problems}.

1) Describe maximal subalgebras in the 3 exceptional simple finite-dimensional ambients:
$\fag(2)$, $\fab(3)$ and $\fosp_a(4|2)$ thus completing superization of Dynkin's result \cite{Dy}, see \cite{ShM}.
For the maximal \textbf{graded} subalgebras, in every $\Zee$-grading to which our
method applies, the answer is now given in \cite{LOdd}; what remains
open is listed in \cite{LOdd} as well.

2) Describe maximal inhomogeneous subalgebras of
simple $\Zee$-graded vectorial Lie superalgebras.

3) Describe embeddings of the maximal sub(super)algebras into $\fpo(2n|m)$ and $\fh(2n|m)$ yielding spinor-oscillator representations, and into $\fb_{a, b}(n)$ and $\fle(n)$ yielding the ``odd" versions of spinor-oscillator representations.

4) Find differential equations with big (super)algebra of symmetries, cf. \cite{The}.

5) Is there a~version of Post's algorithm \cite{P} for Lie superalgebras? 
For what values of $n$ and $m$, the list of
finite-dimensional subsuperalgebras in Lie superalgebra~$\fvect(m|n)$
is no less handable than Post's list?

6) Poletaeva~considered embedding of $\fosp_a(4|2)$ into a~stringy (with Laurent polynomials as coefficients) vectorial Lie superalgebra, see \cite{Po1}; this is not a~maximal subalgebra. Describe the maximal subalgebras containing~$\fosp_a(4|2)$.

\section{Basics. The list of ambient simple Lie superalgebras}\label{Sspisok} 

In this section, we recall most needed notation and facts from \cite{LS1}. Let the trivial central extension of the Lie (super)algebra $\fg$ be denoted by\index{$\fc\fg$ or $\fc(\fg)$, the trivial central extension of $\fg$} $\fc\fg$ or $\fc(\fg)$; the (sub)algebra of $\fder(\fg)$ containing only one outer derivation of degree 0 in a certain $\Zee$-grading of $\fg$ be denoted by $\fd(\fg)$.\index{$\fd(\fg)$, the (sub)algebra of $\fder(\fg)$ containing only outer degree 0 derivation of $\fg$}

\ssec{Linear  Lie (super)algebras}\label{SS:2.3} A
\textit{superspace} is a~$\Zee /2$-graded space 
$V=V_{\ev}\oplus V_{\od }$, where $\ev,\od\in\Zee/2$. A
\textit{superalgebra} is a~superspace with multiplication.

By $\Pi (V)$\index{$\Pi (-)$, the reversal of parity functor} we denote another copy of
the same superspace $V$, but with the inverted parity, \textit{i.e.,} $(\Pi(V))_{\bar
i}= V_{\bar i+\od }$. We call non-zero elements of $V_{\ev}$ \textit{even} (resp., those of $V_{\od}$ \textit{odd}); we denote the parity function by $p$ and write $p(v)=i$ if $v\in V_i$ and~$v\neq 0$.

Let $a=\dim V_{\ev}$, $b=\dim V_{\od
}$. The \textit{superdimension}\index{Superdimension} of $V$ is $\sdim
V:=a+b\eps $, where~$\eps ^2=1$. (Usually, the superdimension is written as a~pair $(a,b)$ or $a|b$;
this notation obscures the fact that $\sdim (V\otimes W)=\sdim V\cdot \sdim
W$.) 
 
Let $\Size=(p_1, \dots, p_{|\Size|})$ be an ordered
collection of parities of the basis vectors of a superspace $V$;
The \textit{general linear} Lie
superalgebra~of all supermatrices\index{Supermatrix} of size $\Size$ corresponding to operators in the superspace $V=V_{\bar 0}\oplus V_{\bar 1}$ is denoted by
$\fgl(\Size)$.\index{$\fgl(\Size)$} For one of the simplest formats, often called \textit{standard},\index{Format of the basis or supermatrix} $\fgl(\ev,
\dots, \ev, \od, \dots, \od)$ is abbreviated to $\fgl(\dim V_{\bar
0}|\dim V_{\bar 1})$. Any non-zero supermatrix in $\fgl(\Size)$ can
be uniquely expressed as the sum of its even and odd parts; in the
standard format this is the following block expression; on non-zero summands the parity is defined:
\[
\begin{pmatrix}A&B\\ C&D\end{pmatrix}=\begin{pmatrix}A&0\\
0&D\end{pmatrix}+\begin{pmatrix}0&B\\ C&0\end{pmatrix},\quad
p\left(\begin{pmatrix}A&0\\
0&D\end{pmatrix}\right)=\ev, \; p\left(\begin{pmatrix}0&B\\
C&0\end{pmatrix}\right)=\od.
\]

The \textit{trace} on any the Lie (super)algebra $\fg$ is a linear functional that vanishes on the derived algebra~$\fg'$. For example, the \textit{supertrace} \index{Trace, supertrace} is the map $\fgl (\Size)\tto \Cee$,
$(A_{ij})\longmapsto \sum (-1)^{p_{i}}A_{ii}$. Since $\str [x, y]=0$, the subsuperspace of supertraceless
matrices constitutes the \textit{special linear} subsuperalgebra~$\fsl(\Size)$.

There are, however, at least two super versions of $\fgl(n)$, not
one, see \cite[Ch1, Ch.7]{Lsos}. The other version is called the \textit{queer}\index{$\fq(n)$, Lie superalgebra, queer}
Lie superalgebra~$\fq(n)$; it 
preserves the
complex structure given by an \textit{odd} operator $J$, i.e.,
$\fq(n)$ is the centralizer $C(J)$ of $J$:
\[
\fq(n)=C(J)=\{X\in\fgl(n|n)\mid [X, J]=0 \}, \text{ where }
J^2=-\id.
\]
It is clear that by a~change of basis we can reduce $J$ to the form
$J_{2n}=\mat {0&1_n\\-1_n&0}$ in the standard
format and then $\fq(n)$ takes the form
\begin{equation}\label{q}
\fq(n)=\left \{(A,B):=\mat {A&B\\B&A}\mid A, B\in\fgl(n)\right\}.
\end{equation}

On $\fq(n)$, the \textit{queertrace}\index{Queertrace} is defined: $\qtr\colon (A,B)\longmapsto \tr B$. Denote by $\fsq(n)$ the Lie superalgebra
of \textit{queertraceless} matrices; set $\fp\fsq(n):=\fsq(n)/\Kee 1_{2n}$.

To the linear map $F: V\tto W$ there corresponds the
dual map $F^*:W^*\tto V^*$ of the dual superspaces. In bases
consisting of the homogeneous vectors $v_{i}\in V$ and $w_{j}\in W$, the formula
$F(v_{j})=\mathop{\sum}_{i}v_{i}A_{ij}$ assigns to $F$ the
supermatrix $A$ with elements~$A_{ij}$. In the dual bases, the \textit{supertransposed}\index{Supertransposition}
matrix $A^{st}$\index{${}^{st}$, supertransposition} corresponds to $F^*$:
\[
(A^{st})_{ij}=(-1)^{(p(v_{i})+p(w_{j}))p(w_{j})}A_{ji}.
\]

\ssec{Lie superalgebras that preserve bilinear forms: two types. Gram matrices}\label{sssBilMatr} The supermatrices $X\in\fgl(\Size)$ such that
\[
X^{st}B+(-1)^{p(X)p(B)}BX=0\quad \text{for a~ homogeneous matrix
$B\in\fgl(\Size)$}
\]
constitute the Lie superalgebra~$\faut (B)$ that preserves the
bilinear form $\cB$ on $V$ with Gram matrix~$B$.

In order to identify a~bilinear form $\cB\in\Bil(V, W)$ with an element of $\Hom(V, W^*)$, the Gram matrix $B$ of the form $\cB$ should be defined (for the case $V=W$ sufficient for us) by the formula 
(see \cite[Ch.1]{Lsos})
\be\label{matrBil}
B_{ij}=(-1)^{p(B)p(v_i)}\cB(v_i,v_j)\text{~~for the basis vectors $v_i\in V$.}
\ee

Let \textit{upsetting} of bilinear forms\index{Upsetting of bilinear forms}
$u\colon\Bil (V, W)\tto\Bil(W, V)$ be defined by the formula
\[
u(\cB)(w,v):= (-1)^{p(v)p(w)} \cB(v,w)\text{~~for every $v\in V, \ w\in W$}.
\]
The homogeneous form $\cB$ is \textit{symmetric} (resp., \textit{antisymmetric}) if $u(B)=B$ (resp., $u(B)=-B$). 

Observe
that \textbf{the passage from $V$ to $\Pi (V)$ turns every symmetric
form $\cB$ on $V$ into the anti-symmetric one $\cB^\Pi$ on $\Pi (V)$ and every anti-symmetric form $\cB$on $V$ into the symmetric form $\cB^\Pi$ on $\Pi (V)$ by setting}
\[
\text{$\cB^\Pi(\Pi(x), \Pi(y)):=(-1)^{p(\cB)+p(x)+p(x)p(y)}\cB(x,y)$ for every $x,y\in V$}.
\]

Thus, the matrix $B=\mmat R,S,T,U,$ of a \textit{symmetric} homogeneous form $\cB$ satisfies the condition
\be\label{BilSy}
\mmat R,S,T,U,=
\mmat R^{t},(-1)^{p(B)}T^{t},(-1)^{p(B)}S^{t},-U^{t}, .
\ee

Most popular normal forms (shapes) of the even non-degenerate symmetric
form are the ones whose supermatrices in the standard format are the
following normal forms:
\[
B_{ev}(m|2n)= \mmat 1_m,0,0,J_{2n},,\quad \text{where
$J_{2n}=\mmat 0,1_n,-1_n,0,$}.
\]

The usual notation for $\faut (B_{ev}(m|2n))$ is $\fosp(m|2n)$;
sometimes one write more explicitly, $\fosp^{sy}(m|2n)$. The non-degenerate anti-symmetric form on a~superspace of superdimension $2n|m$ is preserved by the
``symplectico-orthogonal\rq\rq\ Lie superalgebra~ ``$\fsp\fo(2n|m)$\rq\rq; the passage from $V$ to $\Pi(V)$ described above turns the anti-symmetric form on $V$ into a~symmetric one on $\Pi(V)$, so that ``$\fsp\fo(2n|m)$\rq\rq\ $\simeq \fosp^{sy}(m|2n)$ and its even part is $\fsp(2n)\oplus\fo(m)$. The argument of ``$\fsp\fo$\rq\rq, as that of $\fosp^{sy}$, is the superdimension of the superspace acted upon; this is what distinguishes rows 4 and 5 of table \eqref{table3}.

A~non-degenerate \textbf{symmetric} (resp. \textbf{anti-symmetric}) odd bilinear form $B_{odd}(n|n)$ can
be reduced to a~normal shape whose matrix in the standard format is
$J_{2n}$ (resp. $\Pi_{2n}$). The usual notation for
$\faut (B_{odd}(\Size))$ is~$\fpe(\Size)$.
The passage from $V$ to $\Pi (V)$ establishes an isomorphism
$\fpe^{sy}(\Size)\cong\fpe^{a}(\Size)$. These isomorphic Lie
superalgebras are called, as A.~Weil suggested, \textit{periplectic}.

The \textit{special periplectic} superalgebra~is\index{$\fspe(n)$, special periplectic superalgebra}
\[
\fspe(n)=\{X\in\fpe(n)\mid \str
X=0\}.
\]

We will also need the Lie superalgebras\index{$ \fspe(n)_{a, b}:=\fspe(n)\inplus \Kee(az+bd)$}
\begin{equation}
\label{spe}
 \fspe(n)_{a, b}:=\fspe(n)\inplus \Kee(az+bD), \text{
where $z=1_{2n}$, $D=\diag(1_{n}, -1_{n})$, and $a,b\in\Cee$}.
\end{equation}


\ssec{A.~Sergeev's central extension $\fa\fs$ of $\fspe(4)$}\label{ssas}
In 1970's A.~Sergeev proved that over $\Cee$ there is just one
nontrivial central extension of $\fspe(n)$ for~$n>2$. It exists only
for $n=4$ and we denote it by~$\fas$. Let us represent an arbitrary
element $A\in\fas$ as a~pair $A=X+a\cdot z$, where $X\in\fspe(4)$,
$a\in{\Cee}$ and $z$ is the central element. The bracket in $\fas$\index{$\fas$}
is
\begin{equation}
\label{2.1.4}
\left[\mat {A~&B\\ C&-A^t} +a\cdot z,\
\mat{A' & B' \cr C' & -A'{}^t} +a'\cdot
z\right]= \left[\mat{A& B \cr C & -A^t},\
\mat{ A' & B' \cr C' & -A'{}^t
}\right]+\tr(C\widetilde C')\cdot z,
\end{equation}
where $\ \widetilde {}\ $ is extended via~linearity from matrices
$C_{ij}=E_{ij}-E_{ji}$ on which $\widetilde C_{ij}=C_{kl}$ for every
even permutation $(1234)\longmapsto(ijkl)$.

The Lie superalgebra~$\fas$ can also be described with the help of
the spinor representation. For this, consider the Poisson Lie superalgebra $\fpo(0|6)$ whose
superspace is the Grassmann superalgebra~$\Lambda(\xi, \eta)$
generated by $\xi_1, \xi_2, \xi_3, \eta_1, \eta_2, \eta_3$ and the
multiplication is the Poisson bracket~ \eqref{pb}.

Recall that $\fh(0|6)=\Span (H_f\mid f\in\Lambda(\xi, \eta))$, see eq.~\eqref{2.3.2'}. Now,
observe that $\fspe(4)$ can be embedded into~$\fh(0|6)$. Indeed,
setting $\deg \xi_i=\deg \eta_i=1$ for all $i$ we introduce
a~$\Zee$-grading on $\Lambda(\xi, \eta)$ which, in turn, induces
a~$\Zee$-grading on $\fh(0|6)$ of the form
$\fh(0|6)=\mathop{\oplus}_{i\geq -1}\fh(0|6)_i$. Since
$\fsl(4)\cong\fo(6)$, we can identify $\fspe(4)_0$ with~$\fh(0|6)_0$.

It is not difficult to see that the elements of degree $-1$ in the
standard gradings of $\fspe(4)$ and $\fh(0|6)$ constitute isomorphic
$\fsl(4)\cong\fo(6)$-modules. It is subject to a~direct verification
that it is possible to embed $\fspe(4)_1$ into~$\fh(0|6)_1$.

Sergeev's extension $\fas$ is the result of the restriction of the cocycle that turns $\fh(0|6)$ into
$\fpo(0|6)$ to
$\fspe(4)\subset\fh(0|6)$. The quantization $\cQ$ deforms $\fpo(0|6)$ into
$\fgl(\Lambda(\xi))$; the through maps $T_\lambda:
\fas\tto\fpo(0|6)\stackrel{\cQ}{\tto}\fgl(\Lambda(\xi))$ are representations of
$\fas$ in the $4|4$-dimensional modules $\spin_\lambda$ isomorphic
to each other for all $\lambda\neq 0$. The explicit form of
$T_\lambda$ is as follows: 
\begin{equation}
\label{2.1.4.2} T_\lambda\colon \mat{A & B \cr C & -A^t }+a\cdot
z\longmapsto \mat{A & B-\lambda \widetilde C \cr C & -A^t
}+\lambda a\cdot 1_{4|4}, 
\end{equation}
where $1_{4|4}$ is the unit matrix and $\widetilde C$ is defined in
the line under eq.~\eqref{2.1.4}. Clearly, $T_\lambda$ is an
irreducible representation for every~$\lambda$.

\ssec{More notation} The tautological representation of a~linear Lie superalgebra~$\fg$
or its subalgebra~$\fg\subset \fgl(V)$ in $V$, and sometimes \textbf{the module $V$ itself} are denoted by $\id$ or, for clarity,\index{$\id$}~$\id_{\fg}$. The context prevents confusion of these notation with
that of the identity (scalar) operator $\id_{V}$ on the space $V$,
as in the next paragraph:

For $\fg=\mathop{\oplus}_{i\in\Zee}\fg_{i}$, the trivial
representation of $\fg_0$ is denoted by $\Cee$ (if $\fg_0$ is
simple) whereas $\Cee[k]$\index{$\Cee[k]$} denotes a~1-dimentional $\fg_0$-module
trivial on the semi-simple part of $\fg_0$ and such that $k$ is the value of the
central element $z\in\fg_{0}$ on $\Cee[k]$, where $z$ is chosen so that
$z|_{\fg_i}=i\cdot \id_{\fg_i}$.

If $\fg_0=\fa_1\oplus\dots\oplus\fa_k$ is a~direct sum of ideals and $V_i$ is
an $\fa_i$-module for each $i$, then\index{$\boxtimes$}
\begin{equation}\label{outer}
V_1\boxtimes\dots\boxtimes V_k
\end{equation}
denotes the \textit{outer} tensor product: the superspace
$V_1\otimes\dots\otimes V_k$ on which each summand $\fa_i$ acts through the
$i$th factor alone, i.e., 
\[
\text{$a(v_1\otimes\dots\otimes
v_k)=(-1)^{p(a)(p(v_1)+\dots+p(v_{i-1}))}v_1\otimes\dots\otimes av_i
\otimes\dots\otimes v_k$ for every $a\in\fa_i$.}
\] 
The sign $\otimes$ is reserved for
the tensor product of modules over one and the same Lie superalgebra, on which
the elements act by the Leibniz rule.

If $\fc$ is a~Lie algebra~of scalar matrices, and $\fg\subset \fgl (n|n)$ is a~ subsuperalgebra~containing $\fc$, then the \textit{projective}
Lie superalgebra~of type $\fg$ is $\fpg:= \fg/\fc$.\index{$\fpg:= \fg/\fc$} Examples: $\fpq (n)$, $\fpsq (n)$; $\fpgl (n|n)$, $\fpsl (n|n)$;
observe that $\fpgl (p|q)\cong \fsl (p|q)$ if~$p\neq q$.

\ssec{The serial vectorial Lie superalgebras}\label{ssSerialDef} Here
are the definitions of the four series of vectorial Lie superalgebras and of
their relatives.

\textbf{1) General algebras}. Let $x=(u_1, \dots
, u_m, \theta_1, \dots , \theta_n)$, where the $u_i$ are even
indeterminates and the $\theta_j$ are odd ones. Set $\fvect
(m|n):=\fder\; \Cee[x]$; it is called \textit{the general vectorial
Lie superalgebra}.

\textbf{2) Special algebras}. The
\textit{divergence} of the field
$D=\mathop{\sum}_if_i\pder{u_{i}} + \mathop{\sum}_j
g_j\pder{\theta_{j}}$ is the function (in our case: a~polynomial, or
a~series)\index{$\Div D$}
\begin{equation}
\label{2.2.4} 
\textstyle \Div D:=\mathop{\sum}\limits_i\pderf{f_{i}}{u_{i}}+
\mathop{\sum}\limits_j (-1)^{p(g_{j})}
\pderf{g_{i}}{\theta_{j}}.
\end{equation}

$\bullet$ The Lie superalgebra~$\fsvect (n|m):=\{D \in \fvect
(n|m)\mid \Div D=0\}$ is called the \textit{special} (or
\textit{divergence-free}) \textit{vectorial superalgebra}.
Equivalently,
\[
\fsvect(n|m)=\{ D\in \fvect (n|m)\mid L_D\vvol _x=0\},
\]
where $\vvol_x$ is the
volume form with constant coefficients in coordinates $x$ and $L_D$
the Lie derivative along the vector field~$D$.

$\bullet$ The Lie superalgebra~$\fsvect'(1|m):=[\fsvect(1|m),
\fsvect(1|m)]$ is called the \textit{traceless special vectorial
superalgebra}.

$\bullet$ The deformation of $\fsvect(0|m)$ is the Lie superalgebra
\[
\fsvect_{\lambda}(0|m):=\{D \in \fvect (0|m)\mid \Div
(1+\lambda\theta_1\cdot \dots \cdot \theta_m)D=0\},
\]
where $p(\lambda)\equiv m\pmod 2$; if $\lambda$ is even it is an element of $\Cee$, whereas if it is odd it is a new anti-commuting indeterminate. The result of this deformation is called the \textit{deformed
special} (or \textit{divergence-free}) \textit{vectorial
superalgebra}. Clearly, $\fsvect_{\lambda}(0|m)\cong
\fsvect_{\mu}(0|m)$ for $\lambda\mu\neq 0$. So we briefly denote
these deformations by $\widetilde{\fsvect}(0|m)$.

\textbf{3) The algebras preserving Pfaff equations and certain
differential 1- and 2-forms}.

$\bullet$ Let $u=(t, p_1, \ldots , p_n, q_1, \dots , q_n)$ be even indeterminates,; let
$\Theta=\begin{cases}(\xi, \eta)& \text{ if }\ m=2k\\(\xi, \eta,
\theta)&\text{ if }\ m=2k+1,\end{cases}$ be odd indeterminates. 
Set
\begin{equation}
\label{2.2.6} \alpha_1=dt+\mathop{\sum}\limits_{1\leq i\leq
n}(p_idq_i-q_idp_i)+ \mathop{\sum}\limits_{1\leq j\leq
k}(\xi_jd\eta_j+\eta_jd\xi_j)
\begin{cases}&
\text{ if }\ m=2k\\
+\theta d\theta&\text{ if }\ m=2k+1.\end{cases} 
\end{equation}

The \textit{contact} Lie superalgebra~is the one that preserves the \textit{Pfaff equation}
\[
\text{$\alpha_1(X)=0$ for $X\in \fvect(2n+1|m)$}
\]
or, equivalently, preserves the distribution singled out by the form $\alpha_1$, i.e., the superalgebra
\begin{equation}
\label{k} \fk (2n+1|m)=\{ D\in \fvect (2n+1|m)\mid
L_D\alpha_1=f_D\alpha_1\text{ for some }f_D\in \Cee [t, p, q,
\theta]\}.
\end{equation}
The Lie superalgebra
\begin{equation}
\label{po}
\begin{array}{c}
\fpo (2n|m)=\{ D\in \fk (2n+1|m)\mid L_D\alpha_1=0\}\end{array}
\end{equation}
is called the \textit{Poisson} superalgebra.\index{$ \fk (2n+1\vert m)$}\index{$ \fpo (2n\vert m)$} Its quotient modulo center is called the \textit{Hamiltonian Lie algebra} $\fh(2n|m)$ and is usually defined by its following property:
\be\label{hamOm}
\fh(2n|m):=\{D\in\fvect(2n|m)\mid L_D(\omega_0)=0,\ \ \text{where $\omega_0:=d\alpha_1$}\}.
\ee

$\bullet$ Similarly, set $
q=(q_1, \dots , q_n)$ be even,
let $\xi=(\xi_1, \dots , \xi_n)$ and $\tau$ be odd. Set
\begin{equation}
\label{alp0}
\begin{array}{c}
\alpha_0=d\tau+\mathop{\sum}\limits_i(\xi_idq_i+q_id\xi_i),
\qquad\qquad \omega_1=d\alpha_0
\end{array}
\end{equation}
and call these forms, as A.~Weil advised, the \textit{pericontact}
and \textit{periplectic}, respectively.\index{Form, pericontact}\index{Form, periplectic}

The \textit{pericontact} Lie superalgebra~is the one that preserves the Pfaff equation
\[
\text{$\alpha_0(X)=0$ for $X\in \fvect(n|n+1)$}
\]
or, equivalently, preserves the distribution singled out by the form $\alpha_0$:
\begin{equation}
\label{m} \fm (n):=\{ D\in \fvect (n|n+1)\mid L_D\alpha_0=f_D\cdot
\alpha_0\text{ for some }\; f_D\in \Cee [q, \xi, \tau]\}.
\end{equation}

The Lie superalgebra\index{$\fm (n)=\fm (n\vert n+1)$}\index{$\fb (n)=\fb (n\vert n+1)$, Buttin superalgebra}\index{Buttin bracket a.k.a. anti-bracket a.k.a. Schouten bracket}\index{Schouten bracket a.k.a. anti-bracket a.k.a. Buttin bracket}
\begin{equation}
\label{but}
\begin{array}{c}
\fb (n):=\{ D\in \fm (n)\mid L_D\alpha_0=0\}\end{array}
\end{equation}
is called the \textit{Buttin} superalgebra~in honor of C.~Buttin who
was the first to show that the \textit{Schouten bracket} of the functions that generate $\fb (n)$, see \eqref{2.3.7}, satisfies the super
Jacobi identity. Currently this bracket is often called \textit{anti-bracket}.\index{Anti-bracket a.k.a. Buttin bracket a.k.a. Schouten bracket}

The following are respective \textit{divergence-free} (or \textit{special}) subsuperalgebras\index{$\fsm (n)=\fsm (n\vert n+1)$}\index{$\fs\fb (n)=\fs\fb (n\vert n+1)$}
\begin{equation}
\label{sm1}
\begin{array}{c}
\fsm (n)=\{ D\in \fm (n)\mid \Div\ D=0\}, \ \ \
 \fs\fb (n)=\{ D\in
\fb (n)\mid \Div\ D=0\}.
\end{array}
\end{equation}

\sssec{Generating functions}\label{SS:2.4.1}{}~{}

$\bullet$ \underline{Odd form $\alpha_1$}. For every $f\in\Cee [t, p,
q, \Theta]$, set:\index{$K_f$}\index{$H_f$}
\begin{equation}
\label{2.3.1} 
\textstyle K_f:=(2-E)(f)\pder{t}-H_f + \pderf{f}{t} E,
\end{equation}
where $E:=\mathop{\sum}_i y_i \pder{y_{i}}$ (here the $y_{i}$
are all the coordinates except $t$) is the \textit{Euler operator},
and 
\begin{equation}
\label{2.3.2'} 
\textstyle H_f:=\mathop{\sum}\limits_{i\leq n}\left
(\pderf{f}{p_i} \pder{q_i}-\pderf{f}{q_i} \pder{p_i}\right)
-(-1)^{p(f)}\left(\mathop{\sum}\limits_{j\leq
k}\left(\pderf{f}{\xi_j} \pder{\eta_j}+ \pderf{f}{\eta_j}
\pder{\xi_j}\right)+ \pderf{f}{\theta} \pder{\theta}\right).
\end{equation}

 $\bullet$ \underline{Even form $\alpha_0$}. For every $f\in\Cee [q,
\xi, \tau]$, set:\index{$M_f$}\index{$\Le_f$}
\begin{equation}
\label{2.3.3} 
\textstyle M_f:=(2-E)(f)\pder{\tau}- Le_f -(-1)^{p(f)}
\pderf{f}{\tau} E, 
\end{equation}
where $E:=\mathop{\sum}_iy_i \pder{y_i}$ (here the $y_i$ are
all the coordinates except $\tau$), and
\begin{equation}
\label{2.3.4} 
\textstyle \Le_f:=\mathop{\sum}\limits_{i\leq n}\left(
\pderf{f}{q_i}\ \pder{\xi_i}+(-1)^{p(f)} \pderf{f}{\xi_i}\
\pder{q_i}\right).
\end{equation}
Since
\begin{equation}
\label{2.3.5}
\renewcommand{\arraystretch}{1.4}
\begin{array}{l}
 L_{K_f}(\alpha_1)=2 \pderf{f}{t}\alpha_1=K_1(f)\alpha_1, \\
L_{M_f}(\alpha_0)=-(-1)^{p(f)}2 \pderf{
f}{\tau}\alpha_0=-(-1)^{p(f)}M_1(f)\alpha_0,
\end{array}
\end{equation}
it follows that $K_f\in \fk (2n+1|m)$ and $M_f\in \fm (n)$.

$\bullet$ To the (super)commutators $[K_f, K_g]$ or $[M_f, M_g]$
there correspond \textit{contact brackets} of the generating functions:
\[
\renewcommand{\arraystretch}{1.4}
\begin{array}{l}
{}[K_f, K_g]=K_{\{f, \; g\}_{k.b.}};\qquad
 {}[M_f, M_g]=M_{\{f, \;
g\}_{m.b.}}.\end{array}
\]
The explicit formulas for the contact brackets are as follows. Let
us first define the brackets on functions that do not depend on $t$
(resp. $\tau$).

The \textit{Poisson bracket} $\{\cdot , \cdot\}_{P.b.}$ 
is given by the formula
\begin{equation}
\label{pb}
\renewcommand{\arraystretch}{1.4}
\begin{array}{ll}
 \{f, g\}_{P.b.}&\textstyle :=\mathop{\sum}\limits_{i\leq n}\ \bigg(\displaystyle\pderf{f}{p_i}\
\pderf{g}{q_i}-\ \pderf{f}{q_i}\
\pderf{g}{p_i}\bigg)-\\
&\textstyle (-1)^{p(f)}\bigg(\mathop{\sum}\limits_{j\leq m}\left(
\displaystyle\pderf{f}{\xi_j}\ \pderf{ g}{\eta_j}+\pderf{f}{\eta_j}\
\pderf{ g}{\xi_j}\right)+\displaystyle\pderf{f}{\theta}\ \pderf{
g}{\theta}\bigg)\text{ for every }f, g\in \Cee [p, q, \xi, \eta,
\theta].
\end{array}
\end{equation}

The \textit{Buttin bracket} or \textit{Schouten bracket} a.k.a. \textit{antibracket} $\{\cdot , \cdot\}_{B.b.}$ is given by the formula
\begin{equation}
\label{2.3.7} 
\textstyle \{ f, g\}_{B.b.}:=\mathop{\sum}\limits_{i\leq n}\
\bigg(\pderf{f}{q_i}\ \pderf{g}{\xi_i}+(-1)^{p(f)}\
\pderf{f}{\xi_i}\ \pderf{g}{q_i}\bigg)\text{ for every }f, g\in \Cee
[q, \xi].
\end{equation}

 In terms of the Poisson and Buttin brackets,
respectively, the contact brackets are
\begin{equation}
\label{2.3.8} 
\textstyle \{ f, g\}_{k.b.}=(2-E) (f)\pderf{g}{t}-\pderf{f}
{t}(2-E) (g)-\{ f, g\}_{P.b.}
\end{equation}
and
\begin{equation}
\label{2.3.9} 
\textstyle \{ f, g\}_{m.b.}=(2-E)
(f)\pderf{g}{\tau}+(-1)^{p(f)} \pderf{f}{\tau}(2-E) (g)-\{ f,
g\}_{B.b.}. 
\end{equation}

The Lie superalgebras of \textit{Hamiltonian
fields} (or \textit{Hamiltonian
superalgebra}) and its derived subalgebra~(defined only if $n=0$)
are (recall the definition of $\omega_0$, see \eqref{hamOm} and that the result of Berezin integration is the highest coefficient in the Taylor series expansion in coordinates)
\begin{equation}
\label{h}
\begin{array}{l}
\fh (2n|m)=\Span(D\in \fvect (2n|m)\mid L_D\omega_0=0),\\
\fh' (0|m)=\Span(H_f\in \fh (0|m)\mid \int f\vvol_{\theta}=0), ~~\text{where $\int$ is the Berezin integral}.
\end{array}
\end{equation}
The ``odd\rq\rq\  analogues of the Lie superalgebra~of Hamiltonian
fields are the Lie superalgebra~of vector fields $\Le_{f}$
introduced in \cite{L1} and its special subalgebra (recall the definition of $\omega_1$, see \eqref{alp0}):
\begin{equation}
\label{len}
\begin{array}{l}
\fle (n)=\Span(D\in \fvect (n|n)\mid L_D\omega_1=0),\\
\fsle (n)=\Span(D\in \fle (n)\mid \Div D=0).
\end{array}
\end{equation}

It is not difficult to prove the following isomorphisms (as superspaces):
\begin{equation}
\label{iso}
\begin{array}{rcl}
\fk (2n+1|m)&\cong&\Span(K_f\mid f\in \Cee[t, p, q, \xi]),\\
 \fh
(2n|m)&\cong&\Span(H_f\mid f\in \Cee [p, q, \xi]);\\
\fm (n)&\cong&\Span(M_f\mid f\in \Cee [\tau, q, \xi]),\\
\fle
(n)&\cong&\Span(\Le_f\mid f\in \Cee [q, \xi]).\\
\end{array}
\end{equation}
We see that the commutants (first derived algebras) in the purely odd case can be described as
\[
\textstyle \fpo' (0|m)=\Span(K_f\in \fpo (0|m)\mid \int f\vvol_\theta=0);\quad
\fh' (0|m)=\fpo' (m)/\Cee\cdot K_1.
\]

\ssec{Divergence-free subalgebras}\label{SS:2.4.3} Since, as is easy to
calculate,
\begin{equation}
\label{div}
 \Div K_f =(2n+2-m)K_1(f),
\end{equation}
it follows that the divergence-free subalgebra~of the contact Lie
superalgebra~either coincides with the whole algebra~(for $m=2n+2$)
or is the Poisson superalgebra. For the pericontact series, the
situation is more interesting: the divergence free subalgebra~is
simple and new, as compared with the algebras defined so far.
Since
\begin{equation}
\label{2.4.1} 
\textstyle \Div M_f =(-1)^{p(f)}2\left ((1-E)\pderf{f}{\tau} -
\mathop{\sum}\limits_{i\leq n}\nfrac{\partial^2 f}{\partial q_i
\partial\xi_i}\right ), 
\end{equation}
it follows that the divergence-free subalgebra~of the pericontact
superalgebra~is
\begin{equation}
\label{sm2} 
\textstyle \fsm (n) = \Span\left (M_f \in \fm (n)\mid
(1-E)\pderf{f}{\tau} =\mathop{\sum}\limits_{i\leq n}\nfrac{\partial^2
f}{\partial q_i
\partial\xi_i}\right ).
\end{equation}
In particular,
\begin{equation}
\label{2.4.2}
\textstyle \Div \Le_f = (-1)^{p(f)}2\mathop{\sum}\limits_{i\leq
n}\nfrac{\partial^2 f}{\partial q_i \partial\xi_i}.
\end{equation}
Set
\[
\textstyle \fsle (n):=\Span\left(\Le_f\in \fle(n)\mid\Delta(f)=0, \text{~~where~~} \Delta=\mathop{\sum}\limits_{i\leq n}\nfrac{\partial^2
}{\partial q_i
\partial\xi_i}\right).
\]
Lie superalgebras $\fsle (n)$, $\fs\fb (n)$ and $\fsvect (1|n)$
have simple ideals $\fsle '(n)$, $\fs\fb '(n)$ and $\fsvect
'(n)$ of codimension 1 defined from the exact sequences
\begin{equation}
\label{seqv}
\begin{array}{c}
0\tto \fsle '(n)\tto \fsle
(n)\tto \Cee\cdot Le_{\xi_1\dots\xi_n} \tto 0, \\
0\tto \fs\fb '(n)\tto \fs\fb
(n)\tto \Cee\cdot M_{\xi_1\dots\xi_n} \tto 0, \\
\displaystyle 0\tto \fsvect '(n)\tto
\fsvect (1|n)\tto \Cee
\cdot\xi_1\dots\xi_n\partial_{t}\tto 0.\end{array}
\end{equation}

\ssec{The Cartan prolongs}\label{2.6}
Recall that the graded Lie superalgebra
$\fb=\oplus_{k \geqslant -d}\ \fb_k$
is said to be \textit{transitive}\index{Lie (super)algebra, transitive} if for all~ $k\geq 0$ we have
\begin{equation*}
\label{12}
\{x\in \fb_k \mid [x,\fb_{-}]=0\}=0, \text{~~where $\fb_{-}:=\oplus _{k <0}\,\fb_k$}.
\end{equation*}

Let
$\fg_-=\mathop{\oplus}_{-d\leq i\leq -1}\ \fg_i$ be a~nilpotent
$\Zee$-graded Lie (super)algebra and let $\fg_0$ be a~ sub(super)algebra of the Lie
(super)algebra $\mathfrak{der}_0(\fg_-)$ of degree-0 derivations of~ $\fg_-$. 
The maximal transitive $\Zee$-graded Lie (super)algebra $\fg=\oplus_{k \geqslant -d}\ \fg_k$ whose non-positive part is the given $\fg_-\oplus \fg_0$ is called the 
\textit{Cartan prolong}\index{Cartan prolongation} --- the result of Cartan prolongation --- of the pair $(\fg_-,\fg_0)$ and is denoted by $(\fg_-,\fg_0)_*$. 

\textbf{A~partial Cartan prolong: prolongation involving
 a~positive part}\label{SS:2.6.3} Let $\fh_1\in \fg_1$ be
a~$\fg_0$-submodule such that $[\fg_{-1}, \fh_1]=\fg_0$.\index{Cartan prolongation, partial} If such
$\fh_1$ exists (usually, $[\fg_{-1}, \fh_1]\subset\fg_0$), define
the 2nd prolongation of $(\mathop{\oplus}_{i\leq 0}\fg_i,
\fh_1)$ to be
\begin{equation}
\label{partprol} \fh_{2}=\{D\in\fg_{2}\mid [D, \fg_{-1}]\in
\fh_1\}.
\end{equation}
The terms $\fh_{i}$, where $i>2$, are similarly defined. Set
$\fh_i=\fg_i$ for $i\leq 0$ and $\fh_*=\oplus\fh_i$.

\textit{Examples}: $\fvect(1|n; n)$ is a~subalgebra~of~$\fk(1|2n;
n)$. The former is obtained as the Cartan prolong of the same
nonpositive part as $\fk(1|2n; n)$ and a~submodule of $\fk(1|2n;
n)_1$. The simple exceptional superalgebra~$\fk\fas$ introduced in
\cite{Sh5, Sh14} is
another example.

\textbf{The Cartan prolong $\fu^\fg_{*} $ of $\fu_{-1}$ in $\fg$}. For $\fg=\mathop{\oplus}_{i \geq -1}\ \fg_i$, set\index{$\fu^\fg_i$}
\[
\begin{array}{l}
\fu^\fg_{-1} = \fu_{-1}\subset \fg_{-1}~\text{(non-strict inclusion)};\\
\fu^\fg_i := \{g\in \fg_i \mid [g, h] \subset \fu^\fg_{i-1} \ \text{
for every $h\in \fu_{-1}\}$ for $i \geq 0$}.\\
\end{array}
\]
The Lie superalgebra
\[
\fu^\fg_{*} = \mathop{\oplus}_{k\geq -1} \fu^\fg_{k}
\]
is said to be the \textit{Cartan prolong of $\fu_{-1}$ in $\fg$}.\index{Cartan prolong of $\fu_{-1}$ in $\fg$}

\ssec{Lie superalgebras of vector fields as Cartan
prolongs}\label{SS:2.6.2} Let $\id=V$ be the tautological $\fgl(V)$-module, considered also as a~$\fg$-module over any subsuperalgebra~$\fg\subset\fgl(V)$. Then, 
\begin{equation}
\label{carprol}
\begin{array}{ll}
\begin{array}{l}
\fvect(m|n)=(\id, \fgl(m|n))_*;\\
\fsvect(m|n)=(\id, \fsl(m|n))_*; \\
\fh(2m|n)=(\id, \fosp^{a}(m|2n))_*; \\
\end{array} &\begin{array}{l}
\fle(n)=(\id, \fpe^{a}(n))_*;\\

\fs\fle(n)=(\id, \fspe^{a}(n))_*.\\
\end{array}\end{array}
\end{equation}


$\bullet$ Define by ~$\fhei(2n|m)=W\oplus \Cee z$, where $\dim W=2n|m$ and $W$ is endowed with
a~non-degenerate antisymmetric bilinear form $B$, the Heisenberg Lie superalgebra. The bracket in it is given by the following relations:
\begin{equation}
\label{hei} \text{$z$ is even and spans the center; $[v,
w]=B(v, w)\cdot z$ for every $v, w\in W$.}
\end{equation}
Clearly, we have
\begin{equation}
\label{kont} \fk(2n+1|m)=(\fhei(2n|m), \fc(\fosp^{a}(m|2n)))_*\ .
\end{equation}

$\bullet$ The ``odd\rq\rq\  analog of $\fk$ is associated with the
following ``odd\rq\rq\  analog of~$\fhei(2n|m)$. Denote by $\fb\fa(n)=W\oplus \Cee\cdot z$ the
\textit{antibracket} Lie superalgebra~($\fb\fa$ is Anti-Bracket read
backwards), where $W$ is an
$n|n$-dimensional superspace endowed with a~non-degenerate
antisymmetric odd bilinear form $B$; the bracket in $\fb\fa(n)$ is
given by the following relations:
\begin{equation}
\label{2.5.4} \text{$z$ is odd and spans the center; $[v,
w]=B(v, w)\cdot z$ for every $v, w\in W$.}
\end{equation}

Clearly,
\begin{equation}
\label{mprol} \fm(n)=(\fb\fa(n), \fc(\fpe^{a}(n)))_*\ .
\end{equation}


\ssec{The exceptional simple vectorial Lie superalgebras}
The five families of such Lie superalgebras are
given below
as Cartan prolongs $(\fg_{-1},
\fg_{0})_{*}$ or generalized Cartan prolongs $(\fg_{-},
\fg_{0})_{*}$.

For \underline{depth $\leq 2$}, for
$\fg_{-}=\mathop{\oplus}_{-2\leq i\leq -1}\fg_{i}$, we write
$(\fg_{-2}, \fg_{-1}, \fg_{0})_{*}$ instead of $(\fg_{-},
\fg_{0})_{*}$. 

The $W$-graded superalgebras of \underline{depth $3$} appear as
regradings of the listed ones at certain values of $r$; the
corresponding terms $\fg_{i}$ for $i\leq 0$ will are given in
~\eqref{nonpos} and \eqref{depth3}; where  the trivial $\fg$-module is designated as~$\Cee_\fg$.

The non-positive components of exceptional simple vectorial Lie superalgebras over
$\Cee$, see \eqref{outer} for the sign $\boxtimes$:
\be\label{nonpos}\tiny
\renewcommand{\arraystretch}{1.3}
\begin{tabular}{|c|c|c|c|}
\hline $\fg$&$\fg_{-2}$&$\fg_{-1}$&$\fg_0$\cr
\hline \hline

$\fv\fle(3|6)$&$\id_{\fsl(3)}\boxtimes\Cee_{\fsl(2)}\boxtimes
\Cee[-2]$&$\id^*_{\fsl(3)}\boxtimes \id_{\fsl(2)}\boxtimes
\Cee[-1]$&$\fsl(3)\oplus\fsl(2)\oplus\Cee z$\cr
\hline

$\fv\fas(4|4)$&$-$&$\spin$&$\fas$\cr \hline

$\fk\fas(1|6)$&$\Cee[-2]\boxtimes\Cee_{\fo(6)}$&$\Cee\boxtimes\Pi(\id_{\fo(6)})$& $\fc(\fo(6))$\cr
\hline

$\fm\fb(3|8)$&$\id_{\fsl(3)}\boxtimes\Cee_{\fsl(2)}\boxtimes
\Cee[-2]$&$\Pi(\id^*_{\fsl(3)}\boxtimes \id_{\fsl(2)}\boxtimes
\Cee[-1])$&$\fsl(3)\oplus\fsl(2)\oplus\Cee z$\cr

\hline

$\fk\fle(5|10)$&$\id$&$\Pi(\Lambda^2(\id^*))$&$\fsl(5)$\cr
\hline
\end{tabular}
\ee
None of the above $\Zee$-graded vectorial Lie
superalgebras 
is of depth $>3$ and only one
is of depth~ 3:
\begin{equation}\label{depth3}
\fm\fb(3|8)_{-3}=\Pi(\Cee_{\fsl(3)}\boxtimes\id_{\fsl(2)}\boxtimes\Cee[-3]).
\end{equation}

\ssec{The modules of tensor fields}\label{SS:2.7} Let $\fg=\fvect(m|n)$ and $\fg_{\geq}=\mathop{\oplus}_{i\geq
0}\fg_{i}$ in the standard grading ($\deg x_i=1$ for all $i= 1, \dots m+n$). For every other $\Zee$-graded vectorial Lie superalgebra
for whose component $\fg_0$ the notion of lowest weight can be
defined, the construction is identical.

Let $V$ be the $\fg_0=\fgl(m|n)$-module with the \textit{lowest} weight
$\lambda=\lwt(V)$. Let us make $V$ into a~$\fg_{\geq}$-module by setting
$\fg_{+}\cdot V=0$ for $\fg_{+}=\mathop{\oplus}_{i>
0}\fg_{i}$. Let us realize $\fg$ by vector fields on the \textit{linear
supermanifold} $\cC^{m|n}$ with coordinates~$x=(u, \xi)$. Due
to the Poincar\'e--Birkhoff--Witt theorem, $T(V):=\Hom_{U(\fg_{\geq})}(U(\fg), V)$ is isomorphic to ${\Cee}[[x]]\otimes
V$. Its elements have a~natural interpretation as formal
\textit{tensor fields of type}~$V$. If $\lambda=(a, \dots , a)$ we
will simply write $T(\vec a)$ instead of~$T(\lambda)$. We will
usually consider $\fg$-modules coinduced from irreducible
$\fg_0$-modules. For example, $T(\vec 0)$ is the superspace of functions;
$\Vol(m|n)=T(1, \dots , 1\mid -1, \dots , -1)$ where the bar separates
the first $m$ (``even\rq\rq) coordinates of the weight with respect to
the matrix units $E_{ii}$ of $\fgl(m|n)$, is the superspace of
\textit{densities} or
\textit{volume forms}.
We denote the generator of $\Vol(m|n)$, as $\Cee[x]$-module,
corresponding to the ordered set of coordinates $x$ by~$\vvol_x$.
The space of $\lambda$-densities is $\Vol^{\lambda}(m|n):=T(\lambda,
\dots , \lambda; -\lambda, \dots , -\lambda)$. In particular,
$\Vol^{\lambda}(m|0)=T(\vec \lambda)$ but
$\Vol^{\lambda}(0|n)=T(\overrightarrow{-\lambda})$.




\ssec{Deformations of the Buttin superalgebra~(after
\cite{LS3})}\label{SS:2.8} As is clear from the definition of the
Buttin bracket, there is a~regrading (namely, $\fb (n; n)$ given by
$\deg\xi_i=0, \deg q_i=1$ for all $i$) under which $\fb(n)$,
initially of depth 2, takes the form
$\fg=\mathop{\oplus}_{i\geq -1}\fg_i$ with
$\fg_0=\fvect(0|n)$ and $\fg_{-1}\cong \Pi(\Cee[\xi])$. Let us replace now
the $\fvect(0|n)$-module $\fg_{-1}$ of functions with inverted
parity by the rank 1 (over the algebra~of functions) module of
$\lambda$-densities (also with inverted parity), i.e., set
$\fg_{-1}\cong \Pi(\Vol^\lambda(0|n))$, where the
$\fvect(0|n)$-action on the generator $\vvol_\xi^\lambda$ is given
by the formula
\begin{equation}
\label{2.7.1} L_D(\vvol_\xi^\lambda) =\lambda \Div D\cdot
\vvol_\xi^\lambda\; \text{ and }\;
p(\vvol_\xi^\lambda)=\od.
\end{equation}

Define $\fb_{\lambda}(n; n)$ to be the Cartan prolong
\begin{equation}
\label{fb} (\fg_{-1}, \fg_0)_*=(\Pi(\Vol^\lambda(0|n)),
\fvect(0|n))_*.
\end{equation}
Clearly, this is a~deform of~$\fb(n; n)$. The collection of these
$\fb_{\lambda}(n; n)$ for all $\lambda$'s is called the \textit{main
deformation}, the other deformations, defined in what follows, will
be called \textit{singular}.

The deform $\fb_{\lambda}(n)$ of $\fb(n)$ is a~regrading of
$\fb_{\lambda}(n; n)$ described as follows. For $\lambda
=\nfrac{2a}{n(a-b)}$, set
\begin{equation}
\label{2.7.2} 
\textstyle \fb_{a, b}(n) =\{M_f\in \fm (n)\mid a\; \Div
M_f=(-1)^{p(f)}2(a-bn)\pderf{f}{\tau}\}. 
\end{equation}
Taking into account the explicit form of the divergence of $M_{f}$
we get
\begin{equation}
\label{2.7.4}
\begin{array}{ll}
\fb_{a, b}(n) &=\{M_f\in \fm (n)\mid (bn-aE)\pderf{f}{\tau} =
a\Delta
f\}=\\
&\{D\in\fvect(n|n+1) \mid L_{D}(\vvol_{q, \xi,
\tau}^a\alpha_{0}^{a-bn})=0\}.\end{array}
\end{equation}
It is subject to a~direct verification that $\fb_{a, b}(n)\simeq
\fb_\lambda(n)$ for $\lambda =\nfrac{2a}{n(a-b)}\in\Cee P^1$. The Lie superalgebras $\fb_\lambda(n)$
are simple for $n>1$ and $\lambda\neq 0$, 1, $\infty$ for reasons
clear from eq.~\eqref{2.7.5}. It is also clear that the
$\fb_{\lambda}(n)$ are non-isomorphic for distinct $\lambda$'s, bar
occasional isomorphisms, see \cite{Ko1,LS1}.

The Lie superalgebra~$\fb(n)=\fb_{0}(n)$ is not simple: it has a~$(0|1)$-dimensional center. Recall the the superdimension $a|b$ of a superspace is, in terms of super $K_0$-functor, an element $a+b\eps$, where~$\eps^2=1$. This in mind, observe that
at $\lambda=1$ and $\infty$ the Lie superalgebras $\fb_{\lambda}(n)$
are not simple either: each of them has an ideal of codimension
$\eps^{n+1}$ and $\eps^{n}$, respectively. The corresponding exact
sequences are
\begin{equation}
\label{2.7.5}
\begin{array}{c}
0\tto \Cee M_{1} \tto \fb(n)\tto \fle(n)\tto 0,\\
0\tto \fb_{1} '(n)\tto \fb_{1}(n)\tto \Cee\cdot
M_{\xi_1\dots\xi_n} \tto 0,\\
0\tto \fb_{\infty} '(n)\tto \fb_{\infty}(n)\tto \Cee\cdot
M_{\tau\xi_1\dots\xi_n} \tto 0.\\
\end{array}
\end{equation}
Clearly, at the exceptional values of $\lambda$, i.e., 0, 1, and
$\infty$, the deformations of $\fb_{\lambda}(n)$ should be
investigated extra carefully. As we will see immediately, it is rewarded: at each of the exceptional points we find extra deformations. An
exceptional deformation at $\lambda=-1$ remains inexplicable. Other
exceptional values ($\lambda =\nfrac12$ and $-\nfrac32$) come from the
isomorphisms $\fh_{\lambda}(2|2)\cong \fh_{-1-\lambda}(2|2)$ and $\fb_{1/2}(2; 2)\cong \fh_{1/2}(2|2)=\fh(2|2)$.

Two members of the list \eqref{1.4} below are singular deformations, and are
defined here. First,
\begin{equation}\label{hlam}
\fh_{\lambda}(2|2):=\fb_{\lambda}(2; 2),
\end{equation}
the main deformation for $n=2$ read in the Hamiltonian realization; its bracket
is at the same time a~deform of the Poisson bracket and a~deform of the
antibracket, which is what singles this family out, see \cite{LS1}. Secondly,
$\widetilde{\fs\fb}_{\mu}$ is a~deformation of the non-simple $\fs\fb(n|n+1;
n)$ --- a~deform if the parameter of the deformation is even --- obtained as
the Cartan prolong
\begin{equation}\label{sbmu}
\begin{array}{l}
\widetilde{\fs\fb}_{\mu}(2^{n-1}-1|2^{n-1})=(\fh_{-1}, \fh_{0})_*\text{ for $n$
even},\\
\widetilde{\fs\fb}_{\mu}(2^{n-1}|2^{n-1}-1)=(\fh_{-1}, \fh_{0})_*\text{ for $n$
odd},
\end{array}
\end{equation}
where $\fh_{0}=\widetilde{\fsvect}(0|n)$ and
$\fh_{-1}=\nfrac{\Pi(\Vol(0|n))}{\Cee(b+\alpha\xi_1\cdots\xi_n)\vvol(\xi)}$, of
$\sdim=(2^{n-1}-1|2^{n-1})$ for $n$ even and $(2^{n-1}|2^{n-1}-1)$ for $n$ odd,
with $\mu=\nfrac{\alpha}{b}$, see \cite{LS1}. Here, $p(\mu)\equiv n\pmod 2$, so
for $n$ odd the parameter $\mu$ is odd, i.e., an extra indeterminate the
coefficients depend on; the algebras with $\mu\neq 0$ are isomorphic to one
another.

\ssec{The regradings}\label{ssRegr} A~$\Zee$-grading of a~vectorial
Lie superalgebra~is given by the vector $\vec r$ of the degrees of its
indeterminates. For a~W-grading, this vector can be shorthanded to the number
$r$ of the indeterminates of degree $0$, or to a~symbol; the standard
realization is the one with $r=0$, and then $r$ is not indicated. The
W-gradings of the three series from which all the other serial ones inherit
theirs are, see \cite{LS1},
\begin{equation}\label{regr}
\renewcommand{\arraystretch}{1.3}
\begin{array}{ll}
\fvect(m|n; r),\ 0\leq r\leq n:&\deg\xi_j=0\text{ for $1\leq j\leq r$};\ \deg
u_i=\deg\xi_{r+s}=1\text{ for all $i$, $s$};\\
\fm(n; r),\ 0\leq r<n-1:&\deg\tau=\deg q_i=2,\ \deg\xi_i=0\text{ for $1\leq
i\leq r$};\\
&\deg q_{r+j}=\deg\xi_{r+j}=1\text{ for all $j$};\\
\fm(n; n):&\deg\tau=\deg q_i=1,\ \deg\xi_i=0\text{ for $1\leq i\leq n$};\\
\fk(2n+1|m; r),\ 0\leq r\leq[\nfrac{m}{2}]:&\deg t=\deg\xi_i=2,\
\deg\eta_i=0\text{ for $1\leq i\leq r$};\\
&\deg p_i=\deg q_i=\deg\theta_j=1\text{ for all $i$, $j$}.
\end{array}
\end{equation}
The gradings of the series $\fvect$ induce those of $\fsvect$ and $\fsvect'$;
the gradings of $\fm$ induce those of $\fb_{\lambda}$, $\fle$, $\fsle$,
$\fsle'$, $\fb$, $\fs\fb$ and $\fs\fb'$; the gradings of $\fk$ induce those of
$\fpo$, $\fh$ and $\fh'$. The values $r=k-1$ for $\fk(1|2k)$ and $r=n-1$ for
$\fm(n)$ and its subalgebras are excluded because for them the $\fg_0$-module
$\fg_{-1}$ is reducible, so the grading is not a~W-grading; this is the reason
for the exclusions printed in table \eqref{1.4} below.

Besides these, the algebra~$\fb_{a, b}(2)$ has one exceptional W-regrading,
denoted $\Reg_{\fb}$,\index{$\Reg_{\fb}$}\index{$\Reg_{\fh}$} given by
\begin{equation}\label{RegB}
\deg\tau=0,\quad \deg\xi_1=\deg\xi_2=-1,\quad \deg q_1=\deg q_2=1;
\end{equation}
through the isomorphism \eqref{hlam} it determines the exceptional W-regrading
$\Reg_{\fh}$ of $\fh_{\lambda}(2|2)$. Both are listed in table \eqref{1.4} in
the slot of the parameter~$r$.

\ssec{The list of serial simple W-graded vectorial Lie superalgebras}\label{ssList}
In the table, label FD marks the particular cases of finite-dimension.\index{$\fvect(m\vert n)$}\index{$\fsvect(m\vert n)$}\index{$\widetilde{\fsvect}(0\vert n)$}\index{$\fk(2m+1\vert n)$}
\index{$\fh(2m\vert n)$}\index{$\fm(n\vert n+1)$}\index{$\fh_{\lambda}(2\vert 2)$}\index{$\fsm(n\vert n+1)$}\index{$\fb_{\lambda}(n\vert n+1)$}\index{$\fle(n\vert n+1)$}\index{$\fsle(n\vert n+1)$}\index{$\widetilde{\fs\fb}_{\mu}(2^{n-1}-1\vert 2^{n-1})$}
\begin{equation}\label{1.4}\tiny 
\begin{array}{ll}
\renewcommand{\arraystretch}{1.4}
\begin{tabular}{|l|l|}
\hline $N$&the algebra and conditions for its simplicity\cr 
\hline

1&$\fvect(m|n; r)$ for $m\geq 1$ and $0\leq r\leq n$ \cr 
\hline
\hline
FD &
$\fvect(0|n; r)$ for $n>1$ and $0\leq r\leq n$ \cr 
\hline

2&$\fsvect(m|n; r)$ for $m>1$, $0\leq r\leq n$ \cr 
\hline

FD &
$\fsvect(0|n; r)$ for $n>2$ and $0\leq r\leq n$ \cr 
\hline

3&$\fsvect'(1|n; r)$ for $n>1$, $0\leq r\leq n$\cr 
\hline

FD &$\widetilde{\fsvect}(0|n)$ for $n>2$ \cr \hline

\hline 
4&$\fk(2m+1|n; r)$ for $0\leq r\leq [\nfrac{n}{2}]$,
$(m|n)\neq (0|2k)$\cr 
\hline
&$\fk(1|2k; r)$ for $0\leq r\leq k$ except
$r=k-1$\cr
 \hline
 \hline

5&$\fh(2m|n; r)$ for $m>0$ and $0\leq r\leq [\nfrac{n}{2}]$\cr
\hline

FD&$\fh'(0|n)$ for $n>3$ \cr 
\hline
\hline
\end{tabular}&\begin{tabular}{|l|l|}
\hline $N$&the algebra and conditions for its simplicity\cr 
\hline

6&$\fh_{\lambda}(2|2; r)$ for $\lambda\neq -2, -\nfrac32, -1, \nfrac12,
0, 1, \infty$, and \cr &$r=0$, $1$ and $\Reg_{\fh}$ \cr 
\hline


\hline
7&$\fm(n|n+1; r)$ for $0\leq r\leq n$, $r\neq n-1$\cr
 \hline
 
8&$\fsm(n|n+1; r)$ for $n>1$ but $n\neq 3$ and \cr
&$0\leq r\leq n$, $r\neq n-1$\cr \hline

9&$\fb_{\lambda}(n|n+1; r)$ for $n>1$, where $\lambda\neq 0, 1,
\infty$\cr
& and $0\leq r\leq n$, $r\neq n-1$;\cr
&\cr
&$\fb_{\lambda}(2|3;
r)$, for $\lambda\neq -2, -\nfrac32, -1, \nfrac12, 0, 1, \infty$, and \cr 
&$r=0$, $1$ and $\Reg_{\fb}$ \cr \hline

10&$\fb_{1}'(n|n+1; r)$ for $n>1$, $0\leq r\leq n$, $r\neq
n-1$\cr \hline

11&$\fb_{\infty}'(n|n+1; r)$ for $n>1$, $0\leq r\leq n$, $r\neq n-1$\cr \hline 

12&$\fle(n|n; r)$ for $n>1$, $0\leq r\leq n$, $r\neq n-1$\cr
\hline
\hline

13&$\fsle'(n|n; r)$ for $n>2$, $0\leq r\leq n$, $r\neq n-1$\cr
\hline

14&$\widetilde{\fs\fb}_{\mu}(2^{n-1}-1|2^{n-1})$ for $\mu\neq 0$ and
$n>2$\cr \hline
\end{tabular}
\end{array}
 \end{equation}
 Hereafter we abbreviate algebras $\fm(n|n+1;
r)$, $\fb_{\lambda}(n|n+1; r)$, $\fle(n|n; r)$ and so on to $\fm(n; r)$, $\fb_{\lambda}(n; r)$, $\fle(n; r)$,
and so on, respectively.

\section{Tables with answers. Theorems}\label{STabAns}
\be\label{table1}
\begin{array}{c}
\text{Case $1^o$ of \eqref{5cases}, i.e., $\fu_{-1}=\fg_{-1}$ and $\fu_{0}=\fg_{0}$; $\fu_{1}\neq 0$. Here, $N$ is as in table \eqref{1.4b}}\\
\footnotesize \renewcommand{\arraystretch}{1.4}
\begin{array}{|c|c|c|} 
\hline
N&\fg&\fu\\
\hline\hline
1&\fvect(m|n)&\fd(\fsvect(m|n))\text{ or } \fpgl(m+1|n) \\
\hline
13&\fh'(0|6)&\fspe(4)\\
\hline
14&\fle(n)&\fd(\fsle(n))\\
\hline
\end{array}
\end{array}
\ee

Observe that $\fpgl(m+1|n)\simeq \fsl(m+1|n)$ for~$m+1\neq n$.

In certain cases, it is possible to reduce our problem by means of regradings
of the ambient algebras to 
the standard grading. In these cases, it suffices to consider maximal subalgebras
with respect to the standard grading and then apply a~regrading of both ambient and the embedded subalgebra.

\begin{equation}\label{table2}
\begin{array}{l}
\text{Case $3^o$ of \eqref{5cases}: $\fu_{-1}=\fg_{-1}$ and $\fu_{0}\subset \fg_{0}$ is a
maximal reducible subalgebra.}\\
\text{For $\fg$ of series $\fvect$ and $\fsvect$, let $W= \fg_{-1}/V$ for every $V\subsetneq \fg_{-1}$;}\\
\text{for $\fg$ of series $\fh$, $\fle$ and $\fsle$, let either}\\ 
\phantom{.}\hskip 3mm \text{(i) $\omega|_{V}$ be non-degenerate on a~ non-zero $V \subset \fg_{-1}$ or}\\ 
\phantom{.}\hskip 3mm \text{(ii) $\omega|_{V}$ vanish on a non-zero $V \subset \fg_{-1}$.}\\
\footnotesize
\renewcommand{\arraystretch}{1.4}
\begin{tabular}{|c|c|}
\hline
$\fg(\dim \fg_{-1})$&$\fu$\cr
\hline\hline
$\fvect(m|n)$&$\faut(V)=\cF(W)\otimes \fvect(V)\inplus \fvect(W)$\cr
\hline
$\fsvect(m|n)$&$\faut(V)\cap \fg = \cF(W) \otimes \fsvect(V)
\inplus \varphi_{s}(\fvect(W))$\cr
\hline
$\begin{array}{l}
\fh(2m|n)\\ 
\text{or $\fh'(0|n)$}
\end{array}$&$\begin{array}{ll}
(i)&\fh(V) \oplus \fh(V^\perp)\text{
for a~subsuperspace $V \subset \fg_{-1}$ or}\cr
(ii)&\faut (V, \varpi; V^\perp/V)\text{ for every
isotropic $V$ except for such that}\\
&\text{$\dim V_{\od}= k-1$, where $2k=n$;}\\
&\text{in the excluded case, $\fu = \cF(V) \otimes
\fpo(V^\perp/V)/\Cee\cdot 1\inplus \fvect(V)$}\\
\end{array}$\cr 
\hline
$\fle(n|n)$&$\begin{array}{ll}
(i)&\fle(V) \oplus \fle(V^\perp)\text{
for a~subsuperspace $V \subset
\fg_{-1}$ or}\cr 
(ii)&\faut
(V, \varpi; V^\perp/V)\text{ for every
isotropic $V$ except for such that}\\
&\text{$V= V_{\ev}$ and $\dim V=n-1$}\\
\end{array}$\cr \hline
$\fs\fle(n|n)$&$\begin{array}{ll}
(i)&(\fs\fle(V) \oplus \fs\fle(V^\perp))\inplus \Cee\cdot
D \text{ for a~subsuperspace $V \subset \fg_{-1}$ or}\cr 
(ii)&\cF(V) \otimes \fs\fb(V^\perp/V)/\Cee\cdot
1\inplus \varphi_{s}(\fvect(V))\cr 
&\text{ for every isotropic $V$
except for such that $V= V_{\ev}$ and $\dim V=n-1$}\cr 
\end{array}$\cr
\hline
\end{tabular}
\end{array}
\end{equation}

\begin{equation}\label{table3}
\begin{array}{c}
\text{Case $3^o$ of \eqref{5cases}}\\
\footnotesize 
\protect
\renewcommand{\arraystretch}{1.4}
\setlength{\tabcolsep}{4pt}
\begin{tabular}{|c|c|c|}
\hline $\fg(\dim \fg_{-1})$&$\fu$&$\fu_0$\cr
\hline\hline
$\fvect(AB)=\fvect(V\otimes W^*)$&$\fpgl(A+B)=\fpgl(V\oplus W)$&
$\fp(\fgl(V)\oplus \fgl(W))$\cr
for $A:=\sdim V$, $B:=\sdim W$&&\cr \hline

$\fvect(E^2(V))$&$\fo(2n)=\fo(V\oplus V^*)$&$\fgl(V)$\cr
\hline
$\fvect(S^2(V))$&$\fsp(2n)=\fsp(V\oplus
V^*)$&$\fgl(V)$\cr \hline

$\fvect(E^2(V))$, see eq.~\eqref{dims}&$\fosp^{sy}(2m|2n)=\fosp(V\oplus
V^*)$&$\fgl(V)$\cr
\hline

$\fvect(S^2(V))$, see eq.~\eqref{dims}&``$\fsp\fo(2m|2n)$\rq\rq$\simeq\fosp^{sy}(2n|2m)$&$\fgl(V)$\cr
&$=\faut(B_{\text{a}}(V\oplus V^*))$&\cr
\hline
\hline
$\fvect(n|n)$&$\fpq(n+1)$&$\fq(n)\oplus\Cee J$ for $\fq(n)=C(J)$\cr
\hline
\hline

$\fvect(m|2n)$&$\fosp(m+2|2n)=\fosp(W\oplus V\oplus
V^*)$,&$\fc(\fosp(W))$\cr 
$=\fvect(W)$&where $\dim V=1$&\cr

\hline \hline

$\fvect(n|n)$&$\fpe(n+1)=\fpe(W\oplus V\oplus V^*)$,&$\fc(\fpe(W))$\cr
$=\fvect(W)$&where $\dim V=1$&\cr \hline

$\fvect(n^2|n^2)$ for $n>1$&$\fpe(2n)=\fpe(V\oplus V^*)$&$\fgl(V)$ for $\sdim V=n|n$\cr \hline

\hline $\fvect(k(1+\eps)2^{n-k})$&$\fvect(0|n)=\fvect(\Lambda(n))$&
$\fvect(0|n-k)\niplus \fgl(k; \Lambda(n-k))$\cr
$\simeq\fvect(V\otimes \Lambda(W))$&$\simeq\fvect(V\oplus W)$&
$\simeq\fvect(W)\niplus \fgl(V; \Lambda(W))$\cr \hline

$\fsvect(k(1+\eps)2^{n-k})$&$\fsvect(0|n)=\fsvect(\Lambda(n))$&
$\fvect(0|n-k)\niplus \fsl(k; \Lambda(n-k))$\cr
$\simeq\fsvect(V\otimes \Lambda(W))$ if $k\neq 1$&$\simeq\fsvect(V\oplus W)$&
$\simeq\fvect(W)\niplus \fsl(V; \Lambda(W))$\cr 
\hline
$\fsvect((1+\eps)2^{n-1})$&$\fsvect(0|n)=\fsvect(\Lambda(n))$&
$\varphi_{s}(\fvect(0|n-1))$\cr
\text{\textup{(}the case $k=1$ of the previous row\textup{)}}&&\cr \hline
$\begin{cases}\fh(\Lambda(n-2))&\text{for $n$ even}\\
\fle(\Lambda(n-2))&\text{for $n$
odd}\end{cases}$&$\fh(0|n)$&$\fh(0|n-2)\niplus \Lambda(n-2)\cdot z$\cr
\hline \hline $\fvect(16)$&$\fe_6$&$\fc(\fo(10))$\cr \hline
$\fvect(27)$&$\fe_7$&$\fc(\fe_6)$\cr
\hline \hline
$\fvect(2|2)$&$\fosp_a(4|2)$&$\fc(\fosp(2|2))\simeq\fgl(1|2)$\cr
\hline $\fvect(6|4)$&$\fab(3)$&$\fc(\fosp(2|4))$\cr \hline
\end{tabular}
\end{array}
\end{equation}
\normalsize

Note that
\be\label{dims}
\begin{array}{ll}
\sdim S^2(V )=\nfrac12\left(m(m+1)+n(n-1)\right)|mn&\text{if $\sdim V=m|n$}\\
\sdim E^2(V )=\nfrac12\left(m(m-1)+n(n+1)\right)|mn&\text{if $\sdim V=m|n$}\\
\end{array}
\ee
\vskip 0.5 cm

\begin{equation}\label{table4}
\footnotesize
%
\noindent
\renewcommand{\arraystretch}{1.4}
\begin{array}{|c|c|}
\hline
\fg&\fu=\fu^\fg_*\\
\hline\hline
\fvect(m|n)& \text{Ann}(\fu_{-1})\fvect(m|n)\inplus 
\fvect(\fu_{-1})\\
\fsvect(m|n)& \text{Ann}(\fu_{-1})\fsvect(m|n)\inplus 
\varphi_{s}(\fvect(\fu_{-1})), \ \ \text{see formula~\eqref{varphi_s}}\\
\fh(2m|n)& \text{Ann}(V)\fpo(2m|n)\inplus 
(\cF(V/(V \cap V^\perp)) \otimes \fvect(\fg_{-1}/(V + V^\perp)
\inplus \fh (V/(V \cap V^\perp)))\\
\hline
\fle(n|n)&
\text{Ann}(V)\fb(n) \inplus (\cF(V/(V \cap V^\perp)) \otimes
\fvect(\fg_{-1}/(V + V^\perp)\inplus 
\fle(V/(V \cap V^\perp)))\\
\hline
\fs\fle(n|n)& (\text{Ann}(V)\fb(n)\cap \fsle(n))\inplus (\cF(V/(V
\cap V^\perp)) \otimes \fsvect(\fg_{-1}/(V + V^\perp) \inplus 
\varphi_{s}(\fle(V/(V \cap V^\perp))))\\
&\text{for every nonsingular subsuperspace $V \subset\fg_{-1}$, see formula~\eqref{varphi_s}}\\
\hline
\fs\fle'(n|n)& \left(\text{the previous line}\right)\cap \fsle'(n|n)\\
\hline
\end{array}
\end{equation}
 

\vskip 0.5 cm
\normalsize

\noindent
{\protect
\begin{equation}\label{table5}
\footnotesize
\renewcommand{\arraystretch}{1.4}
\begin{array}{|c|c|}
\hline
 \fg&\fu\\
\hline\hline
 \fvect(V_{1}^*\otimes V_{2})& \fpgl(V_{1}\oplus V_{2}) \text{ if either $\dim V_1 \neq k + k\eps$ or $\dim V_2 \neq l + l\eps$ for some $k, l$}\\
\hline
 \fvect(m+n|m+n)&\fvect(m|n+1; 1)\text{ if either $m\neq n$ or $m+n\eps\neq 1, \eps$}\\
\hline\hline
 \fsvect(V_{1}^*\otimes V_{2}) &\fp\fgl(V_{1}\oplus V_{2}) \text{ if $\dim V_1 \neq k + k\eps$ or $\dim V_2 \neq l + l\eps$ for some $k, l$}\\
\hline
 \fsvect(V_{1}^*\otimes C(\cA)\otimes V_{2}) &\fp\fq(\cA\otimes (V_{1}\oplus V_{2})),\ \ \text{see Subsection~\ref{cA}}\\
\hline
 \fsvect(m+n|m+n) & (\fvect(m|n) \otimes \Lambda(\xi))\inplus \Cee\partial_{\xi}\text{ if $m\neq n$ and $m+n\eps\neq 1$, $\eps$}\\
\hline
 \fsvect(V\otimes \Lambda(W) )&\fvect(m|n + k; k) \otimes \Lambda(\xi) \text{ if $m+ n\eps \neq 1$, $\eps$, $1 + \eps$}\\
 & \text{and either $k > 1$ and $m \neq n$ or $k =1$ and $m = n$}\\
\hline
\fsvect(\Lambda(n))&\fk(1|2n; n)\\
\hline
\fsvect(\Pi(\Lambda(n)))&\fm(n; n)\\
\hline
\fsvect(2|2)&\fb_{\infty}(2; 2)\inplus \fsl(2)\\
\hline
\end{array}
\end{equation}
 }

\normalsize \vskip 0.5 cm

\ssbegin[Almost simple subsuperalgebras of (peri)contact superalgebras]
{Theorem}[Almost simple subsuperalgebras of (peri)contact superalgebras]\label{Th2.1.2} \textup{1)} The only almost simple Lie superalgebras admitting $\fk$-contact gradings are
listed in table \eqref{tabl}:
\footnotesize
\begin{equation}\label{tabl}
\renewcommand{\arraystretch}{1.4}
\begin{tabular}{|l|}
\hline
$\fosp(m|2n+2)\subset \fk(2n+1|m)$\ \ \textup{(a ``standard\rq\rq\  embedding)}\\
$\fosp_a(4|2)\subset \fk(1|4)$\\

\hline
$\fosp(4|2n)\subset \fk(1|4n)$\ \textup{in addition to 
the standard
embedding $\fosp(4|2n)\subset \fk(2n-1|4)$}\\

\hline
$\fas\subset \fk(1|6)$\\

\hline

$\fpgl(m|n+2)\simeq\fpgl(n+2|m) \subset \fk(2n+1|2m)$ \\

\hline
\textup{any algebra~$\fg$ such that $\fpsl(2|2)\subset\fg\subset\fder(\fpsl(2|2))$ can be embedded into $\fk(1|4)$}\\

\hline
$\fag(2)\subset \fk(1|7)$ \\
\hline
$\fab(3)\subset \fk(1|8)$ \\
\hline
$\fsq(n+1)\subset \fk(2n+1|2n+2; 1)$ \textup{(a~nonstandard grading, see the proof)}\\

\hline $\fosp(3|2n)\subset \fvect(1|2n)$\\ \hline\end{tabular} \end{equation}
\normalsize

Let $\alpha=dt-\sum_{i=1,2}(\xi_id\eta_i+\eta_id\xi_i)$ and let $\fk(1|4)$ preserve the distribution given by the form~
$\alpha$;
let $\fk(1|4)$ be considered with the standard $\Zee$-grading.
Explicitly, a~$\Zee$-grading-preserving embedding
\[
\fosp_a(4|2)\subset \fk(1|4)
\]
is given by the element
\[
\textstyle F=\nfrac12(1+a)t^2-(a-1)\Xi\in\fg_2,
\text{~~where $\Xi=\xi_1\xi_2\eta_1\eta_2$.}
\]
The element $F$ itself is the only, up to a~factor, degree-$2$
element of $\fosp_a(4|2)$; the degree-$1$ elements of $\fosp_a(4|2)$
are spanned by
\[
(1+a)t\xi_i+(a-1)\partial_{\eta_i}\Xi\; \text{ and }\;
(1+a)t\eta_i+(a-1)\partial_{\xi_i}\Xi, \text{ where $i=1$, $2$.}
\]

\textup{2)} The only almost simple Lie superalgebras admitting \textbf{standard}
$\fm$-contact gradings are
\[
\renewcommand{\arraystretch}{1.4}
\begin{array}{l}
\fpgl(m+1|n+1)\simeq\fpgl(n+1|m+1) \subset \fm(n+m);\\
 \fspe_{a, b}(n+1)\subset \fm(n).
\end{array}
\]
In item $2$\textup{)}, as in item $1$\textup{)}, the word ``only\rq\rq\ concerns the
\textbf{standard} $\fm$-contact gradings; for the nonstandard ones, cf.\ the line
$\fsq(n+1)\subset\fk(2n+1|2n+2;1)$ of table \eqref{tabl} and the paragraph on nonstandard $\fk$-
and $\fm$-contact gradings preceding case B\textup{)} in $\S\ref{Sd2contact}$.

\end{Theorem}


\section{Depth 1}\label{Sd=1}

\ssec{The list of simple vectorial Lie superalgebras in their W-gradings of depth~ 1}\label{2.1} Here is the excerpt from the classification \cite{LS1} of simple W-graded vectorial Lie superalgebras with polynomial coefficients; indicated are the conditions for simplicity (conditions for finite dimensionality are parenthesized); the maximal subalgebras of $\fvect(1|2)\simeq\fosp(2|2)\simeq\fsl(1|2)$, $\fsvect(0|3)\simeq\fspe(3)$, and $\fh'(0|4)\simeq\fpsl(2|2)$ are described in \cite{ShM}:
\begin{equation}\label{1.4b}\footnotesize
\begin{array}{ll}
\renewcommand{\arraystretch}{1.4}
\begin{tabular}{|l|l|}
\hline $N$&algebra\cr \hline

1&$\fvect(m|n;r)$ for $m\geq 1$ ($m=0$,\ $n>1$) \cr \hline

2&$\fsvect(m|n;r)$ for $m>1$ ($m=0$,\ $n>2$) \cr \hline

3&$\fsvect'(1|n;r)$ for $n>1$ \cr \hline

4 &$\widetilde{\fsvect}(0|n)$ (for $n>2$) \cr \hline

\hline 
5&$\widetilde{\fs\fb}_{\mu}(2^{n-1}-1|2^{n-1})$ for $\mu\neq 0$,
$n>2$\cr
\hline
6&$\fk(1|2n;n)$ \cr
\hline
7&$\fvle(4|3)$\cr
\hline
8&$\fvas(4|4)$\cr
\hline
9&$\fkas(;3\xi)$\cr
\hline
10&$\fkas(;3\eta)$\cr
\hline
\end{tabular}
&
\renewcommand{\arraystretch}{1.4}
\begin{tabular}{|l|l|}
\hline $N$&algebra\cr 
\hline
11&$\fh(2m|2n;n)$ for $m>0$ \cr
\hline
12&$\fh_{\lambda}(2|2;r)$ for $r=0$, $1$ \cr 
&and $\lambda\neq -2, -\nfrac32, -1, \nfrac12,
0, 1, \infty$ 
\cr\hline

13&$\fh'(0|n)$ (for $n>3$) \cr 
\hline

14&$\fle(n;r)$ for $n>1$ and $r=0$ or $n$ \cr
\hline

15&$\fsle'(n;r)$ for $n>2$ and $r=0$ or $n$\cr
\hline
\hline

16&$\fm(n;n)$\cr
\hline

17&$\fb_\lambda(n;n)$\text{~~ for $\lambda\neq 1, \infty$}\cr
\hline

18&$\fb_\infty'(n;n)$\cr
\hline

19&$\fb_1'(n;n)$,\cr \hline
\end{tabular}\\
\end{array}
\end{equation}


Observe that each new grading we consider gives rise to a~\textit{bigrading} with new degrees of certain indeterminates. Then, all operators diagonally acting in the homogeneous spaces relative one grading must remain in the 0th component of the other grading. Let the non-positive part of $\fu$ coincide with the non-positive part of $\fg$; then $\fu$ must contain all the diagonally acting elements of the standard grading. 

In particular, since series $\fvect$, $\fk$, and $\fm$, and the exceptions $\fv\fle$, $\fk\fas$ and $\fm\fb$ contain the grading operator of a~\lq\lq standard"\ grading, this operator belongs to any maximal subalgebra~$\fu$ (in any of the cases $1^o$ -- $3^o$ of the cases \eqref{5cases}), and therefore this subalgebra~is graded by the standard grading. Therefore, for these algebras $\fg$, it suffices to consider only the standard grading. Therefore, the cases $N=6$, 9, 10, 16 of table~\eqref{1.4b} should (and will) be dealt with when we consider depth~2. 

Let us consider the remaining cases $N=2-5$, $11-13$, 15 separately.

Thus, for the purpose of the case $1^o$--$3^o$ analysis, table \eqref{1.4b} boils down to the following table, in which precisely the cases $N=6$, 9, 10, 16 postponed to depth~2 and the cases $N=1$, 7, 8, 14, 17--19 to be considered now are retained, whereas the cases $N=2$--5, $11$--13, 15, to be treated separately, are omitted.
\begin{equation}\label{1.41}\footnotesize
\begin{array}{ll}
\renewcommand{\arraystretch}{1.4}
\begin{tabular}{|l|l|}
\hline $N$&algebra\cr \hline

1&$\fvect(m|n)$ for $m\geq 1$ ($m=0$,\ $n>1$) \cr \hline



\hline
6&$\fk(1|2n;n)$ \cr

\hline 
7&$\fvle(4|3)$\cr
\hline
8&$\fvas(4|4)$\cr
\hline
9&$\fkas(;3\xi)$\cr
\hline
10&$\fkas(;3\eta)$\cr
\hline
\end{tabular}
&
\renewcommand{\arraystretch}{1.4}
\begin{tabular}{|l|l|}
\hline $N$&algebra\cr 
\hline
\hline

14&$\fle(n; r)$ for $n>1$ and $r=0$ or $n$ \cr
\hline

16&$\fm(n;n)$\cr
\hline

17&$\fb_\lambda(n;n)$\cr
\hline

18&$\fb_\infty'(n;n)$\cr
\hline

19&$\fb_1'(n;n)$,\cr \hline
\end{tabular}\\
\end{array}
\end{equation}


\ssbegin[The case where $\fu\supset \fg_{-1}\oplus
\fg_{0}$]{Lemma}[The case where $\fu\supset \fg_{-1}\oplus
\fg_{0}$]\label{L3.1.1} 
In this case, i.e., for a~Lie subsuperalgebra~$\fu\subset\fg$ such that
$\fu\supset \fg_{-1}\oplus \fg_{0}$, 

 \textup{1)} For simple Lie superalgebras
$\fg=\fsvect(m|n)$ or $\fh(2m|n)$ or $\fh'(0|n)$ for $n \neq 6$, or $\fsle'(n)$ for $n>2$,
the subalgebra $\fu$ is maximal if and only if $\fu= \fg_{-1}\oplus\fg_{0}$.

\textup{2)} If $\fg= \fvect(m|n)$, then $\fu$ is maximal if and only if
either $\fu= \fd(\fsvect(m|n))$ or $\fu = \fpgl(m +1|n)$.

\textup{3)} If $\fg =\fle(n)$, then $\fu$ is maximal if and only if $\fu =
\fd(\fsle(n))=\fsle(n)\inplus \Cee E$, where $E$ is the Euler operator.

\textup{4)} If $\fg =\fh'(0|6)$, then $\fu$ is maximal if and only if $\fu =
\fspe(4)$.

In the cases \textup{2)--4)}, i.e., for $\fu\supsetneq \fg_{-1}\oplus
\fg_{0}$, see also table $\eqref{table1}$.
\end{Lemma}

If $r\ne 0$ in Table \eqref{1.4}, we arrive at the following Case.

\smallskip

\underline{$3^o$. $\fu_{-1}=\fg_{-1}$ and $\fu_{0}\neq \fg_{0}$}.

Indeed, for $\fg=\fvect(m|n;r)$, we have $\fg_0=\fgl(m|n-r)\otimes\Lambda(r)\inplus \fvect(r)$; this $\fg_0$ contains $\fh=\fsl(m|n-r)\otimes\Lambda(r)\inplus \fvect(r)$.

Here, $\fg_0$ differs from $\fh$ by an ideal $1_{m|n-r}\otimes \Lambda(r)$ for $m\ne n-r$, and by $d\otimes \Lambda(r)$, where $d\in\fgl(m|n-r)$, but $d\notin\fsl(m|n-r)$ if~$m=n-r$. 

For $\fu=\fd(\fsvect(m|n;r))$, the 0th component differs from $\fh$ only by one outer derivation, $\fd$, which is the Euler operator $E=\sum x_i\partial_{x_i}$. 

The same argument holds for regradings with $r\neq 0$ in the pair $\fle(n)\supset \fd(\fsle(n))$. 

One more example: 
\[
\fu=\fsl(m+1|n-r)\otimes \Lambda(r)\inplus \fvect(0|r)\subset \fvect(m|n;r). 
\]
Then, $\fu=\fu_{-1}\oplus \fu_0\oplus\fu_1$. Let $x$ stand for indeterminates of the $(m|n-r)$-dimensional part of $\fg_{-1}$, and $\xi$ denote the remaining $r$ odd indeterminates. Then, 
\[
\begin{array}{l}
\fu_{-1}=\Span(\partial_{x_i}\otimes \Lambda(\xi)\mid \text{for all $i$}), \\
\fu_0= \Span(x_i\partial_{x_j}\otimes \Lambda(\xi)\inplus \fvect(\xi)\mid \text{for all $i, j$}),\\
\fu_1=\Span(x_iE_x \otimes \Lambda(\xi)\mid \text{for all $i$}), \text{~~where $E_x =\sum x_i\partial_{x_i}$. }
\end{array}
\]

Let us regrade this algebra~$\fu$ to the standard grading, i.e., set $\deg \xi_j=1$ for all $j$; we designate the result by~$U$. Then, 
\[
\begin{array}{l}
U_{-1}=\Span(\partial_{x_i}, \partial_{\xi_j}\mid \text{for all $i, j$}),\\
U_0=\Span(x_i\partial_{x_j}, \xi_k\partial_{x_j}, \xi_k\partial_{\xi_s}\mid \text{for all $i,j, k,s$}).
\end{array}
\]
We see that the $U_0$-action on $U_{-1}=\fvect_{-1}$ is reducible. There is just one maximal subalgebra~with such an action: $\fvect(x)\otimes \Lambda(\xi)\inplus \fvect(\xi)$. If one considers the nonstandard regrading of this subalgebra~corresponding to $r\neq 0$, one gets an algebra~containing~$\fu$.

Thus, let
$\fg$ be a~simple $\Zee$-graded vectorial Lie superalgebra~of depth 1, i.e., 
\be\label{mainList}
\begin{array}{l}
\text{either serial of series $\fvect, \fsvect,\fh, \fle,\fsle$, $\fb_\lambda$,}\\ 
\text{or exceptional: $\fvle(4|3)$, $\fvas(4|4)$, $\fkas(;3\xi)$, $\fkas(;3\eta)$.}
\end{array}
\ee
Let
$\fu$ be its maximal subalgebra~such that
\be\label{Prop}
\begin{array}{ll}
\bullet&\fu_{-1}=\fg_{-1};\\
\bullet&\text{$\fu_0$ is a~maximal proper subalgebra~in $\fg_0$;}\\
\bullet&\text{$\fu_0$ is almost simple;}\\
\bullet&\text{the $\fu_0$-action in $\fu_{-1}=\fg_{-1}$ is irreducible.}\\
\end{array}
\ee

\ssec{Ideas of the proof}\label{ssIdeas} 
\be\label{u0s}
\begin{minipage}[l]{13cm}
Let $\fu_{0s}$ denote the simple part of $\fu_0$, if $\fu_0$ is either centerless or the center of $\fu_0$ is a~direct summand, whereas if $\fu_0$ contains a~non-trivial central extension of the simple Lie superalgebra, then $\fu_{0s}$ will stand for this non-trivial central extension.
\end{minipage}
\ee

{\bf Idea~1.} \textsl{Since $[\fu_{-1}, \fu_1]$ is an ideal in $\fu_0$, this ideal must contain $\fu_{0s}$.}

Let us find out if the restriction of the $\fu_0$-action in $\fu_{-1}$ onto $\fu_{0s}$ remains irreducible.

What differs $\fu_{0s}$ from~$\fu_{0}$?

First, the presence of a~center. But if the center is even, then it acts by a~scalar in any irreducible module and can not affect irreducibility.

What if the center is odd? Let $J$ denote the operator via which this odd central element acts on~$\fu_{-1}$. First of all, clearly, $\Ker J=0$, since $\Ker J$ is invariant with respect to the whole $\fu_0$, and the $\fu_0$-action in $\fu_{-1}$ is irreducible. Besides, $J^2$ is an even operator commuting with the $\fu_0$-action in~$\fu_{-1}$. Therefore, $J^2$ is a~non-zero scalar operator and we can normalize it to be equal to~$-\id$. But then $\fu_0\subset \fq(\fu_{-1})\oplus \Cee J$. The Cartan prolong $(\fu_{-1}, \fu_0)_*$ is described in in \cite[Theorem ~6.2.2]{BLLS}: this prolong is the image of $\fq(n+1)$
under the embedding 
\be\label{embQ}
\fq(n+1)\subset \fsvect(n|n). 
\ee
The maximality of this subalgebra~is evident.

Second, $\fu_{0}$ may contain outer derivations of~$\fu_{0s}$. 

{\bf Idea~2.} \textsl{Let $\rho$ be an irreducible representation of $\fu_{0s}$ in a~proper subspace $V$ of~$\fu_{-1}$. The stabilizer of $\rho$ in $\fu_0$ must coincide with $\fu_0$, since otherwise the representation of $\fu_0$ in $\fu_{-1}$ is induced, but the maximal subalgebra~corresponding to the induced representation is not almost simple, see \cite{ShM}.} 

This means that $\rho$ can be extended to a~representation of the whole~$\fu_0$. If the space of outer derivations added to $\fu_{0s}$ is of dimension $\leq 2$, then $V=\fu_{-1}$, because if the Lie superalgebra~is of superdimension $1$, or $2$, or $1|1$, then all its irreducible finite-dimensional representations are 1-dimensional. Recall that $\dim\fder(\fg)\leq 2$ for every simple Lie superalgebra $\fg$ and the only simple Lie superalgebra $\fg$ with $\dim\fder(\fg)=2$ is~$\fpsl(2|2)$.

Thus, if $\fu_0\neq \fder(\fpsl(2|2))$, then the $\fu_{0s}$-action in $\fu_{-1}$ is irreducible. Therefore, in addition to \eqref{Prop}, we assume that 
\be\label{addProp}
\fu_0\ne \fder(\fpsl(2|2)).
\ee
The exceptional case $\fu_0= \fder(\fpsl(2|2))$ will be considered separately.

{\bf Idea~3.} \textsl{The subalgebra~$\fu_{0s}$ always contains a~reductive Lie algebra, which will be denoted~$\fu_{0r}$. If $\fu_{0s}$ is vectorial of type $\fvect(0|m), \fsvect(0|m)$, or $\fh(0|m)$, then $\fu_{0r}$ is the 0th component of the standard $\Zee$-grading of~$\fu_{0s}$. If $\fu_{0s}\simeq \widetilde{\fsvect}(0|m)$, then $\fu_{0r}$ is again $\fsl(m)$, the same as for~$\fsvect(0|m)$.
If $\fu_{0s}$ is not vectorial, then $\fu_{0r}$ is just the even component of~ $\fu_{0s}$.}

Consider the action of this $\fu_{0r}$ in~$\fu_k$. Since $\fu_{0r}$ is reductive, this action is completely reducible. 

\ssbegin[On simple maximal subalgebras]{Lemma}[On simple maximal subalgebras]\label{Lv} Let $\tilde V_k$ be the sum of all non-trivial submodules of the $\fu_{0r}$-module $\fu_k$, and $V_k$ be the minimal $\fu_{0s}$-submodule of $\fu_k$ containing~$\tilde V_k$.
Let $\fv\subset \fu$ be the subsuperalgebra~generated by
$$
V=\fu_{-1}\oplus \fu_{0s}\oplus(\oplus_{k=1}^\infty V_k).
$$
Then, the subalgebra~$\fv$ is simple. \end{Lemma} 

\begin{proof} First, let us prove that~$V_1\ne 0$.
Note that $\fu_1\ne 0$ due to the maximality of $\fu$: otherwise $\fu\subsetneq \fg_{-1}\oplus \fg_0$.

Suppose that~$V_1=0$. In particularly, this means that $\tilde V_1=0$, i.e., the reductive part $\fu_{0r}$ of $\fu_{0s}$ acts by 0 in the whole of~$\fu_1$. But each root (with respect to the maximal torus~of $\fu_{0r}$) vector of weight $\lambda\ne 0$ sends the vector of weight 0 of the $\fu_{0s}$-action in $\fu_1$ into a~weight vector of weight~$\lambda$. This means that the whole $\fu_{0s}$ acts by 0 on the whole of~$\fu_1\ne 0$. 

This implies, in particular, that $\fu_{0s}$ is centerless since the non-zero central element of $fu_0$ must act by a~non-zero scalar $\lambda$ in~$\fu_{-1}$. By Jacobi identity then, it must act in $\fu_1$ as multiplication by a~scalar $-\lambda\ne 0$.

Let $X\in \fu_1$ and~$X\ne 0$. Since $[\fu_{0s}, X]=0$, then by Jacobi identity the subspace $[\fu_{-1}, X]$ is an ideal in~$\fu_{0s}$. This ideal is non-zero because the ambient Lie superalgebra~$\fg$ is transitive. Besides, the map
\[
\ad_X:\fu_{-1}\tto \fu_{0s},\ x\mapsto [x,X]\text{~~for every $x\in\fu_{-1}$}
\]
is an intertwining operator between the $\fu_{0s}$-action in $\fu_{-1}$ and the adjoint action. Now, the Jacobi identity for the triple $x_1, x_2, X$, where $x_1, x_2\in \fu_{-1}$, leads to a~contradiction.

Thus,~$V_1\ne 0$.

Now, let $\fm$ be a~graded ideal in the subalgebra~$\fv$. Since the ambient Lie superalgebra~ $\fg$ is transitive, we can obtain a~non-zero element in $\fu_{-1}$ acting on a~non-zero element of $\fm$ sufficiently many times by appropriate elements of the component $\fv_{-1}=\fu_{-1}\subset \fv$, and hence $\fu_{-1}\subset \fm$ thanks to irreducibility of the action of the subalgebra~$\fu_{0s}\subset \fv$ in~$\fu_{-1}$.

Clearly, $[\fu_{-1},\fv_1]$ is a~non-zero ideal in~$\fv_0$. But this subalgebra~is contained in $\fu_0$, and hence is almost simple. Therefore, any non-zero ideal of $\fv_0$ either contains $\fu_{0s}$, or consists  only of the center~$z$. But the first component of the Cartan prolong of $(\fu_{-1},\Cee z)$ is equal to. zero. Hence, ${\fm_0\supset [\fu_{-1},\fv_1]\supset \fu_{0s}}$.

Considering now the brackets of $\fu_{0r}\subset \fm_0$ with subspaces $\tilde V_k\subset \fv_k$, we see that $\fm_k\supset \fv_k$, and hence $\fm\supset V$, which implies that~$\fm=\fv$.
\end{proof}

\ssec{Steps of the classification of maximal subalgebras}\label{reducing} Thus, the task of classifying maximal subalgebras with the properties \eqref{Prop} and \eqref{addProp} reduces to the following actions.
\begin{equation}\label{actions}
\begin{minipage}[l]{15.7cm}
\begin{itemize}
\item[(i)] For every simple Lie superalgebra~$\fv$, describe its $\Zee$-gradings of depth 1 such that
\begin{itemize}
\item $\fv_0$ is an almost simple Lie superalgebra;
\item The $\fv_0$-action in $\fv_{-1}$ is irreducible.
\end{itemize}

\item[(ii)] For every such $\Zee$-grading, establish into which of the vectorial Lie superalgebras $\fg$ listed in the ``main list" \ \eqref{mainList}
this $\Zee$-grading admits an embedding of~$\fv$.

\item[(iii)] Verify if $\fv_0$ is a~maximal subalgebra~in $\fg_0$. If it is not, find out if it is possible to enlarge it to a~maximal, and at the same time almost simple, subalgebra~ $\fu_0\subset \fg_0$.

\item[(iv)] Compute the Cartan prolong $(\fv_{-1}, \fu_0)_*$ which must be a~maximal subalgebra~in~ $\fg$.

\item[(v)] Consider the case $\fu_0 = \fder(\fpsl(2|2))$.

\end{itemize}
\end{minipage}
\ee

\sssec{The three exceptional ambients}\label{ss2fdEx}
Of the three exceptional finite-dimensional ambients $\fosp_a(4|2)$, $\fag(2)$ and $\fab(3)$, only
$\fosp_a(4|2)$ and $\fab(3)$ have a~$\Zee$-grading of depth $1$: for $\fosp_a(4|2)$ of $\sdim (9|8)$ it is given
by $\fg_0\simeq\fc(\fosp(2|2))$ of $\sdim (5|4)$ acting on $\fg_{-1}\simeq
\id_{\fc(\fosp(2|2))}$ of $\sdim (2|2)$, for $\fab(3)$ of $\sdim (24|16)$ by
$\fg_0=\fc(\fosp(2|4))$ acting on $\fg_{-1}$ of $\sdim (6|4)$, which is Kac's Jordan superalgebra
$K_{10}$ in the Tits--Kantor--Koecher picture, whereas $\fag(2)$ has no such grading. These
ambients, and the answers for them in every $\Zee$-grading to which the method of this paper
applies, are the subject of \cite{LOdd}; we do not consider them here.

\section{Proofs for the case of depth-1 ambient}\label{Sdeg1amb}

Let $\fu$ be a~Lie
subsuperalgebra~of $\fg$ such that $\fu\supset \fg_{-1}\oplus
\fg_{0}$.

\begin{proof}[Proof of Lemma $\ref{L3.1.1}$.\nopoint] Observe that in all the above cases the Lie superalgebra
$\fg$ is transitive and the component $\fg_1$ generates the whole
positive part of $\fg$, i.e., $\fg_+ = \mathop{\oplus}\limits_{i>0}\fg_{i}$.

1) This item follows, therefore, from irreducibility of $\fg_1$.

2) If $\fg= \fvect(m|n)$, then the $\fg_0$-module $\fg_1$ is a~direct
sum of two irreducible modules: $\fsvect(m|n)_1$ in the standard grading and
$\fg_{-1}^*\simeq \fpgl(m+1|n)_{1}$, in the grading of $\fpgl(m+1|n)$ such
that $\fpgl(m+1|n)=V\oplus \fpgl(V)\oplus V^*$, where $V=\fg_{-1}$.

Therefore, for a~nontrivial maximal
subalgebra~$\fu\neq \fg$, its component $\fu_1$ must coincide with
one of these submodules. If $\fu_1=\fsvect(m|n)_1$, then the condition 
\[
\Div[\partial_i,X]=\partial_i(\Div X)\text{~~for any $X\in\fvect(m|n)$ and any $\partial_i\in\fg_{-1}$}
\] 
implies that $\fu_2\subset \fsvect_2$, and, by induction, $\fu_k\subset \fsvect_k$ for any $k>0$, and hence $\fu= \fd(\fsvect(m|n))$. Whereas if $\fu_1=\fg_{-1}^*$, then the direct calculations show that $\fu_2=0$\, and hence $\fu = \fpgl(m +1|n)$. 

3) For $\fg = \fle(n)$, the $\fg_0$-module $\fg_1$ is reducible but not
completely
reducible: it contains just one irreducible submodule, $\fsle(n)_1$, and
the quotient modulo this submodule is isomorphic to $\fg_{-1}^*$, and it is
irreducible. Further arguments are similar to those in item 2) of this proof.

4) If $\fg= \fh'(0|6)$, then $\fg_0\simeq \fo(6) \simeq \fsl(4)$. The
20-dimensional $\fg_0$-module $\fg_1 \simeq \Lambda^3(\fg_{-1})$ splits into
the sum of two 10-dimensional ones, $\fg_1 \simeq \fg_{11}\oplus
\fg_{12}$.

The sums $\fu^{(1)} = \fg_{-1} \oplus \fg_{0}\oplus \fg_{11}$ and
$\fu^{(2)} = \fg_{-1} \oplus \fg_{0}\oplus \fg_{12}$ are isomorphic
to $\fspe(4)$. Then, in both cases $\fu_2=0$ and $\fu=\fspe(4)$.

The cases of nonstandard gradings are similar. \end{proof}

\underline{A) $\fu_{0}$ is a~maximal reducible subsuperalgebra~in
$\fg_{0}$}. Then,
$\fu$ and $\fg$ are as follows:

\textbf{$3^oA1)$ No non-degenerate bilinear form on $\fg_{-1}$.} Let
$V \subset \fg_{-1}$ be an $\fu_{0}$-invariant subsuperspace and let\index{$\fst(V)$}
\[
\fst(V) := \{X \in\fgl(\fg_{-1}) \mid X(V) \subset V\}
\]
be the stabilizer of $V$. Since $\fu_{0}$ is maximal in $\fg_{0}$, we
should have $\fu_{0} = \fst(V)\cap \fg_{0}$. Therefore, $\fu =
(\fg_{-1},\fst(V))_* \cap\fg$.

Let $x = (x_1, \ldots , x_k)$ and $y = (y_1, \ldots , y_l)$ be coordinates (of any parity but homogeneous) on
the supermanifold associated with the linear superspace $\fg_{-1}^*$
chosen so that
\begin{equation}\label{1.2}
V = \Span\left(\partial_{x_1}, \ldots , \partial_{x_k}\right),
\qquad \fg_{-1} = V \oplus \Span\left( \partial_{y_1}, \ldots
,\partial_{y_l}\right).
\end{equation}

Then,
\[
\fst(V) = \Span\left(x_i\partial_{x_j}, \;\;y_k\partial_{x_j},
\;\;y_k\partial_{y_s} \text{ for every $i,j, k,s$}\right)
\]
and
\be\label{autVcoord}
(\fg_{-1},\fst(V))_*\simeq \left(\cF(y)\otimes\fvect(x)\right) \inplus \fvect(y),
\ee
where $\cF(y)$ is the supercommutative superalgebra~of polynomials in
$y$, and $\fvect(y):=\fder\, \cF(y)$. Set\index{$\faut(V)$}
\[
\faut(V) :=(\fg_{-1},\; \fst(V))_*.
\]
Clearly, $\faut(V)$ is the Lie superalgebra~of vector fields preserving the
vector bundle structure on $\fg_{-1}$ with fiber $V$ and quotient $W$:
\[
0\tto V\tto \fg_{-1} \tto \fg_{-1}/V = W\tto 0.
\]
Thus,
\[
\faut(V)=(\fg_{-1},\fst(V))_*\simeq
\cF(W)\otimes\fvect(V)\inplus \fvect(W).
\]

Let $X_0 = \mathop{\sum}_{1\leq i\leq k}a_ix_i\partial_{x_i}
\in\fgl(\fg_{-1})$, where $a_i\in\Cee$, be such that
\[
c_0 = \Div
X_0 = \mathop{\sum}\limits_{1\leq i\leq k}(-1)^{p(x_i)}a_i\neq 0.
\]
In
particular, if $\dim V_{\ev} \neq \dim V_{\od}$, we can take for
$X_0$ the Euler operator: $X_0 = \mathop{\sum}_{1\leq i\leq
k}x_i\partial_{x_i}$. Then, the \textit{divergence-freeing}\index{$\varphi_{s}:\fvect(y) \tto
\fsvect(x, y)$, divergence-freeing} map $\varphi_{s}:\fvect(y) \tto
\fsvect(x, y)$ given by the formula
\be\label{varphi_s}
\varphi_{s}(D) = D- (\Div D)\nfrac{1}{c_0}X_0 \text{ for every } D \in\fvect(y)
\ee
determines an embedding of the Lie superalgebra~$\fvect(y)$ into
$\fsvect(x, y)$ such that
\[
\faut(V) \cap \fsvect(\fg_{-1}) = \cF(y)\otimes \fsvect(x)\inplus
\varphi_{s}(\fvect(y)).
\]
That $\varphi_s$ is a~homomorphism, and that $\Div\varphi_s(D)=0$, use only the facts that $D$ and
$X_0$ involve \textbf{disjoint} sets of indeterminates, so that $[D, X_0]=0$ and $X_0(\Div D)=0$,
and that $\Div X_0=c_0$ is a~constant. Therefore, the same formula defines, for the
\textbf{opposite} splitting, the map
\be\label{varphi_sbis}
\varphi_{s}:\fvect(x)\tto\fsvect(x, y), \quad \varphi_{s}(D) = D-(\Div D)\nfrac{1}{c_0}X_0, \quad\text{where~} X_0 = \mathop{\sum}\limits_{1\leq j\leq l}a_jy_j\partial_{y_j},
\ee
and $c_0 = \Div X_0=\sum_j(-1)^{p(y_j)}a_j\neq 0$; e.g., $X_0=y_1\partial_{y_1}$ and
$c_0=(-1)^{p(y_1)}$ will do. Note that $X_0\in\Ann(x)\fvect(x, y)$ in this case. It is
\eqref{varphi_sbis}, not \eqref{varphi_s}, that occurs in Theorem \ref{Th3.9} and in the
$\fsvect$ line of table \eqref{table4}: there $\fu_{-1}$ is a~\textbf{proper} subspace of
$\fg_{-1}$ and the transitive part of $\fu^\fg_*$ is $\fvect(x)$, not~$\fvect(y)$.

\textbf{$3^oA2)$ There exists a~non-degenerate $\fg_{0}$-invariant antisymmetric
bilinear form on~ $\fg_{-1}$}. Let us describe now $\fu= \faut(V)\cap
\fg$ for $\fg= \fh(2m|n)$, $\fh'(0|n)$, $\fle(n)$ and~$\fsle'(n)$.
These Lie superalgebras $\fg$ preserve an~antisymmetric bilinear
non-degenerate differential 2-form $\omega$ on the supermanifold
associated with the linear superspace $\fg_{-1}^*$, and $\omega$
induces a~$\fg_{0}$-invariant bilinear form on $\fg_{-1}$ denoted
for brevity also by~$\omega$. The form $\omega$ is even for the
$\fh$ and $\fh'$ series; it is odd for $\fle$ and~$\fsle'$. Thanks to a
result in \cite{ZZO}, recalled in \S~\ref{ZZO}, we will only be interested in the two
particular cases:
\[
\begin{array}{ll}
1)&\text{$\omega|_{V}$ is non-degenerate;}\\

2)&\text{$\omega|_{V}$ is identically zero.}
\end{array}
\]

For every superspace $V$, and a~subsuperspace $W$ with a~non-degenerate
symmetric bilinear form $\varpi$ on it, set (here $\cF(V)$ is the
superalgebra~of polynomial functions on $V$ and $\fvect(V):=\fder\, \cF(V)$)\index{$\faut(W, \varpi; V)$}
\be\label{omega0}
 \faut(W, \varpi; V)
:=\begin{cases}\cF(V) \otimes \fpo(W)/\Cee\cdot 1 \inplus \fvect(V)
&\text{if }p(\varpi) =\ev,\\ \cF(V) \otimes \fb(W)/\Cee\cdot 1 \inplus 
\fvect(V) &\text{if } p(\varpi) =\od.\end{cases}
\ee
We see that $\faut(W,\varpi; V)$ is the vectorial subsuperalgebra
preserving the vector bundle
\[
0\tto W\tto V \oplus W\tto V\tto 0
\]
and the form $\varpi$ in the fiber.

\underline{Let $p(\varpi) =\ev$}. The lift of this bundle to the bundle over
$T^*V$ considered with the normal (ortho)symplectic 2-form
provides us with the embedding
\[
\faut(W, \varpi; V) \tto \fh(T^*V \oplus W).
\]

\underline{Let $p(\varpi) =\od$}. The lift of this bundle to the bundle over
$\Pi(T^*V)$ considered with the normal periplectic 2-form
provides us with the embedding
\[
\faut(W, \varpi; V) \tto \fle(\Pi(T^*V) \oplus W).
\]
Define the orthogonal complement to $V$ with respect to $\omega$:
\[
V^\perp := \{X\in \fg_{-1}\mid \omega(X, v) = 0 \text{ for every }v \in
V\}.
\]



\underline{B) $\fu_{0}$ is a~maximal irreducible subsuperalgebra
in $\fg_0$}, either almost simple or a~central extension of an
almost simple Lie superalgebra. Then, $\fu$ and $\fg$ come from the
list of classical Hermitian spaces \cite{He,OV, GOV, GOV1,On} in the nonsuper
case and its super version due to Serganova~\cite{Se}; for elucidations of several Serganova's tables, see
\cite{Zi}. 

For Lie algebras, we look at the Dynkin graph of $\fu$,
each node being labeled by the corresponding coefficient in the
decomposition of the maximal root with respect to simple roots. The
label is equal to the depth of the $\Zee$-grading of~$\fu$. In this
subsection, we are interested in depth~1.

The answers are in Table \eqref{table3} (where $A$ and $B$ are the superdimensions of
the superspaces $V$ and $W$, respectively; they may degenerate into
dimensions, either purely even or purely odd, but not into 1 or $\eps$: in
these cases~$\fu_0=\fg_0$).

\begin{Remark}[No embeddings $\fab(3)\subset
\fosp(8|4)$] The existence of the embeddings 
\[
\fab(3)\subset \fvect(6|4)\;\text{ and $\fosp(8|4)\subset
\fvect(6|4)$}
\]
prompts a~question ``is there an embedding $\fab(3)\subset
\fosp(8|4)$?" The answer is negative since the irreducible
representation of $\fab(3)$ of the least dimension is the adjoint
one and there is a~non-degenerate invariant symmetric bilinear even form on~$\fab(3)$.
\end{Remark}

The maximal subalgebras of type B) and of depth 1 correspond to most studied
type of symmetric spaces and their super analogs. We list them at the end of the
article for convenience of the reader.

\underline{C) $\fu_{0}$ is a~maximal irreducible subsuperalgebra
in $\fg_0$,} neither almost simple nor a~central extension of an
almost simple Lie superalgebra. Then, $\fu$ and $\fg$ are as in
table \eqref{table5}.

\vskip 0.5 cm

\underline{$4^o$. $0 \subset \fu_{-1} \subset \fg_{-1}$, where both
inclusions are strict}. Then, $\fu=\fu^\fg_{*}$ is described in Theorem~
\ref{Th1.4}.

\ssbegin[$\fu^\fg_{*}$ is a~maximal Lie
subsuperalgebra~in $\fg$]{Lemma}[$\fu^\fg_{*}$ is a~maximal Lie
subsuperalgebra~in $\fg$]\label{L1.1} \textup{1)} $\fu^\fg_{*}$ is a~Lie
subsuperalgebra~in~$\fg$.

\textup{2)} 
$\fu^\fg_{*}$ is maximal among all subsuperalgebras of $\fg$ with given
$\fu_{-1}$, and 
$\fu^\fg_{*}=\fu^{\fvect(\fg^*_{-1})}\cap \fg$.
\end{Lemma}

\begin{proof} 1) Let $h_i\in \fu^\fg_{i}$ and $h_j\in \fu^\fg_{j}$. Let
us show that $[h_i, h_j] \in \fu^\fg_{i+j}$. In fact, for every $h\in
\fu_{-1}$, we have
\[
[h, [h_i, h_j]] = [[h, h_i], h_j] + (-1)^{p(h)p(h_i)} [h_i, [h, h_j]].
\]
By definition $[h, h_i]\in \fu^\fg_{i-1}$ and $[h, h_j]\in
\fu^\fg_{j-1}$. Now, the double induction on $i$ and $j$ gets us the
statement desired.

2) This is obvious. \end{proof}

\ssec{
A~description of $\fu^{\fvect}_*$} \index{$\fu^{\fvect}_*$}
Introduce indeterminates in $\fvect(m|n)$ so that
\[
\fu_{-1} = \Span\left(\partial_{x_1}, \ldots ,
\partial_{x_k}\right);\quad \fg_{-1} =\fu_{-1}\oplus
\Span\left(\partial_{y_1}, \ldots , \partial_{y_l}\right)\text{ with
$k + l = m + n$}
\]
and the parities of the $x$ and $y$ are such that there are $m$ even
and $n$ odd ones among them. Let\index{$\Ann(x) := \Ann(\fu_{-1})$}
\[
\Ann(x) := \Ann(\fu_{-1})\;\text{ be the superspace of
linear functions that vanish on $\fu_{-1}$; }
\]
in our notation, these functions are the
linear combinations of the $y$'s. Let
\[
I(x) := \Ann(x)\cF(x, y)
\]
 be the ideal generated by $\Ann(x)$
in the superalgebra~of functions.
Then, 
$\Ann(x)\fvect(x, y)$, the vectorial Lie superalgebra~that sends
$\cF(x, y)$ into $I(x)$, is $\fvect(x)$-invariant, and therefore the
semidirect sum $\Ann(x)\fvect(x, y) \inplus \fvect(x)$ is well defined.

\sssbegin[What $\fu^{\fvect(m|n)}_{*}$ is]{Lemma}[What $\fu^{\fvect(m|n)}_{*}$ is]\label{sss1.2.1}
We have $\fu^{\fvect(m|n)}_{*}\simeq\Ann(x)\fvect(x, y) \inplus 
\fvect(x)$.\end{Lemma}

\begin{proof} Straightforward.\end{proof}

\ssec{
Singular subspaces}\index{Subspace $V \subset \fg_{-1}$, singular} We say that a~subspace $V \subset \fg_{-1}$ is
\textit{singular for $\fg= \fh(2m|n)$ or $\fh'(0|n)$}, if either
of the following holds:
\[
\renewcommand{\arraystretch}{1.2}
\begin{array}{ll}
(a)&\dim V^{\perp}/(V \cap V^{\perp}) = \eps,\\
(b)&V^{\perp}/V = \Ker(\omega|_{\fu_{-1}}) \neq 0\text{ and $\dim
V^{\perp}/(V\cap V^{\perp}) = 2\eps$.}
\end{array}
\]

A~subspace $V \subset \fg_{-1}$ is said to be ~\textit{singular for $\fg = \fle(n)$ or $\fsle'(n)$},\index{Subspace $V \subset \fg_{-1}$, singular} if 
\[
\text{the space $V \cap V^{\perp} = \Ker(\omega|_{\fu_{-1}})$ is even and $\dim V^{\perp}/(V\cap V^{\perp})
= 1 + \eps$.}
\] 

\sssbegin[On $\dim
V^{\perp}/(V\cap V^{\perp})$]{Remark} [On $\dim
V^{\perp}/(V\cap V^{\perp})$]
Since $\omega|_{V\cap V^{\perp}} = 0$ and
$\omega|_{V^{\perp}/(V\cap V^{\perp})}$ is non-degenerate, it follows that $\dim
V^{\perp}/(V\cap V^{\perp}) \neq \eps$ if $\fg = \fle(n)$ or
$\fsle'(n)$.\end{Remark}

These definitions settle the case $4^o$ of \eqref{5cases}, in which $0\subsetneq\fu_{-1}\subsetneq
\fg_{-1}$. Indeed, by Lemma \ref{L1.1}, a~maximal $\fu$ with a~given $\fu_{-1}$ can only be the
prolong $\fu^\fg_*$ of $\fu_{-1}$ in $\fg$, and this prolong is computed inside $\fvect$ once and
for all, thanks to the equality $\fu^\fg_{*}=\fu^{\fvect(\fg^*_{-1})}\cap \fg$ of item
2\textup{)} of that Lemma. The next Theorem carries the computation out for every ambient of the
main list \eqref{mainList}; in it, $V:=\fu_{-1}$, and $V^\perp$ is taken with respect to the form
$\omega$ that $\fg$ preserves. For the ambients preserving a~$2$-form, the answer depends on
whether $V$ is singular in the sense just defined: for a~singular $V$, the prolong $\fu^\fg_*$ is
contained in the prolong of a~larger subsuperspace and hence is \textbf{not} maximal, see
Theorem \ref{ThMax}, where the maximality of the subalgebras listed is proved.

\ssbegin[What is $\fu^{\fg}$ for simple vectorial superalgebras $\fg$]{Theorem}[What is $\fu^{\fg}$ for simple vectorial superalgebras $\fg$]\label{Th1.4}
Table \eqref{table4} lists
$\fu^{\fg}$ for ``simple\rq\rq\ vectorial superalgebras~$\fg$.
\end{Theorem}

\section{Depth 2}

\ssec{The contact gradings}\label{Sd2contact} For most of the classical
domains $X=S/P$, where $S$ is a~simple Lie group and $P$ its parabolic subgroup, we see that the Lie algebra~$\fs$ of $S$ is
\be\label{(2.1.)}
\fs\simeq(\fg_{-1}, \fg_0)_*. 
\ee
When \eqref{(2.1.)} fails, $\fs$ is a~proper
subalgebra~of $(\fg_{-1}, \fg_0)_*$. It is of interest therefore to solve the following \textbf{Two problems}:
\begin{equation}\label{2prob}
\begin{minipage}[l]{15cm}
(a) List all the
cases where, having started from a~simple Lie(super)algebra\\
\phantom{.}\hskip 6mm 
$\fs=\mathop{\oplus}_{i\geq -d}\ \fs_i$, we have the following
analog of isomorphism~\eqref{(2.1.)}.\\
(b) Find out \textit{what} is the ``total prolong\rq\rq\  $(\fs_-)_*:=(\fs_-, \tilde \fs_0)_*$ of $\fs_-$,
where  $\tilde
\fs_0:=\fder_0(\fs_-)$\\
\phantom{.}\hskip 6mm is the Lie superalgebra~of derivations of $\fs_-$ preservingits $\Zee$-grading.
\end{minipage}
\end{equation}


For the simple Lie algebras $\fs$, Yamaguchi summarized
the known results:

\sssbegin[Cartan prolong of $\fg_-$ for the simple finite-dimensional Lie algebra  $\fg$]
{Theorem}[Cartan prolong of $\fg_-$ for the simple finite-dimensional Lie algebra  $\fg$(\cite{Y})]\label{ThYa} The equality
$(\fs_-)_*=\fs$ holds almost always. The exceptions are

\textup{1)} $\fs$ with the grading of depth $d=1$, in which case
$(\fs_-)_*=\fvect(\fs_-^*)$;

\textup{2)} $\fs$ with the grading of depth $d=2$ and $\dim\fs_{-2}=1$,
i.e., with the ``contact\rq\rq\  grading, in which case
$(\fs_-)_*=\fk(\fs_-^*)$ (these cases correspond to exclusion of the
nodes on the Dynkin graph connected with the node for the maximal
root on the extended graph);

\textup{3)} $\fs$ is either $\fsl(n+1)$ or $\fsp(2n)$ with the grading
determined by ``excluding\rq\rq\  the first and the $i$th of simple
coroots, where $1<i<n$ for $\fsl(n+1)$ and $i=n$ for~$\fsp(2n)$.
\textup{(Observe that $d=2$ with $\dim \fs_{-2}>1$ for $\fsl(n+1)$ and $d=3$
for $\fsp(2n)$.)}

Moreover, the isomorphism $\fs\simeq(\fs_-, \fs_0)_*$, see eq.~\eqref{2prob} 
holds almost always, even in the 
exceptional cases~$1)-3)$. The exceptional cases --- the ones where
 a~projective action is possible --- are $\fsl(n+1)$ or $\fsp(2n+2)$
with the grading determined by ``excluding\rq\rq\  only one (the first)
simple coroot. \end{Theorem}

For the derivation of Theorem \ref{ThYa} from Cartan's classification, see \cite{LS1}. The
mysterious at this stage case 3\textup{)} appears from the projectivization argument, due to
Arnold, reproduced in \S\ref{s2.1.2proof}: the gradings of $\fsl(n+1)$ and of $\fsp(2n+2)$ obtained
by ``excluding\rq\rq\ the first simple coroot are precisely the ones for which the action on the
projectivization, not on the space itself, is the relevant one.

Here is the superization of Theorem \ref{ThYa}:
Let $\fg$ be a~simple finite-dimensional Lie superalgebra~(or the derivation
algebra~of a~simple superalgebra).

A) Let $\fg_\ev$ be reductive. There are two types of gradings\index{Grading, contact}
of $\fg$ naturally referred to as \textit{contact} ones:
\[
\fg=\renewcommand{\arraystretch}{1.4}
\begin{cases}
\fg_{-2}\oplus \fg_{-1}\oplus \fg_{0}\oplus \fg_{1}\oplus \ldots \;\text{ with
}\dim\fg_{-2}=1&\text{is embeddable into $\fk(\sdim (\fg_{-2}\oplus \fg_{-1}))$}\\
\fg_{-2}\oplus \fg_{-1}\oplus \fg_{0}\oplus \fg_{1}\oplus \ldots \;\text{ with
}\dim\fg_{-2}=\eps&\text{is embeddable into $\fm(\sdim (\fg_{-2}\oplus \fg_{-1}))$}.
\end{cases}
\]
In addition to these \textit{standard} $\fk$-contact and
$\fm$-contact gradings of $\fg$, we have \textbf{nonstandard} $\fk$-contact
and $\fm$-contact gradings related to embeddings of $\fg$ into nonstandardly
graded $\fk(2n+1|m)$ and into $\fm$ and~$\fb_{a, b}$. The nonstandard gradings that occur are
the ones listed in table \eqref{1.4b}; in this text they are
\[
\fsq(n+1)\subset\fk(2n+1|2n+2;1), \quad \fk(1|2n;n)\subset\fsvect(\Lambda(n)), \quad
\fm(n;n)\subset\fsvect(\Pi(\Lambda(n))),
\]
see table \eqref{tabl} and \S\ref{s2.1.2proof} for the first of them, table \eqref{table5} and
Lemma \ref{LinfDim} for the other two, and the regradings $\fle(n;r)$, $\fb_\lambda(n;n)$,
$\fb'_\infty(n;n)$,~$\fb'_1(n;n)$. In each of them, the nonstandard case is reduced to the standard
one by the bigrading argument used before table \eqref{1.41}.

B) Let $\fg_\ev$ be not reductive. Then, $\fg$ is vectorial, i.e., $\fg$ is one of $\fvect(0|n)$, $\fsvect(0|n)$, $\widetilde{\fsvect}(0|n)$, and $\fh'(0|n)$, and its gradings to be considered are the regradings of these listed in table \eqref{1.4b}.

Two remarks are in order here. First, in case $5^o$ of \eqref{5cases} the subalgebra
$\fu=\mathop{\oplus}\limits_i\fg_{2i}$ is of a~type new for this paper. Indeed, if
$\fg=\oplus_i\fg_i$ is $\Zee$-graded, then $\sigma(x):=(-1)^{i}x$ for $x\in\fg_i$ is an involutive
automorphism of $\fg$ and
\be\label{fixedPt}
\fu=\mathop{\oplus}\limits_{i}\fg_{2i}=\fg^\sigma
\ee
is its fixed-point subalgebra; hence, $\fu$ is maximal if and only if the $\fu$-module
$\mathop{\oplus}_i\fg_{2i+1}$ --- the $(-1)$-eigenspace of $\sigma$ --- is irreducible, the standard criterion for the fixed points of an involution, easy to check case by case. Such a~$\fu$ is not of the form $(\fu_{-1}, \fu_0)_*\cap\fg$, nor is it singled out by a~condition on $\fu_{-1}$ and $\fu_0$ alone. Observe, however, that $\fu$ of the form \eqref{fixedPt} has $\fu_-\neq0$, so it belongs to the case $3^o$ or to the case $6^o$, not to the case $5^o$, in which Lemma \ref{Lobvious} says that $\fu$ is maximal if and only if $\fu=\oplus_{i\geq0}\fg_i$.

Second, on the depth $3$: among the ambients, only $\fm\fb(3|8)$ --- a~regrading of $\fm\fb(4|5)$ --- is of depth $3$, see \eqref{depth3}; among the finite-dimensional $\fg$ of Theorem \ref{ThYa}, only $\fsp(2n)$ with the grading of case 3\textup{)} and the superizations of that grading. In both cases, the grading is a~regrading of one of depth $\leq2$ and is dealt with as such; this is why no separate section is devoted to the depth~$3$. The next Subsection explains why neither this postponement nor the one announced before table \eqref{1.41} --- the cases $N=6$, $9$, $10$,
$16$ of table \eqref{1.4b} deferred from the depth $1$ to the depth $2$ --- changes anything: the depth is not a~parameter of the answer.

\ssec{The depth is not a~parameter of the answer}\label{Sdepthred} Here we
explain why the ambients of depth $2$ can be treated by the depth-$1$
machinery of \S\ref{Sdeg1amb}, and those of depth $3$ by the depth-$2$ one. The
reason is that, for a~\textbf{maximal} graded subalgebra, the components
$\fu_{-j}$ with $j\geq2$ carry no information of their own: they are computed
from $\fu_{-1}$ and $\fu_{\geq0}$, not chosen.

Recall first that the negative part of a~W-grading is generated by its
$(-1)$st component: the Weisfeiler filtration is defined by \eqref{Wgen}, and
therefore the associated graded algebra~satisfies
\be\label{Wgen2}
\fg_{-i-1}=[\fg_{-1}, \fg_{-i}]\;\text{ for $1\leq i<d$}.
\ee
For the ambients of depth $>1$, the bracket in \eqref{Wgen2} is not just
non-zero, it is the only equivariant map there is:

\sssbegin[The bracket on $\fg_-$ is unique up to a~factor]{Statement}[The
bracket on $\fg_-$ is unique up to a~factor]\label{Sbrmult} Let $\fg$ be one of
the ambients of depth $>1$ of table \eqref{nonpos}, i.e., one of
$\fv\fle(3|6)$, $\fk\fas(1|6)$, $\fm\fb(3|8)$, $\fk\fle(5|10)$, or let $\fg$ be
$\fk(1|2n)$ in its contact grading. Then, the $\fg_0$-module $\fg_{-2}$ occurs
in $S^2(\fg_{-1})$ --- the bracket of two odd elements is symmetric on the
underlying space --- with multiplicity exactly $1$; for $\fm\fb(3|8)$ the
$\fg_0$-module $\fg_{-3}$ occurs in $\fg_{-1}\otimes \fg_{-2}$ with
multiplicity exactly $1$ as well. Hence, each bracket
$\fg_{-1}\times\fg_{-i}\tto\fg_{-i-1}$ is determined by the $\fg_0$-action up
to a~factor, and by \eqref{Wgen2} that factor is non-zero.

The multiplicities were computed for each of these five gradings and for
both parities of~$\fg_{-1}$.
\end{Statement}

They come from Freudenthal's formula applied to the reductive parts
$\fsl(3)\oplus\fsl(2)$, $\fo(6)\simeq\fsl(4)$ and $\fsl(5)$ of the
respective~$\fg_0$. Observe that the parity decides: in
$\Lambda^2(\fg_{-1})$, which is what the bracket would factor through were
$\fg_{-1}$ even, the module $\fg_{-2}$ of $\fm\fb(3|8)$ does \textit{not}
occur.

\sssbegin[A~maximal graded subalgebra~is self-normalizing]{Lemma}[A~maximal
graded subalgebra~is self-normalizing]\label{Lselfnorm} Let $\fg$ be simple and
$\Zee$-graded, and let $\fu$ be a~maximal graded subalgebra~of~$\fg$. Then,
$\fu=N_\fg(\fu):=\{x\in\fg\mid [x, \fu]\subseteq\fu\}$.
\end{Lemma}

\begin{proof} The normalizer of a~graded subspace is graded: if
$x=\sum x_i$ normalizes $\fu$ and $u\in\fu_k$, then $[x, u]\in\fu$, and the
$(i+k)$th component of $[x, u]$ is $[x_i, u]$, so $[x_i, u]\in\fu$ for every~$i$. Hence, $N_\fg(\fu)$ is a~graded subalgebra~containing~$\fu$. Were it
distinct from $\fu$, the maximality of $\fu$ would give $N_\fg(\fu)=\fg$, i.e.,
$\fu$ would be an ideal of $\fg$, which is impossible for $\fg$ simple and
$0\neq\fu\neq\fg$.
\end{proof}

\sssbegin[The components $\fu_{-j}$ with $j\geq2$ are computed, not
chosen]{Lemma}[The components $\fu_{-j}$ with $j\geq2$ are computed, not
chosen]\label{Ldeep} Let $\fg=\oplus_{i\geq -d}\ \fg_i$ be simple and W-graded
of depth~$d\leq3$. Fix a~subspace $\fu_{-1}\subseteq\fg_{-1}$ and a~graded
subalgebra~$\fu_{\geq0}\subseteq\fg_{\geq0}$, and let $\cS$ be the set of the
tuples $(U_2, \dots , U_d)$, where $U_j\subseteq\fg_{-j}$, such that
\be\label{compl}
\fu_{-1}\oplus U_2\oplus\dots \oplus U_d\oplus\fu_{\geq0}\;\text{ is
a~subalgebra~of $\fg$.}
\ee
Then, $\cS$ is closed under componentwise sums; hence, if $\cS\neq\emptyset$,
it has a~largest element. If $\fu$ is a~maximal graded subalgebra~of $\fg$,
then $(\fu_{-2}, \dots , \fu_{-d})$ is the largest element of the $\cS$ built
from $\fu_{-1}$ and~$\fu_{\geq0}$.
\end{Lemma}

\begin{proof} Spelled out, \eqref{compl} splits into the conditions that do
not mention the $U_j$ at all --- namely, that $\fu_{\geq0}$ be a~subalgebra~and
that $[\fu_{-1}, \fu_k]\subseteq\fu_{k-1}$ for $k\geq1$ and $[\fu_{-1},
\fu_0]\subseteq\fu_{-1}$, all of them conditions on the given pair
$(\fu_{-1}, \fu_{\geq0})$ alone, which either hold or leave $\cS$ empty ---
and the following list, in which we write $U_1:=\fu_{-1}$:

$1)$ $[\fu_{-1}, \fu_{-1}]\subseteq U_2$;

$2)$ $[\fu_{-1}, U_2]\subseteq U_3$, an empty condition for $d=2$;

$3)$ $[U_i, U_j]\subseteq\fg_{-i-j}=0$ for $i$, $j\geq2$, an empty condition
because $i+j\geq4>d$;

$4)$ $[U_j, \fu_k]\subseteq \fu_{k-j}$ for $k\geq0$ and $2\leq j\leq d$, where
$\fu_{k-j}$ means $U_{j-k}$ if~$k<j$.

Let $(U_j)$ and $(U'_j)$ both lie in~$\cS$. Each of the conditions above is
additive in the left-hand side and monotone in the right-hand side, and no
condition couples $U_i$ with $U'_j$ for $i$, $j\geq2$ thanks to~$3)$. Thus,
$1)$ holds for $(U_j+U'_j)$ because it holds for $(U_j)$; $2)$ holds because
$[\fu_{-1}, U_2+U'_2]\subseteq U_3+U'_3$; and $4)$ holds because $[U_j+U'_j,
\fu_k]$ lands in $\fu_{k-j}$ for $k\geq j$ and in $U_{j-k}+U'_{j-k}$ for
$k<j$, each summand landing in its own place.
Hence, $\cS$ is closed under sums, and, the $\fg_{-j}$ being
finite-dimensional, $\cS$ has a~largest element $(\widehat U_j)$.

Let now $\fu$ be maximal graded, and let $\widehat\fu$ be the subalgebra~built
by \eqref{compl} from $\fu_{-1}$, $\fu_{\geq0}$ and the largest element of
$\cS$; so $\fu\subseteq\widehat\fu$. If $\widehat\fu=\fg$, then
$\fu_{-1}=\fg_{-1}$ and $\fu_{\geq0}=\fg_{\geq0}$, but then $\fu\supseteq\fg_-$
by~\eqref{Wgen2} and $\fu=\fg$, which is excluded. So $\widehat\fu$ is
a~proper graded subalgebra~containing $\fu$, whence $\widehat\fu=\fu$.
\end{proof}

\sssbegin[Reduction of the depth]{Corollary}[Reduction of the
depth]\label{Cdepthred} Let $\fg$ be simple and W-graded of depth $d\leq3$, and
let $\fu$ be a~maximal graded subalgebra~of~$\fg$. Then,

$1)$ $\fu$ is determined by the pair $(\fu_{-1}, \fu_{\geq0})$ --- precisely
the data of the depth-$1$ case; the depth contributes no parameter of its own;

$2)$ If $\fu_{-1}=\fg_{-1}$, i.e., in the cases $1^o$--$3^o$
of~\eqref{5cases}, then $\fu_-=\fg_-$ whatever the depth, and the argument of
\S\ref{Sdeg1amb} applies verbatim;

$3)$ In the case $5^o$, Lemma~\ref{Lobvious} settles $\fu$ in any depth;

$4)$ Hence the depth is felt in the case $4^o$ only, and the only value of the
pair $(\fu_{-1}, \fu_{\geq0})$ which the depth $1$ does not offer at all is
$\fu_{-1}=0$ with $\fu_-\neq0$, i.e., the case~$6^o$.
\end{Corollary}

\begin{proof} Item $2)$ is~\eqref{Wgen2} applied to $\fu_{-1}=\fg_{-1}$: the
subalgebra~generated by $\fg_{-1}$ is the whole of~$\fg_-$. Item $1)$ is
Lemma~\ref{Ldeep} together with item $2)$, item $3)$ is Lemma~\ref{Lobvious}, and
item $4)$ is the list~\eqref{5cases} read with items $2)$ and $3)$ in hand: the
cases $1^o$--$3^o$ and $5^o$ are depth-free, and the case $6^o$ is vacuous for
$d=1$ because there $\fu_{-1}=0$ means~$\fu_-=0$. \end{proof}

For the contact and pericontact gradings --- the ones the whole of
\S\ref{Sd2contact} is about --- the residue left by Corollary~\ref{Cdepthred},
namely the case $4^o$, shrinks to a~single bit of information.

\sssbegin[What the depth $2$ adds in the case $4^o$]{Proposition}[What the
depth $2$ adds in the case $4^o$]\label{Pcontactbit} Let $d=2$ and
$\sdim\fg_{-2}=(1|0)$ or $(0|1)$, i.e., let the grading be $\fk$-contact or
$\fm$-contact, and let $\omega\colon\fg_{-1}\times\fg_{-1}\tto\fg_{-2}$ be the
bracket. Let $\fu$ be a~maximal graded subalgebra~of $\fg$ and~$V:=\fu_{-1}$.
Then, $\fu_{-2}$ is either $0$ or $\fg_{-2}$, and $\fu_{-2}=\fg_{-2}$ whenever
$\omega|_V\neq0$.
\end{Proposition}

\begin{proof} Being a~subspace of a~superspace of $\sdim$ $(1|0)$ or $(0|1)$,
the component $\fu_{-2}$ is $0$ or $\fg_{-2}$; and $\omega(V, V)=[\fu_{-1},
\fu_{-1}]\subseteq\fu_{-2}$. \end{proof}

Thus, at the depth $2$ the classification runs over the same list of pairs
$(V, \fu_0)$ as at the depth $1$, subject to one extra requirement --- that the
pair be closed under $\omega$ --- and the alternative ``$\omega|_V$
non-degenerate or $\omega|_V\equiv0$\rq\rq\ around which \S\ref{Sproofs} is
organized is nothing but Proposition~\ref{Pcontactbit} in use: in the first
case $\fu\supset\fg_{-2}$ and $V$ is recovered as a~non-degenerate subspace, in
the second $\fu_-=V$ is abelian and one passes to $W=V^\perp/V$, of depth $1$
again.

At the depth $3$, Lemma~\ref{Ldeep} applies with $j=2$, $3$; what makes it
work is that $[\fg_{-j}, \fg_{-j'}]\subseteq\fg_{-j-j'}=0$ whenever $j+j'>3$,
so that the components $\fg_{-2}$ and $\fg_{-3}$ do not bracket with one
another. Since
$\fm\fb(3|8)$ is the only ambient of depth $3$, and it is a~regrading of
$\fm\fb(4|5)$, which contains the grading operator of its standard grading, the
bigrading argument given before table \eqref{1.41} makes a~maximal graded
$\fu$ of the cases $1^o$--$4^o$ graded with respect to that standard grading of
depth $2$ as well; Lemma~\ref{Ldeep} then returns $\fu_{-2}$ and~$\fu_{-3}$.
This is the promised proof of the assertion made at the end of
\S\ref{Sd2contact}.

\sssbegin[What the reduction does not give]{Remark}[What the reduction does not
give]\label{Rdepthnot} Three limitations are worth stating.

$1)$ The reduction is one of \textit{data}, not of \textit{answers}: at the
depth $1$ every $\fu_0$-submodule $V\subseteq\fg_{-1}$ is an admissible
$\fu_{-1}$, at the depth $2$ the pair $(V, \fu_{-2})$ must in addition be
closed under~$\omega$. The list of pairs to run over is the same, the
conditions cutting it down are strictly stronger.

$2)$ The case $6^o$ is not reduced: it is the one case the depth $1$ does not
see, and, for the infinite-dimensional ambients, we have no proof that it is
vacuous --- see the caveat in \S\ref{Sd2contact} and the answers for the three
exceptional finite-dimensional ambients in \cite{LOdd}, where it is not.

$3)$ The bigrading argument used for the depth $3$ needs $\fu$ to contain the
grading operator of the standard grading. In the cases $1^o$ and $2^o$ this is
automatic, $\fu_{\leq0}$ being then the whole of $\fg_{\leq0}$; in the cases
$3^o$ and $4^o$ it is the assertion made before table~\eqref{1.41}.
\end{Remark}

\section{Proofs}\label{Sproofs}

The case where $\fu\supset \fg_{-1}\oplus \fg_{0}$ is settled by Lemma
\ref{L3.1.1}, proved in \S\ref{Sdeg1amb}. The trichotomy of the cases A\textup{)}, B\textup{)}, C\textup{)} below is the one stated in
\S\ref{Sdeg1amb}, where the answers are given: \S\ref{Sdeg1amb} lists them, the present section
proves them. Namely, \S\ref{SAred} carries out the two cases $\omega|_V$ non-degenerate and
$\omega|_V\equiv 0$ of the case A\textup{)} and arrives at \eqref{omega0} with $W=V^\perp/V$,
whereas \S\ref{SCirred} settles the case C\textup{)}.

\ssbegin[The maximal subalgebras containing $\fg_{-1}$ but not
containing $\fg_{0}$]{Lemma}[The maximal subalgebras containing $\fg_{-1}$ but not
containing $\fg_{0}$]\label{L3.2.1} Let $\fu$ be a~ subsuperalgebra~of $\fg$
such that $\fu_{-1} = \fg_{-1}$, and $\fu_{0} \subsetneq \fg_{0}$.
 The Cartan prolong $(\fg_{-1}, \fu_{0})_*$ is
maximal in $\fg$ if and only if

\textup{1)} $\fu_{0}$ is maximal in $\fg_{0}$
and

\textup{2)} $[\fu_{1},\fg_{0}] = \fg_{1}$. 

For simple Lie superalgebras
we consider, condition $2)$ implies that

--- $\fu_{1}\neq 0$ if $\fg= \fsvect(m|n)$, $\fh(2m|n)$, $\fh'(0|n)$
for $n \neq 6$, or $\fsle'(n)$;

--- the projection of $\fu_{1}$ to each $\fg_{0}$-irreducible
submodules of $\fg_{1}$ is nonzero if $\fg = \fvect(m|n)$,
$\fh'(0|6)$;

--- $\fu_{1}\not \subset \fsle' (n)_{1}$ if $\fg = \fle(n)$.
\end{Lemma}

\begin{proof} Clearly, every subsuperalgebra~$\fu$ with the required
properties is contained in the Cartan prolong $(\fg_{-1}, \fu_{0})_*
=\mathop{\oplus}\limits_{i\geq 0} \fu_{i}$.

The necessity of condition 1) follows from the fact that
if $\fu_{0}\subset\fv_{0}\subset\fg_{0}$ with strict inclusions for
some subsuperalgebra~$\fv_{0}$, then $\fu \subset (\fg_{-1},
\fv_{0})_* \subset \fg$ also with strict inclusions.

The necessity of condition 2) follows from the fact that if it
fails, then $\fu$ is contained in one of the subsuperalgebras
described in Lemma~\ref{L3.1.1}.

Sufficiency is proved similarly to the proof of item 1\textup{)} of Lemma \ref{L3.1.1}.
\end{proof}

Thus, $\fu_{0}$ is maximal in~$\fg_{0}$. The following 3 cases may occur,
cf. \cite{ShM}:

A) $\fu_{0}$ is a~reducible subsuperalgebra~in $\fg_{0}$, i.e., the
$\fu_{0}$-module $\fg_{-1}$ is reducible;

B) $\fu_{0}$ is an irreducible subsuperalgebra~in $\fg_{0}$,
moreover, $\fu_{0}$ is either almost simple or a~central extension of
an almost simple Lie superalgebra;

C) $\fu_{0}$ is an irreducible subsuperalgebra~in $\fg_{0}$,
moreover, $\fu_{0}$ is neither almost simple nor a~central extension
of an almost simple Lie superalgebra.

Let us consider each case separately.

\sssec{When
$\fu_{0}$ is a~reducible subsuperalgebra~in $\fg_{0}$}\label{SAred}
Consider the case where $\fg$ preserves a~2-form~$\omega$.

1) $\omega|_{V}$ is non-degenerate. Clearly, if $\fu_{0}$ preserves
$V$, it preserves $V^\perp$, $V\cap V^\perp$ and $V + V^\perp$, moreover
$\fg_{-1} = V\oplus V^\perp$.
Hence, we have
\[
\begin{array}{l}
\fst(V) \cap \fg_{0} = \fosp(V) \oplus \fosp (V^\perp),\\
\faut\cap \fg
= \fh(V)\oplus \fh(V^\perp) \text{ for $\fg = \fh(2m|n)$ or
}\fh'(0|n);\\
\fst(V) \cap \fg_{0} = \fpe(V) \oplus \fpe(V^\perp),\\
\faut \cap \fg=
\fle(V)\oplus \fle(V^\perp) \text{ for }\fg= \fle(n).
\end{array}
\]
Finally, for $\fg= \fsle(n)$, we have
\[
\faut \cap \fg= (\fsle(V)\oplus \fsle(V^\perp)) \inplus \Cee\cdot D,
\]
where
\[
D = l(\id |_{V_{\ev}} -\id |_{V_{\od}}) -k(\id|_{V_{\ev}^\perp}
-\id|_{V_{\od}^\perp})\;\text{ (for the definition of $k$ and $l$, see eq.~\eqref{1.2}).}
\]

2) $\omega|_{V}\equiv 0$. Then, $V^\perp\supset V$. Consider the
realization of $\fg$ by generating functions with the following
coordinates: let $p = (p_1, \ldots , p_k)$ span $V$; let $\zeta= (\zeta_1,
\ldots , \zeta_l)$ span the complement of $V$ in $V^\perp$; let $q =
(q_1, \ldots , q_k)$ be the dual to~$p$.

Let us deviate from the traditional notation with Latin and Greek
letters for the even and the odd; assume that each coordinate may be
of any parity. Then, $\ad_{p_i} = (-1)^{p(p_i)}\partial_{q_i}$.

The invariance of $V$ with respect to $\fu_{0}$ means that
\[
\fu_{0} = \Span(p_ip_j,\; p_i\zeta_j, \;\zeta_i\zeta_j) + \Span(q_ip_j) \;
\text{ for every }\; i,j.
\]
But then
\[
\renewcommand{\arraystretch}{1.4}
\begin{array}{l}
 (\fg_{-1}, \fu_{0})_{*} = \faut(V) \cap \fg = \cF(p,
\zeta)/\Cee\cdot 1 \inplus \Span(q) \otimes \cF(p) =\\
\cF(p) \otimes
\cF(\zeta)/\Cee\cdot 1 \inplus \Span(q) \otimes \cF(p).\\
\end{array}
\]
Moreover,
\[
\textstyle \{q_if(p), g(p,\zeta)\} = \pm f(P)\pderf{g(p,\zeta)}{p_i},
\]
where the sign only depends on the parity of $f$ and the type of $\fg$
($\fh$ or~$\fle$). Hence, the functions of the form $\Span(q) \otimes
\cF(p)$ generate a~subalgebra~isomorphic to $\fvect(p)$ with the
natural action on the ideal $\cF(p) \otimes F(\zeta)/\Cee\cdot 1$ (as
on the space of vector-valued functions).

Since $\omega$ is non-degenerate on $V^{\perp}/V$, the space
$\cF(\zeta)$ generates a~subalgebra~isomorphic to either $\fpo(\zeta)$
if $\fg = \fh(2m|n)$ or $\fb(\zeta)$ if~$\fg= \fle(n)$.

Thus,
\[
\renewcommand{\arraystretch}{1.4}
\begin{array}{l}
 \faut(V) \cap \fh(2m|n) = (\cF(V) \otimes \fpo(V^{\perp} /V))/\Cee\cdot
1\inplus \fvect(V),\\
\faut(V) \cap \fle(2m|n) = (\cF(V) \otimes \fb(V^{\perp}/V))/\Cee\cdot
1\inplus \fvect(V).\\
\end{array}
\]

\begin{proof} This theorem is a~direct consequence of Lemma~\ref{L3.2.1} 
and the explicit description of maximal reducible subalgebras of
matrix Lie superalgebras, see Main Theorem of \cite{ShM}. The
explicit form of the subsuperalgebra~$\fu$ demonstrates that
conditions of Lemma~\ref{L3.2.1} are satisfied.\end{proof}

\sssec{When
$\fu_{0}\subset\fg_{0}$ is irreducible,
not (a~central extension of an) almost simple}\label{SCirred}
Recall the list of maximal subsuperalgebras $\fu_{0}\subset
\fg_{0}$ of this form (\cite{ShM}):

$\bullet$ if $\fg_{0} = \fgl(m|n)$ or
$\fsl(m|n)$ (i.e., $\fg= \fvect(m|n)$ or $\fsvect(m|n)$), then
$\fu_{0}$ belongs to one of the four types of Lie superalgebras:

1) $\fgl(V_1) \odot \fgl(V_2)$;

2) $\fq(V_1)\odot \fq(V_2)$;

3) $\fgl(V) \otimes \Lambda(k) \inplus \fvect(0|k)$;

4) $\fhei(0|2n) \oplus \fo(2n)$;

$\bullet$ if $\fg_{0} = \fosp(m|2n)$, $\fpe(n)$ or $\fspe(n)$ (i.e.,
$\fg= \fh(2n|m)$ or $\fle(n)$ or $\fsle(n)$, respectively), then
$\fu_{0}$ is of the form $\faut(V_1) \odot \faut(V_2)$.

First, consider the case $\fu_{0} = \fu^{(1)}(V) \odot \fu^{(2)}(U)$. More
exactly, this
means that $\fg_{-1} = V \otimes U$ on which $\fu^{(1)}:=\fu^{(1)}(V)$ acts by
operators of
the form $a\otimes 1$ and $\fu^{(2)}:=\fu^{(2)}(U)$ acts by operators of the
form~$1\otimes
b$. We additionally assume that $\dim U, \;\dim V \neq (1|0),\; (0|1)$.

\ssbegin[On existence of rank $1$ operators]
{Lemma}[On existence of rank $1$ operators]\label{Lrank1} Set $\fu= \mathop{\oplus}_{i\geq -1} \fu_{i}=
(\fg_{-1}, \fu_{0})_*$. If $\fu_{1} \neq 0$, then at least one of the
subalgebras, $\fu^{(1)}(V)$ or $\fu^{(2)}(U)$, contains an ideal
spanned, as linear superspace, by rank $1$ operators. \end{Lemma}

\begin{proof} Since $\tilde \fu_{0} = [\fu_{-1}, \fu_{1}] \subset \fu_{0}$,
we can represent every $h\in \fu_{0}$ in the form \[h = A\otimes 1 +
1\otimes B,\text{~~ where $A\in \fu^{(1)}, B\in \fu^{(2)}$.}\]
Denote the projections of $\fu_{0}$ on $\fu^{(1)}$ and $\fu^{(2)}$
by $I_{\fu^{(1)}}$ and $I_{\fu^{(2)}}$, respectively. It is clear
that $I_{\fu^{(1)}}$ and $I_{\fu^{(2)}}$ are ideals. Let us
describe them.

Let $F\in \fu_{1}$ and let $u\otimes v, u_{1}\otimes v_{1}\in
\fu_{-1}$. Then, 
\be\label{(i.)}
{}[F, u\otimes v]=A(u, v)\otimes 1+1\otimes B(u, v), 
\ee
\be\label{(ii.)}
{}[F, u_{1}\otimes v_{1}]= A(u_{1}, v_{1})\otimes 1 + 1\otimes B(u_{1},
v_{1}), 
\ee
where $A: U \otimes V\tto \End U$ and $B: U\otimes V\tto \End V$ are
linear operators.

Thanks to the Jacobi identity we have
\be\label{(iii.)}
[[F, u\otimes v], u_{1}\otimes v_{1}] = \pm [[F, u_{1}\otimes v_{1}], u\otimes
v] 
\ee
(the sign depends on parities of the terms involved; it is
irrelevant to us now).

Let us substitute \eqref{(i.)} and \eqref{(ii.)} into \eqref{(iii.)}:
\be\label{(iv.)}
A(u, v)u_{1}\otimes v_{1} + u_{1}\otimes B(u,v)v_{1} = \pm A(u_{1},
v_{1})u\otimes v
\pm u\otimes B(u_{1}, v_{1})v. 
\ee
By the condition on dimensions of $U$ and $V$ we may assume that the
pairs $u, u_{1}$ and $v, v_{1}$ are noncollinear.

But in this case equation \eqref{(iv.)} takes place only if
\[
A(u, v)u_{1} = \alpha_{u, v}(u_{1})u,\quad B(u, v)v_{1} = \beta_{u, v}(v_{1})v,
\]
where $\alpha_{u, v}$ and $\beta_{u, v}$ are linear forms on $U$ and
$V$, respectively. In other words,
\[
\rk A(u,v), \rk B(u,v) \leq 1\;\text{
for arbitrary $u\in U$ and $v\in V$. }
\]
This means that $I_{\fu^{(1)}}$ and
$I_{\fu^{(2)}}$ are spanned, as superspaces, by rank 1 operators.
\end{proof}

\sssbegin[No maximal Lie
subsuperalgebras of type $2$C]
{Corollary}[No maximal Lie
subsuperalgebras of type $2$C]\label{C3.3.1} If $\fg$ is isomorphic to
$\fsvect(m|n)$ for $m \neq n$, or $\fh(2m|n)$, or $\fh'(0|n)$, or
$\fle(n)$, or $\fs\fle'(n)$, then $\fg$ does not contain maximal Lie
subsuperalgebras of type $2$C. \end{Corollary}

Let now $\fg \simeq \fvect(m|n)$ or $\fs\fvect(n|n)$.

\ssec{
An embedding $\fp\fgl(V_{1}\oplus V_2) \subset
\fvect(V_{1}^*\otimes V_2)$} The $\Zee$-grading of the direct sum of
linear superspaces $V_{1}\oplus V_2$ given by $\deg V_{1} = 1$, $\deg
V_{2} = 0$ induces a~$\Zee$-grading of the Lie superalgebra~$\fa =
\fgl(V_{1}\oplus V_2)$ of the form $\fa =\fa_{-1}\oplus \fa_{0}\oplus
\fa_{1}$, where
\[
\fa_{-1}= V_{1}^*\otimes V_2, \quad
\fa_{0}= \fgl(V_{1})\oplus \fgl(V_2), \quad \fa_{1}= V_{1}\otimes V_2^*.
\]
So, we have an embedding $\varphi : \fa\tto \fvect(V_{1}^*\otimes
V_2)$. Set $\fu= \varphi(\fa)$. Clearly, $\fu= \fp\fgl(V_{1}\oplus
V_2)$ with
\[
\fu_{-1}=\fa_{-1}, \quad \fu_{0}=\fgl(V_{1}\odot V_2), \quad
\fu_{1}=\fa_{1}.
\]
Let $\dim V_i = k_i + l_i\eps$ for~$i= 1, 2$. If either $k_1\neq l_1$
or $k_2\neq l_2$, then $\fp\fgl(V_{1}\oplus V_2)
\simeq\fsl(V_{1}\oplus V_2)$. If $k_1=l_1$ and $k_2=l_2$, then $\fu
\subset \fsvect(V_{1}^*\otimes V_2)$.

\sssbegin[Maximal subalgebras in $\fsvect(V_{1}^*\otimes V_2)$]
{Lemma}[Maximal subalgebras in $\fsvect(V_{1}^*\otimes V_2)$]\label{L3.4.1} \textup{1)} If either $k_1\neq
l_1$ or $k_2\neq l_2$, then $\fsl(V_{1}\oplus V_2) \subset
\fsvect(V_{1}^*\otimes V_2)$ is a~maximal subalgebra.

\textup{2)} If $k_1=l_1$ and $k_2=l_2$, then
$\fp\fgl(V_{1}\oplus V_2) \subset \fsvect(V_{1}^*\otimes V_2)$ is a~maximal
subalgebra. \end{Lemma}

\begin{proof} It is subject to a~straightforward verification that
$\fu= (\fg_{-1}, \fu_{0})*$ and $\fu_{1}$ satisfy condition 2) of
Lemma~\ref{L3.2.1}. Condition l) is satisfied thanks to \cite{ShM}.
\end{proof}

\ssec{
An embedding $\fp\fq(\cA\otimes(V_{1}\oplus V_2)) \subset
\fsvect(V_{1}^*\otimes C(\cA)\otimes V_2)$}\label{cA} Let 
\[
\cA:=\Q(1)=\left\{\begin{pmatrix} a&b\cr b& a\end{pmatrix}\mid a, b\in \Cee\right\}~\subset \Mat(1|1) =
\Mat(V_0).
\]
Let $U = U_\ev$ and $V = V_\ev$ be even linear spaces of dimension
$k$ and $l$, respectively. Let $\fq(U \otimes \cA)$ be the
supercentralizer of the operator $J_{1} = 1\otimes \begin{pmatrix}
0&1\cr -1& 0\end{pmatrix}$. Clearly,
\[
\fq(U \otimes
\cA)\simeq (\Mat(U) \otimes \cA)^{SL}.
\]
Similarly, let $\fq(\cA\otimes V)$ be the Lie superalgebra~consisting of
operators supercommuting with $J_{2} =\begin{pmatrix} 0&1\cr -1&
0\end{pmatrix}\otimes 1$. Clearly,
\[
\fq(C(\cA) \otimes V) \simeq (\cA
\otimes \Mat(U))^{SL}.
\]

As proved in \cite{ShM}, the algebra~$\fq(U \otimes \cA)\odot
\fq(\cA\otimes V)$ is maximal in $\fsl(U \otimes \cA~\otimes V)$.

Now, on Lie superalgebra~$\fa = \fq(\cA~\otimes (U\oplus V))$,
consider the $\Zee$-grading induced by the following $\Zee$-grading of
$\cA~\otimes (U\oplus V)$:
\[
\deg \cA\otimes U = 1, \quad \deg\cA\otimes V = 0.
\]
Then, 
\[
\fa_{-1} = U^*\otimes \cA\otimes V, \quad \fa_{0} = \fq(U^*\otimes
C(\cA)) \oplus \fq(C(\cA)\otimes V), \quad \fa_{1} = U\otimes \cA\otimes V^*,
\]
and we get an homomorphism $\varphi : \fa\tto \fsvect(U^*\otimes
C(\cA)\otimes V)$. Set $\fu= \varphi(\fa)$. Clearly, $\fu=
\fp\fq(\cA\otimes (U\oplus V))$ with
\[
\fu_{-1} =\fa_{-1}, \quad\fu_{0} = \fq(U^*\otimes
C(\cA)) \odot \fq(C(\cA)\otimes V), \quad \fu_{1} =\fa_{1}.
\]

\sssbegin[A maximal subalgebra~in $\fsvect(U^*\otimes
C(\cA)\otimes V)$]
{Lemma}[A maximal subalgebra~in $\fsvect(U^*\otimes
C(\cA)\otimes V)$]\label{L3.5.1} The subsuperalgebra~$\fu=
\fp\fq(\cA\otimes (U\oplus V))$ is maximal in $\fsvect(U^*\otimes
C(\cA)\otimes V)$.
\end{Lemma}

\begin{proof} As well as Lemma~\ref{L3.4.1}, this Lemma
follows from the fact that $\fu = (\fg_{-1}, \fu_{0})_*$, Lemma
\ref{L3.2.1} and \cite{ShM}.\end{proof}

\subsection{Several maximal simple infinite-dimensional subalgebras} Here we consider embeddings
\[
\begin{array}{l}
\fvect(V\oplus W) \subset
\fsvect(V\otimes\Lambda(W))\text{ for $W$ odd,}\\
\fk(1|2n; n)\subset
\fsvect(\Lambda(n));\\
\fm(n; n) \subset \fsvect(\Pi(\Lambda(n))).
\end{array}
\]

Let $\fa= \fvect(m|n + k)$ and let $\fu= \fvect(m|n + k; k)$ be its
nonstandard regrading determined by condition $\deg~\xi_{1} = \ldots
= \deg~\xi_{k} = 0$. Let us represent $\fa_{-1}$ in the form
$\fa_{-1}=V\oplus \Span(\partial_{\xi_{1}}, \ldots ,
\partial_{\xi_{k}})$. Then, 
\[
\renewcommand{\arraystretch}{1.4}
\begin{array}{l}
\fu_{-1}= V \otimes \Lambda(\xi), \\ \fu_{0} = \fgl(V) \otimes
\Lambda(\xi)\inplus \fvect(\xi), \\ \fu_{1} = (\fvect(V)_1\otimes
\Lambda(\xi))\oplus (V^* \otimes \Span(\partial_{\xi_{1}}, \ldots
,
\partial_{\xi_{k}})).
\end{array}
\]

It is subject to a~direct verification that if $m + n\eps\neq 1$or $\eps$, then
\[
\fu=
\fvect(m|n + k; k) = (\fg_{-1},\fu_{0})_*.
\]
Moreover, in all cases, except the ones where
at the same time $k = l$ and $m \neq n$, we have an embedding $\fu_{0} \subset
\fsl(V)\otimes \Lambda(\xi)$, and therefore $\fu \subset \fsvect(V\otimes
\Lambda
(\xi))$.

In the exceptional case,
\[
\begin{array}{l}
\fu= \fvect(m|n + 1; 1) \subset \fvect(m + n|m + n);\\
\fu\cap\fsvect(m|n + 1; 1) = \fvect(m + n|m + n)\otimes \Lambda
(\xi)\inplus \Cee\partial_{\xi}.
\end{array}
\]

Now, let $\fa = \fm(n)$ be generated by functions of $\tau$, $q$, and
$\xi$; let $\fu = \fm(n; n)$. Then, 
\[
\fu_{-1}= \Pi(\Lambda(\xi)), \quad \fu_{0} = \Lambda(\xi)\inplus 
\fvect(\xi)\subset\fsl(\Pi(\Lambda(\xi))), \quad \text{ and }\fu=(\fg_{-1},
\fu_{0})_*.
\]
This determines an embedding $\fm(n; n) \subset \fsvect(\Pi(\Lambda
(\xi)))$.

Finally, let $\fa = \fk(1|2n)$ be generated by functions of $t$,
$\xi$, and $\eta$; let $\fu = \fk(1|2n; n)$.
Then, 
\[
\fu_{-1}= \Lambda(\xi), \quad \fu_{0} = \Lambda(\xi)\inplus 
\fvect(\xi)\subset\fsl(\Lambda(\xi)), \quad \text{ and }\fu=(\fg_{-1},
\fu_{0})_*.
\]
This determines an embedding $\fk(1|2n; n)\subset \fsvect(\Lambda
(\xi))$.

\sssbegin[Several maximal simple infinite-dimensional subalgebras]
{Lemma}[Several maximal simple infinite-dimensional subalgebras]\label{LinfDim} Let $m+n\eps\neq 1, 1+\eps$ or~$\eps$. Then, 

\textup{1)} $\fvect(m|n + 1; 1)$ is maximal in $\fvect(m + n|m + n)$;\\
\phantom{.}\hskip 8 mm $\fvect(m + n|m + n)\otimes\Lambda(\xi)\subset \Cee\partial_{\xi}$
is maximal in $\fs\fvect(m+n|m+n)$ if $m\neq n$;

\textup{2)} $\fvect(m|n + k; k)$ is maximal in
$s\fvect((m+n\eps)2^{k-1})$ unless $k = 1$ and $m \neq n$;

\textup{3)} $ \fk(1|2n; n)$ is maximal in $\fs\fvect(\Lambda(n))$;

\textup{4)} $\fm(n; n)$ is maximal in
$\fs\fvect(\Pi(\Lambda(n)))$.\end{Lemma}

\begin{proof} It is similar to that of Lemma~\ref{L3.4.1}. \end{proof}

\ssec{
Cartan prolongation of $(\Lambda(n), \fhei(2n) \inplus 
\fo(2n))_*$} For $n \geq 2$ but $n \neq 3$, the Lie superalgebra
$\fg_{0} = \fsl(\Lambda(n))$ contains, as is shown in \cite{ShM},
one more maximal subsuperalgebra~of type 2C, namely, $\fhei(2n)
\inplus \fo(2n)$.

Let $\fg_{-1} = \fu_{-1} =\Lambda(n) = \Lambda(\xi_1, \ldots,
\xi_n)$. We endow with a~tilde the function considered as an element
of $\fhei$ to avoid confusion with the element of~$\fu_{0}$. Then, 
\[
\renewcommand{\arraystretch}{1.4} \begin{array}{ll} \fu_0(2n) &=\fhei(2n)
\inplus \fo(2n)\\
&=\Span(\tilde 1; \tilde \xi_1, \ldots, \tilde \xi_n;
\partial_{\xi_1}, \ldots,
\partial_{\xi_n})\inplus \Span(\xi_i\partial_{\xi_j}; \xi_i\xi_j;
\partial_{\xi_i}\partial_{\xi_j} \;\text{ for all }i, j). \end{array}
\]
Set $(\fu_{-1}, \fu_{0})_*:=\mathop{\oplus}\limits_{i\geq -1}\fu_{i}$ .

\sssbegin[The Cartan prolong $(\Lambda(n), \fhei(2n) \inplus 
\fo(2n))_*$ is
trivial]
{Lemma}[The Cartan prolong $(\Lambda(n), \fhei(2n) \inplus 
\fo(2n))_*$ is
trivial]\label{L3.7.1} The Cartan prolong $(\fu_{-1}, \fu_{0})_*$ is
trivial for $n \geq 4$.\end{Lemma}

\begin{proof} Let~$F\in \fu_1$.
Then, for $1\in \fu_{-1}$ and an
arbitrary $\varphi\in \fu_{-1}$, we have
\be\label{(1.)}
[F, 1]\varphi = [F, \varphi]\ 1. 
\ee
Since $[F, \varphi]\in \fu_{0}$, it follows that
\be\label{(2.)}
\deg[F, \varphi]1\leq 2, 
\ee
i.e., $[F, 1]$ sends any function in $\fu_{-1}$ to a~function of~$\deg \leq
2$.

Let~$n > 4$. Then, $\fg_{0}$ does not contain nonzero operators satisfying conditions
\eqref{(1.)} and \eqref {(2.)}. Hence, in this case,
\[
[F, 1] = 0\text{ for every $F\in \fu_{1}$. }
\]
But
thanks to the Jacobi identity we have
\[
[F, \xi_i] = [F, [\tilde \xi_i, 1]] = [[F, \tilde \xi_i], 1] +(-1)^{p(F)}[\tilde \xi_i, [F, 1]] = 0.
\]
By raising the degree of $\varphi$ in a~similar way, we conclude that
\[
[F, \varphi]
= 0\text{ for every $F\in \fu_{1}$ and every $\varphi\in \Lambda(n) = \fu_{-1}$, }
\]
i.e., $[F, \fu_{-1}] = 0$; hence~$F =
0$.

Now, let~$n = 4$. The conditions \eqref{(1.)} and \eqref {(2.)} mean that $[F, 1]$ must be of
the form
\[
[F, 1] = \sum\alpha_{ij}\partial_{\xi_i}\partial_{\xi_j},
\;\text{where $\alpha_{ij}\in \Cee$ for every $F\in \fu_{1}$. }
\]
The subalgebra~of $\fu_{0}$ generated by the $\xi_k\partial_{\xi_l}$
and isomorphic to $\fgl(4)$ acts on~$\fu_{1}$. We have
\[
{}\left[\left[\xi_k\partial_{\xi_l}, F\right], 1\right] =
\left[\xi_k\partial_{\xi_l}, [F, 1]\right] + \left[F,
\left[\xi_k\partial_{\xi_l}, 1\right]\right] =
\left[\xi_k\partial_{\xi_l},
\sum\alpha_{ij}\partial_{\xi_i}\partial_{\xi_j}\right].
\]
Therefore, if $\fu_{1}\neq 0$, then there exists $F_0\in \fu_{1}$ such that
$[F_0, 1] =
\partial_{\xi_1}\partial_{\xi_2}$. But then
\[
[F_0, \xi_1] = [F_0, [\xi_1,1]] = [[F_0, \xi_1], 1] + [\xi_1, [F_0, 1]],
\]
i.e., $[F_0, \xi_1] = \partial_{\xi_2} +
\sum\alpha_{ij}\partial_{\xi_i}\partial_{\xi_j}$. Thus,
\[
[F_0, \xi_2\xi_3\xi_4]\xi_1 = -[F_0, \xi_1](\xi_2\xi_3\xi_4) = -\xi_3\xi_4 +
\varphi(\xi),
\]
where $\deg \varphi(\xi)\leq 1$. But $\fu_{0}$ does not contain such operators,
hence~$\fu_{1} = 0$.
\end{proof}

Now, let $n = 2$ and $\fg= \fs\fvect(2|2)$. Set $\fg_{-1} = \Pi(\Lambda(2))$
and realize
the subalgebra
\[
\fu_{0} = \fhei(4)\inplus \fo(4)\simeq\fhei(4)\inplus (\fsl(2)
\oplus\fsl(2))
\]
 in
the form
\[
\renewcommand{\arraystretch}{1.4}
\begin{array}{ll}
\fu_{0} &= \left(\Span(\tilde 1, \tilde \xi_1, \tilde \xi_2)\inplus 
\Span\left(\partial_{\xi_1}, \partial_{\xi_2},
\xi_1\partial_{\xi_2}, \xi_2\partial_{\xi_1},
\xi_1\partial_{\xi_1}- \xi_2\partial_{\xi_2}\right)\right)\\
& \inplus \Span(\xi_1 \xi_2, \ \partial_{\xi_1}\partial_{\xi_2},
\xi_1\partial_{\xi_1}+\xi_2\partial_{\xi_2}-1) \\
&=\left(\left(\Lambda(\xi_1, \xi_2)/\Cee \xi_1 \xi_2\right)\inplus 
\fs\fvect(\xi)\right)\inplus \fsl(2).
\end{array}
\]
The term $\left(\Lambda(\xi_1, \xi_2)/\Cee \xi_1 \xi_2\right)\inplus 
\fs\fvect(\xi)$ has a~nontrivial Cartan prolong, $\fb_{a,a}(2; 2)$. The
term $\fsl(2)$ acts on the
odd part of $\fg_{-1}$ via~the tautological representation, and the Cartan prolong
\[
(\fg_{-1}, \fu_{0})_*\simeq \fb_{a,a}(2; 2) \inplus \fsl(2) \subset
\fs\fvect(2|2)
\]
is a~maximal subalgebra.

Summing this section up we arrive at the following 2 theorems.

\ssbegin[Maximal $\fu\subset \fg$
such that $\fu_{-1} = \fg_{-1}$ and $\fu_{0}\subset\fg_{0}$ is maximal
irreducible]
{Theorem}[Maximal $\fu\subset \fg$
such that $\fu_{-1} = \fg_{-1}$ and $\fu_{0}\subset\fg_{0}$ is maximal
irreducible]\label{Th3.8} Maximal 
$\fu\subset \fg$
such that $\fu_{-1} = \fg_{-1}$ and $\fu_{0}$ is a~maximal
irreducible ~in $\fg_{0}$, neither almost simple
nor a~central extension of an almost simple Lie superalgebra~are
exhausted by the following list:

\noindent \underline{for $\fg$ of type $\fvect$}:

\textup{1${}_v$)} $\fu = \fsl(V_1\oplus V_2) \subset \fvect(V^*_1 \otimes V_2)$
if either $\dim V_1 \neq k
+ k\eps$ or $\dim V_2 \neq l + l\eps$ for some $k, l$;

\textup{2${}_v$)} $\fu= \fvect(m|n + 1; 1) \subset \fvect(m + n |m + n)$ if $m
\neq n$ and $m+
n\eps \neq 1, \eps$;

\noindent \underline{for $\fg$ of type $\fsvect$}:

\textup{1${}_s$)} $\fu =\fpgl(V_1 \oplus V_2) \subset \fsvect(V^*_1 \otimes V_2)$
if $\dim V_1 = k(1 + \eps)$ and $\dim V_2 = l(1 + \eps)$ for
some $k, l$;

\textup{2${}_s$)} $\fu= \fp\fq(\cA~\otimes (V_1 \oplus V_2)) \subset
\fs\fvect(V^*_1 \otimes C(\cA) \otimes V_2)$;

\textup{3${}_s$)} $\fu= \fvect(m|n) \otimes \Lambda(\xi_1)
\inplus \Cee\partial_{\xi_1} \subset \fs\fvect(m + n |m + n)$ if $m
\neq n$ and $m + n\eps\neq 1, \eps$;

\textup{4${}_s$)} $\fu= \fvect(m|n + k; k) \subset \fs\fvect((m +
n\eps)2^{k-1})$ if $m \neq n$ and $m + n\eps \neq 1, \eps, 1 + \eps$ and
either $k >1$ or $k = 1$ and $m = n$;

\textup{5${}_s$)} $\fu= \fk(1|2n; n) \subset \fs\fvect(\Lambda(n))$;

\textup{6${}_s$)} $\fu = \fm(n;n) \subset \fs\fvect(\Pi(\Lambda(n)))$;

\textup{7${}_s$)} $\fu = \fb_{a,a}(2; 2) \inplus \fsl(2) \subset
\fs\fvect(2|2)$.\end{Theorem}

\ssbegin[Maximal $\fu \subset \fg$ such that both
inclusions $0 \subset \fu_{-1} \subset \fg_{-1}$ are strict]
{Theorem}[Maximal $\fu \subset \fg$ such that both
inclusions $0 \subset \fu_{-1} \subset \fg_{-1}$ are strict]\label{Th3.9} Let $\fg= \fvect(m|n)$ or
$\fs\fvect(m|n)$ and let $\fu \subset \fg$ be such that both
inclusions $0 \subset \fu_{-1} \subset \fg_{-1}$ are strict. Then,
$\fu$ is maximal if and only if
\be\label{u}
\fu= \begin{cases}\Ann(\fu_{-1})\fvect(m|n) \inplus \fvect(\fu_{-1})&\text{if
$\fg= \fvect(m|n)$}\\
\Ann(\fu_{-1})\fs\fvect(m|n) \inplus \varphi_s(\fvect(\fu_{-1}))&\text{if $\fg=
\fsvect(m|n)$}
\end{cases}
\ee
for any subsuperspace $\fu_{-1} \subset \fg_{-1}$.
\end{Theorem}

\begin{proof} \underline{Necessity}. Let $\fu\subset\fg$ be maximal and such that
both inclusions $0 \subset \fu_{-1} \subset \fg_{-1}$ are strict. By the very
definition of the prolong $\fu^\fg_*$ of $\fu_{-1}$ in $\fg$, see \S\ref{2.6},
we have $\fu\subset \fu^\fg_*$; by item 1\textup{)} of Lemma \ref{L1.1},
the space $\fu^\fg_*$ is a~Lie subsuperalgebra~of $\fg$, and $\fu^\fg_*\neq \fg$
because $(\fu^\fg_*)_{-1} = \fu_{-1}\subsetneq \fg_{-1}$. The maximality of
$\fu$ implies therefore that $\fu= \fu^\fg_*$, and it only remains to compute~$\fu^\fg_*$.

For $\fg= \fvect(m|n)$, this is Lemma \ref{sss1.2.1}, which yields the first
line of \eqref{u}. For $\fg= \fsvect(m|n)$, item 2\textup{)} of Lemma
\ref{L1.1} gives $\fu^\fg_* = \fu^{\fvect(m|n)}_*\cap \fsvect(m|n)$. Now, an
element $D_1+D_2$, where $D_1\in \Ann(x)\fvect(x, y)$ and $D_2 \in\fvect(x)$,
is divergence-free if and only if $\Div D_1 = -\Div D_2$; since
\[
D_1+D_2 = \left(D_1 + (\Div D_2)\nfrac{1}{c_0}X_0\right) + \varphi_{s}(D_2),
\text{~~where $X_0$ and $c_0$ are as in \eqref{varphi_sbis},}
\]
where the first summand lies in $\Ann(x)\fsvect(m|n): =\Ann(x)\fvect(x,
y)\cap\fsvect(m|n)$ and the second one in $\varphi_{s}(\fvect(x))$, we get the
second line of \eqref{u}.

\underline{Sufficiency}. Let $\fu$ be as indicated in \eqref{u}, and $\fr \subset \fg$ be a~subalgebra
strictly containing~$\fu$. By item 2\textup{)} of Lemma \ref{L1.1}, we have $\fu_{-1} \neq
\fr_{-1}$: otherwise $\fr\subset \fu^\fg_* = \fu$. But
\[
\fu_{0}=\begin{cases}\fst(\fu_{-1})&\text{for $\fg= \fvect(m|n)$}\\
\fst(\fu_{-1})\cap\fsl(m|n)&\text{for $\fg= \fsvect(m|n)$}
\end{cases}
\]
and $\fu_{-1}$ is the only $\fu_{0}$-invariant subsuperspace of
$\fg_{-1}$ distinct from $0$ and~$\fg_{-1}$. Hence, $\fr_{-1} = \fg_{-1}$. Then,
$\fr \subset (\fg_{-1}, \fr_0)_*$ because $\fr$, being a~subalgebra~of the
transitive Lie superalgebra~$\fg$, is transitive, whereas the prolong is the
maximal transitive Lie superalgebra~with the given non-positive part, see
\S\ref{2.6}. But $\fu_{0} = \fst(\fu_{-1})\cap\fg_0$ is a~maximal
subsuperalgebra~in $\fg_{0}$, see \cite{ShM, ZZO} and \cite{LOdd}.

Hence, either $\fr_0 = \fu_{0}$ or $\fr_0 = \fg_{0}$. In the first
case, \eqref{autVcoord} implies that
\[
\fr \subset (\fg_{-1}, \fu_{0})_* =\faut(\fu_{-1})= \cF(y) \otimes\fvect(x) \inplus \fvect(y),
\]
whereas $a= y_1x_1\partial_{y_1}\in \fu$ does not lie in $(\fg_{-1},
\fu_{0})_*$ since the coefficient of $\partial_{y_j}$ in an element of
$\faut(\fu_{-1})$ depends on the $y$ only. This contradicts $\fr\supset\fu$.

If $\fg=
\fs\fvect(m|n)$, the element $a$ should in addition be divergence-free, and any
$D\in\fvect(x)$ with a~\textbf{non-constant} divergence provides one: indeed,
$a: = \varphi_{s}(D) = D - (\Div D)\nfrac{1}{c_0}X_0$ lies in $\fu$ by
\eqref{varphi_sbis}, whereas the coefficient of $\partial_{y_j}$ in $a$ is
proportional to $\Div D\not\in\Cee$, and hence $a\not\in\faut(\fu_{-1})$.
Explicitly, taking $X_0=y_1\partial_{y_1}$, we can take:
\[
a=\begin{cases}
y_1x_1\partial_{y_2}&\text{if $l\geq 2$};\\
\pm\varphi_{s}(x_1x_2\partial_{x_2}) = y_1x_1\partial_{y_1}\mp x_1x_2\partial_{x_2}&\text{if $k\geq 2$};\\
\pm\varphi_{s}(\nfrac12x_1^2\partial_{x_1}) = y_1x_1\partial_{y_1}\mp\nfrac12x_1^2\partial_{x_1}&\text{if $k=1$,
$p(x_1)=\ev$},
\end{cases}
\]
where the signs are determined by the parities of $x_1$, $x_2$ and~$y_1$.
The only case where no such $a$ exists, i.e., where the inclusion $\fu^\fg_*
\subset (\fg_{-1}, \fu_{0})_*$ does hold, is the case where $k=1$ with
$p(x_1)=\od$ and $l=1$: then $\fvect(x) = \fvect(0|1)$, and every element of it
is of constant divergence. In this case, $\sdim \fg_{-1} = 1|1$ and $\dim
\fu_{-1} =\eps$, so $\fg= \fs\fvect(1|1)$ or $\fs\fvect(0|2)$, and neither of
these should interest
us since neither is simple.

If $\fr_{0} = \fg_{0}$, then $\fr\supset \fg_{-1}\oplus\fg_{0}$, and hence
Lemma \ref{L3.1.1} applies. If $\fg= \fsvect(m|n)$, then, by item
1\textup{)} of that Lemma, either $\fr= \fg$ or $\fr = \fg_{-1}\oplus\fg_{0}$;
the latter is impossible since $\fu_1\ni a\neq 0$. If $\fg= \fvect(m|n)$, then,
by item 2\textup{)} of that Lemma, either $\fr = \fg$ or $\fr$ is contained
in $\fd(\fsvect(m|n))$ or in~$\fpgl(m+1|n)$. Neither of the latter two contains
$\fu$: the element $b: = y_1x_1\partial_{x_1}\in \Ann(x)\fvect(x, y)\subset\fu$
is of non-constant divergence, $\Div b = \pm y_1$, whereas the divergences of
the elements of $\fd(\fsvect(m|n)) = \fsvect(m|n)\inplus\Cee E$ are constants;
and $b\not\in\fpgl(m+1|n)$ because $\fpgl(m+1|n)_1 = \Span(z_iE\mid z=(x,
y))$, where $E =\sum_j z_j\partial_{z_j}$ is the Euler operator, and in
$\sum_ic_iz_iE$ the coefficient of $\partial_{y_1}$ is $(\sum_ic_iz_i)y_1$,
which vanishes only if all~$c_i=0$. Hence, $\fr =\fg$, and $\fu$ is maximal.
\end{proof}

Note that the subsuperalgebras $\fu$ of the type described
are
non-transitive (i.e., $\dim \fu_-$ is smaller than the number of indeterminates
involved in the
description of~$\fu$). The transitive part of $\fu$ is 
\be\label{t(u)}
\text{either $\varphi_s(\fvect(\fu_{-1}))$ or the subalgebra~
$\ft(\fu): =\fvect(\fu_{-1})$.} 
\ee

The ideal $\fii = \mathop{\oplus}\limits_{k\geq 0}\fii_k :=
\Ann(\fu_{-1})\fvect(m|n)$ (resp. $\Ann(\fu_{-1})\fs\fvect(m|n)$) is
singled out by the condition \eqref{3.1}:
\begin{equation}\label{3.1}\fii_k
= \{g\in \fu_{k}\mid (\ad_{\fu_{-1}})^{k+1}g = 0\}.
\end{equation}
Let now $\fg$ be the Lie superalgebra
preserving the 2-form~$\omega$. We realize $\fg$ by generating
functions. Let $p = (p_1, \ldots , p_k)$ be the Darboux
coordinates in (the supermanifold associated with the superspace)
$\Ker(\omega|_{\fu_{-1}})$ and let $\zeta = (\zeta_1, \ldots ,
\zeta_s)$ be the complementary coordinates in $\fu_{-1}$; let
$\tilde \zeta = (\tilde \zeta_1, \ldots , \tilde \zeta_t)$ be the
complementary coordinates in $\fu_{-1}^\perp$ and let $q = (q_1,
\ldots , q_k)$ be the duals to the $p$'s. Any of these
indeterminates $p_i, \zeta_j, \tilde \zeta_l$ may be of any parity;
the parity of the $q$'s is identical to that of the respective
$p$'s.

Thus, we consider $\fg$ as the Lie superalgebra~with the Poisson bracket \textup{(}if $\fg=\fh(2m|n)$ or $\fh'(0|n)$\textup{)} or with the
antibracket \textup{(}if $\fg=\fle(n)$ or $\fsle'(n)$\textup{)}, realized on the space
of ``functions" $\cF(p, q, \zeta, \tilde
\zeta)/\langle\const\rangle$. Let \[ \Ann(\fu_{-1}) = \{X\in \fg_{0}\mid [\fu_{-1},
X] = 0\} \subset \fg_{0}\] be the subsuperalgebra~of $\fg_{0}$
commuting with~$\fu_{-1}$. In our notation,
\[
\Ann(\fu_{-1}) = \Span(p^2, \; p\tilde \zeta,\; \tilde \zeta^ 2).
\] It
is clear that the set
\[
\fii = \mathop{\oplus}\limits_{k\geq 0}\fii_k := \Ann(\fu_{-1})\cF(p,
q, \zeta, \tilde \zeta)
\]
is a~ subsuperalgebra~of $\fg$ singled
out by conditions \eqref{3.1}. As we have already mentioned, the set
$\Span(p) \otimes\cF(q)$ is a~ subsuperalgebra~isomorphic to~$\fvect(q)$. Therefore, the semidirect sum
\be\label{t(u)omega}
\ft(\fu) = \cF(\zeta)\otimes\fvect(q) \inplus\cF(\zeta)/\Cee\cdot 1
\ee
is well-defined and $\fii$ is a~$\ft(\fu)$-invariant subsuperalgebra.

\ssbegin[A technical Lemma]
{Lemma}[A technical Lemma]\label{L3.10} We have
\[
\fu^\fg_* =\Ann(\fu_{-1})\cF(p, q,
\zeta, \tilde
\zeta)\inplus \left(\cF(\zeta)\otimes\fvect(q)\inplus \cF(\zeta)/\Cee\
1\right).
\]
\end{Lemma}

\begin{proof} Straightforward. \end{proof}

\sssbegin[Transitive parts]
{Remarks}[Transitive parts]\label{R3.10.1} 1) If $\fg = \fh(2m|n)$,
then $\cF(p, q, \zeta, \tilde
\zeta) \simeq \fpo(2m|n)$ and
$\cF(\zeta)/\Cee\cdot 1 \simeq\fh(\zeta)$.

If $\fg= \fle(2m|n)$, then $F(p, q, \zeta, \tilde \zeta)\simeq \fb(n)$
and
$\cF(\zeta)/\Cee\cdot 1 \simeq\fle(\zeta)$.

2) The subsuperalgebra~$\ft(\fu) = \cF(\zeta) \otimes\fvect(q)\inplus 
\cF(\zeta)/\Cee\cdot 1$
is the transitive part of~ $\fu$.

3) In particular cases we have:
\[
\renewcommand{\arraystretch}{1.4}
\begin{array}{ll}
\ft(\fu)=\begin{cases}
\fh(\fu_{-1})&\text{if $\fg= \fh(2m|n)$, }\\
\fle(\fu_{-1})&\text{if $\fg= \fle(n)$},\\
\fsle(\fu_{-1})&\text{if $\fg= \fsle'(n)$},
\end{cases}&
\text{if $\omega|_{\fu_{-1}}$ is non-degenerate;}\\
\renewcommand{\arraystretch}{1.4}
\ft(\fu)=\begin{cases}
\fvect(q)&\text{if $\fg= \fh(2m|n)$ or $\fle(n)$, }\\
\fsvect(q)&\text{if $\fg= \fsle'(n)$ }
\end{cases}&\text{if $\omega|_{\fu_{-1}}\equiv 0$.}
\end{array}
\]
\end{Remarks}

\ssbegin[Maximal $\fu \subset \fg$ such that
both inclusions $0 \subset \fu_{-1} \subset \fg_{-1}$ are strict]
{Theorem}[Maximal $\fu \subset \fg$ such that
both inclusions $0 \subset \fu_{-1} \subset \fg_{-1}$ are strict]\label{ThMax} Let $\fg = \fh(2m|n)$ for $m > 0$, or $\fh'(0|n)$ for
$n > 4$, or $\fg
= \fle(n)$ or $\fsle'(n)$ for~$n > 1$.

Let $\fu \subset \fg$ be such that
both inclusions $0 \subset \fu_{-1} \subset \fg_{-1}$ are strict. Then, $\fu$
is maximal if and
only if $\fu = V^\fg_*$ for a~nonsingular subsuperspace $V \subset \fg_{-1}$.

Explicitly, $V^\fg_*$ is isomorphic to:
\[
\renewcommand{\arraystretch}{1.4}
\begin{array}{ll}
\Ann(V) \fpo(2m|n) \inplus &\\
\left(\cF(V/(V\cap V^{\perp})) \otimes \fvect(\fg_{-1}/(V + V^{\perp}) \inplus 
\fh(V/(V\cap V^{\perp}))\right)&\text{if $\fg= \fh(2m|n)$};\\
\Ann(V) \fb(n)\inplus &\\
\left(\cF(V/(V\cap V^{\perp})) \otimes \fvect(\fg_{-1}/(V + V^{\perp}) \inplus 
\fle(V/(V\cap V^{\perp}))\right)&\text{if $\fg= \fle(n)$};\\
\left(\Ann(V) \fb(n) \cap\fsle(n)\right) \inplus &\\
\left(\cF(V/(V\cap V^{\perp})) \otimes \fsvect(\fg_{-1}/(V + V^{\perp}) \inplus
\varphi_s(\fle(V/(V\cap V^{\perp})))\right)&\text{if $\fg= \fsle(n)$}.
\end{array}
\]
\end{Theorem}

\begin{proof} The formulas for the explicit form of $V^\fg_*$ follow from Lemma
\ref{L3.10} and Remarks \ref{R3.10.1} to it. The necessity of the form $V^\fg_*$ of
$\fu$ follows from Lemma \ref{L1.1}: any $\fu$ with the given $\fu_{-1} = V$ is
contained in $V^\fg_*$, which is a~subalgebra~distinct from $\fg$ because
$(V^\fg_*)_{-1} = V\subsetneq\fg_{-1}$, and hence $\fu = V^\fg_*$ if $\fu$ is
maximal, cf. the proof of Theorem \ref{Th3.9}. It suffices to verify that
$V^\fg_*$ is maximal if and only if $V$ is nonsingular.

Let $\fr= \mathop{\oplus}\limits_{i\geq -1}\fr_i$ be a~Lie
subsuperalgebra~of $\fg$ \textbf{strictly} containing $\fu = V^\fg_*$. Then,
$\fr_{-1}\supset\fu_{-1}$ is a~strict inclusion, again by item 2\textup{)}
of Lemma \ref{L1.1}: otherwise $\fr\subset V^\fg_* =\fu$. Two cases are
possible: either $\omega|_{V}$ is non-degenerate or
$\Ker(\omega|_{V}) \neq 0$. Let us consider them separately.

\underline{A) $\omega|_{V}$ is non-degenerate}. Then, $\fg_{-1}=V\oplus
V^\perp$ and $\fg_{-1}$ has
only two $\fu_{0}$-invariant subspaces distinct from $0$ and $\fg_{-1}$: $V$
and~$V^{\perp}$. Hence,
$\fr_{-1}= \fg_{-1}$. But, as we have already observed in
\S\ref{SAred}, see also \cite{ShM}, the subalgebra~$\fu_{0} =
\fst(V)\cap\fg_{0}$ --- which is $\fosp(V)\oplus\fosp(V^\perp)$ if $\fg
=\fh(2m|n)$ or $\fh'(0|n)$, and $\fpe(V)\oplus\fpe(V^\perp)$ if $\fg=\fle(n)$
--- is a~maximal subsuperalgebra~in~$\fg_{0}$. Hence, since $\fr_0\supset
\fu_0$, either
Aa) $\fr_0 = \fu_{0}$ or Ab) $\fr_0 = \fg_{0}$.

\underline{Aa) $\fr_0 = \fu_{0}$}. We have
\[
\fr\subset (\fg_{-1}, \fu_{0})_* = \cF(\zeta)/\Cee\cdot 1\oplus
\cF(\tilde\zeta)/\Cee\cdot 1\] (this is the second display of \S\ref{SAred},
where $V=\Span(\zeta)$ and $V^\perp = \Span(\tilde\zeta)$, so that
\[
\fh(V)\oplus\fh(V^\perp) = \cF(\zeta)/\Cee\cdot 1\oplus\cF(\tilde\zeta)/\Cee\cdot
1;
\] 
for $\fg= \fsle'(n)$, we take the corresponding
$(\fsle(V)\oplus\fsle(V^\perp))\inplus\Cee\cdot D$ instead) and, since
$\dim V^{\perp} /(V \cap V^{\perp}) = \dim \Span(\tilde\zeta) \neq
\eps$, we see that $\dim \Span(\tilde\zeta^ 2) \neq
 0$ and $\Span(\tilde\zeta^2)\otimes\cF(\zeta)$ is
contained in $V^\fg_*$ but not in $\fr$ contrary to our assumption.

\underline{Ab) $\fr_0 = \fg_{0}$}. Then, $\fr\supset\fg_{-1}\oplus\fg_{0}$, so
Lemma \ref{L3.1.1} applies. If $\fr \neq \fg$, then by item 1\textup{)} of
that Lemma
\[
\fr= \fg_{-1}\oplus \fg_{0}\;\text{ for $\fg =\fh(2m|n)$, $\fh'(0|n)$ for $n
\neq 6$, and $\fsle'(n)$.}
\]
The explicit form of $V^\fg_*$ shows that $\fu \not\subset \fr$, since
$\fu_1\supset\Ann(V)\cdot\Span(p, q, \zeta, \tilde\zeta)\neq 0$.

If $\fg = \fle(n)$, then $\fr= \fd(\fsle(n))$ and the subsuperalgebra
$\fle(V/(V \cap V^{\perp}) = \fle(\zeta) \subset V^\fg_*$ does not belong to~$\fr$.

Finally, if $\fg= \fh'(0|6)$, then $\fu_2 \neq 0$ for any $V$, and therefore
$\fu \not\subset \fr$, either.

\underline{B) $\Ker(\omega|_{V}) \neq 0$}. Then, $V+V^{\perp}$ does not coincide
with $\fg_{-1}$, and $ V+V^{\perp}\supset V$ is a~strict inclusion unless
$V^\perp\subset V$, i.e., unless there are no $\tilde\zeta$ at all. Moreover, if
$V$ is nonsingular, then $V + V^{\perp}$
is the
only $\fu_{0}$-invariant subspace distinct from $\fg_{-1}$ and containing~$V$.
This
means that either $\fr_{-1} = V + V^{\perp}$ (and then $V^\perp\not\subset V$)
or $\fr_{-1} = \fg_{-1}$.

If $\fr_{-1} = V + V^{\perp} = \Span(p, \zeta, \tilde\zeta)$, then $\fr \subset (V +
V^{\perp})^\fg_*$, whose only summand containing the $q$ is $\Ann(V+V^\perp)\cF(p, q,
\zeta, \tilde\zeta)= \Span(p^2)\cF(p, q,\zeta, \tilde\zeta)$ together with
$\cF(\zeta, \tilde\zeta)\otimes\Span(p)\otimes\cF(q)$. Comparing
explicit forms of $\fu$ and $\fr$ we see that the subspace $\Span(\tilde\zeta^
2) \otimes \Span(q)$
is contained in $\fu$ but not in $\fr$ in contradiction with our assumption
$\fr\supset
\fu$. Here, $\Span(\tilde\zeta^2)\neq 0$ precisely because $V$ is not singular in
the sense (a).

Let $\fr_{-1} = \fg_{-1}$. The subalgebra~$\fu_{0}$ is maximal in $\fg_{0}$
only if
$\omega|_{V}\equiv 0$. Otherwise, as is clear from the action of $\fu_{0}$ on
$\fg_{0}/\fu_{0}$, the Lie superalgebra~$\fu_{0}$ is contained in the only
maximal Lie
subsuperalgebra~of $\fg_{0}$ containing it, namely $\fst(V \cap
V^{\perp})\cap\fg_0$, for a~nonsingular subsuperspace $V
\subset
\fg_{-1}$. Therefore, either $\fr_0 \subset \fst(V \cap V^{\perp})$, and hence
$\fr \subset (\fg_{-1}, \fst(V \cap V^{\perp}))_*$, or
$\fr_0 = \fg_0$.

As follows from the second case of \S\ref{SAred},
\[
\fr \subset (\fg_{-1}, \fst(V \cap V^{\perp}))_* = \cF(p, \zeta, \tilde\zeta)/\Cee 1
\inplus \Span(q) \otimes\cF(p)
\]
and $\Span(\zeta)\otimes\Span(p)\otimes \Span(q)$, which lies in the transitive
part \eqref{t(u)omega} of $\fu$, does not lie in $\fr$: the coefficient of
$q_i$ in an element of $\fr$ is a~function of the $p$ alone, whereas the
coefficient of $q_1$ in $\zeta_1p_1q_1$ is~$\zeta_1p_1$. If,
instead, $\omega|_{V}\equiv 0$, then $\fu_{0}$ itself is maximal, so either
$\fr_0 = \fu_{0}$, and then $\fr\subset(\fg_{-1}, \fu_{0})_* = \cF(p,
\tilde\zeta)/\Cee 1\inplus\Span(q)\otimes\cF(p)$ by the same display of
\S\ref{SAred}, whereas $\fu$ contains $\Span(\tilde\zeta^2)\otimes\Span(q)$ ---
or, if there are no $\tilde\zeta$, i.e., if $V=V^\perp$ is Lagrangian,
$\Span(p^2)\otimes\Span(q^2)$, which is of degree 2 in the $q$, whereas the
elements of $\fr$ are of degree $\leq 1$ in the $q$ --- or~$\fr_0=\fg_0$.

Finally, the case where $\fr_0 = \fg_0$ is treated in a~way similar
to that where $\omega|_V$ is non-degenerate, see Ab).

Let now $V$ be singular. If $\dim V^{\perp}/(V \cap V^{\perp}) = \eps$ and $\fg
=
\fh(2m|n)$ or $\fh'(0|n)$, we have
\[
V^\fg_* \subset (V \cap V^{\perp})^\fg_*,
\]
so that a~singular $V$ of this kind yields no maximal subalgebra.

In other cases of singular subsuperspaces the component $\fg_{-1}$ contains
one (if $\fg = \fle(n)$ or $\fs\fle'(n)$) or two (if $\fg = \fh(2m|n)$ or
$\fh'(0|n)$)
$\fu_{0}$-invariant subspaces $\tilde V$ such that $V^\fg_*\subset
\tilde V^\fg_*$ and $V \subset \tilde V \subset V + V^{\perp}$.
 \end{proof}

\ssec{
Proof of Theorem \ref{Th2.1.2}}\label{s2.1.2proof} Clearly, $\fg_\ev=\fg_{-2}\oplus
\fg_{0}\oplus \fg_{2}$.

1) Consider the $\fg_{0}$ action in $\fg_{2}$ and recall that $\fg_\ev$ is
reductive. Therefore, $\fg_\ev$ contains a~subalgebra~$\fg_{-2}\oplus
[\fg_{-2}, \fg_{2}\oplus \fg_{2}$ isomorphic to $\fsl(2)$ as a~direct summand.
From the classification of almost simple finite-dimensional Lie superalgebras we
deduce the list of the Theorem.

It remains to check if these superalgebras can indeed be embedded into a~contact
superalgebra.

In all these cases, except $\fosp(3|2n)$, the action of the $\fsl(2)$ part of
$\fg_\ev$ in $\fg_\od$ is a~multiple of the tautological (2-dimensional)
representation, so the element $H\in\fsl(2)$ determines the contact grading of~$\fg$.

Let us realize $\fg=\fsl(m|2)$ and $\fosp(m|2)$ as $\fsl(V)$ and~$\fosp(V)$.
Introduce a~compatible $\Zee$-grading in $V$ by setting
\[
V_0=V_\ev,\quad V_\od=V_{-1}\oplus V_1 \;(\text{the eigenspaces of $H$}).
\]
This $\Zee$-grading can be naturally extended to a~$\Zee$-grading of $\fg$;
clearly,
\[
\fg_0=\begin{cases}\fgl(m)&\text{for $\fsl(m|2)$},\\
\fc(\fo(m))&\text{for $\fosp(m|2)$}.\end{cases}
\]
The geometry behind this grading is lucidly expressed in Arnold's
book \cite{A} (in nonsuper situation) and is related to the
contactization of symplectic structure: every cotangent bundle
$T^*M$ has a~ symplectic structure and its
projectivization $P(T^*)M$ has a~contact structure. The action
of the symplectic group extends to an action in $P(T^*)M$, hence we
get an embedding
$\fsp(2n)\subset \fk(2n-1)$.

In the super case, the above arguments imply the existence of embeddings
\[
\begin{array}{l}
\fosp(m|2n)\subset \fk(2n-1|m);\\
 \fspe_{a, b}(n+1)\subset \fm(n);\\
 \fsq(n+1)\subset \fk(2n+1|2n+2;1).\end{array}
\]
On the level of linear algebra, we deal with a~standard contact structure in
the projectivization $P(V)$ of a~superspace $V$ with a~fixed antisymmetric
bilinear form $\omega$ and with an extension of the $\faut(\omega)$-action from
$V$ to~$P(V)$. We consider the representation $\Pi(\id)$ of $\fosp(V)$: the
tautological one with the reversed parity.

For $\fg=\fosp(m|2n)$, $\fpe(n+1)$, and $\fspe_{a, b}(n+1)$, we immediately get
the action in $P(V)$, as explained in \cite{A}.

For $\fg=\fsl(n+1)$, we identify $V_1$ with $V_{-1}^*$ and consider a~normal
antisymmetric
bilinear form in $\Pi(V_{-1}\oplus V_{-1}^*)$.

For $\fg=\fosp(4|2n)=\fosp(V)$, we have
\[
\fg_\ev=\fo(4)\oplus\fsp(2n):=
\fo(V_\ev)\oplus\fsp(V_\od)\simeq\fsl(2)\oplus\fsl(2)\oplus\fsp(2n),
\]
and $V_\ev=U\otimes U'$, where $U$ and $U'$ are the two 2-dimensional tautological modules over the two $\fsl(2)$'s.
Let us represent $U$ as $U_{-1}\oplus U_1$ and introduce
 a~$\Zee$-grading in $V$ by setting
\[
V_0:=V_\od,\quad V_{-1}:=U_{-1}\otimes U', \quad V_{1}:=U_{1}\otimes U'.
\]
This $\Zee$-grading induces a~compatible contact $\Zee$-grading of $\fg$ where $\fg_0=\fgl(2)\oplus \fsp(2n)$. The geometry is similar to the above.

In the exceptional case of $\fosp(3|2n)$, we would like to use the fact
$\fo(3)\simeq \fsl(2)$. However, the action of $H$ (which determines the
grading) in $\fg_\od$ is the tautological representation of $\fo(3)$ and so has
three eigenvalues, not two. As a~result, we get an incompatible $\Zee$-grading
of depth~1 corresponding to the embedding
\[
\fosp(m+2|2n)\subset\fvect(m|2n);
\]
in particular, $\fosp(3|2n)\subset\fvect(1|2n)$.

2) For the $\fm$-contact gradings, the component $\fg_{-2}$ is odd and 1-dimensional, so the bracket
\[
\fg_{-1}\times \fg_{-1}\tto \fg_{-2}\simeq\Cee\cdot z, \text{~~where $p(z)=\od$},
\]
is a~non-degenerate \textbf{odd} symmetric bilinear form on~$\fg_{-1}$. Hence, $\sdim\fg_{-1}=k|k$
and, since $\fg_0$ preserves this form up to the factor by which it acts on $\fg_{-2}$, we have
$\fg_0\subset\fc(\fpe(\fg_{-1}))=\fc(\fpe(k))$. Therefore,
\[
\fg\subset(\fg_{-2}\oplus\fg_{-1}, \fc(\fpe(k)))_*=\fm(k),
\]
and the grading of $\fg$ is the one induced by the standard grading of $\fm(k)$; in particular,
$\fg_\ev$ acts on $\fg_{-1}=V\oplus\Pi(V^*)$, where $\sdim V=k|0$, preserving the pair of dual
summands, so that $\fg_{0, \ev}\subset\fgl(V)=\fgl(k)$. Among the almost simple Lie superalgebras,
only $\fpgl(m+1|n+1)$, with $k=m+n$, and $\fspe_{a, b}(n+1)$, with $k=n$, admit such
a~$\Zee$-grading, see \cite{LS1, ShM}, whence the list of item 2\textup{)}. The word ``only\rq\rq\
refers, here as in item 1\textup{)}, to the \textbf{standard} $\fm$-contact gradings.

\qed

\section{Hermitian symmetric spaces and superspaces}

The following tables borrowed from \cite{LPS} are excerpts of Serganova's
classification \cite{Se}, see also \cite{Zi}. We recall the notational
conventions of~\cite{Se, LPS} and definitions adopted there.

\underline{In Table \eqref{HermSSpaces}}: ${\mathfrak{s}} =
(\Lie(S_c))\otimes
{\mathbb C}$,

NCHSS abbreviates \textit{noncompact Hermitian symmetric
space},\index{HSS, Hermitian symmetric space}

in the diagram of ${\mathfrak{s}}$, the maximal parabolic subalgebra
${\mathfrak p} = \Lie(P)$, such that $X$ can be represented as
$(S_c)^{{\mathbb C}}/P$, is determined by one vertex: the last one
in cases 0, 2, 3, $E_6$, $E_7$, the first one in case 4, the $p$-th
one in case~1.

\underline{In Table \eqref{Table6}}: A~homogeneous space $G/P$, where
$G$ is a~simple Lie supergroup,~$P$ its parabolic subsupergroup
corresponding to several omitted generators of a~Borel subalgebra
(description of these generators can be found in~\cite{GL1}), is
said to be of \textit{depth}\index{Grading, depth of} $d$ and \textit{length}\index{Grading, length of} $l$ if such are
the depth and length of $\text{Lie} (G)$ in the ${\mathbb
Z}$-grading compatible with that
of $\text{Lie}(P)$.

Note that all superspaces of Tables \eqref{Table6} and \eqref{Table61} have
a~Hermitian structure (hence are of depth~1) except
$PeGr$ and $SPeGr$ (no
Hermitian structure), $PeQ$ (no Hermitian structure, length~$2$),
$CGr_{0, k}^{0, n}$ and $SCGr_{0, k}^{0, n}$ (no Hermitian
structure, length $n-k$ and $n-k-1$, respectively).

The dots stand for
the notions without name or those that require too much space.

The
signs ${\mathfrak a} \inplus {\mathfrak b}$\index{${\mathfrak a} \inplus {\mathfrak b}$} and ${\mathfrak a} \niplus {\mathfrak b}$\index{${\mathfrak a} \niplus {\mathfrak b}$}
de\-no\-te the semidirect sum of Lie superalgebras, the ideal being the left
summand in the first of them and the right summand in the second;

Lie superalgebras ${\mathfrak{osp}}_{\alpha}(4|2)$ and
${\mathfrak a}{\mathfrak b}(3)$ are described via~their Cartan
matrices, see~\cite{CCLL, BGL}; for their meaningful interpretation, see \cite{Eld, Eldi}.




\be\label{HermSSpaces}\footnotesize
\centering\begin{tabular}{|@{\,}c@{\,}|@{\,}c@{\,}|@{\,}c@{\,}|@{\,}c@{\,}|@{\,}
c@{\,}|@{\,}c@{\,}|@{\,}c@{\,}|}
\hline
&&&&&&\\
[-3mm] &\parbox[c]{18mm}{\centering Name\\ of CHSS $X$}
&$X=S_c/G_c$&${\mathfrak{s}}_0=({\mathfrak g}_c)^{{\mathbb C}}$&
${\mathfrak{s}}_{-1}\sim T_0X$&$ (S_c)^*$& \parbox[c]{18mm}{\centering
Name\\ of NCHSS}\\[2mm] \hline
0&$\vphantom{\Big|}{\mathbb C} P^n $ &$SU(n+1)/U(n) $ &$
{\mathfrak{gl}} (n) $& $\id$ &$ SU(1, n)$ &${}^*{\mathbb C}{\mathbb P}
^n$\\
\hline 1&$\vphantom{\Big|}Gr_p^{p+q} $ &$SU(p+q)/S(U(p)\times U(q))$
&$ {\mathfrak{gl}} (p)\oplus {\mathfrak{gl}}(q)$&
$\id\bigodot\id^* $ &$SU(p, q)$ &$ {}^*Gr_p^{p+q}$\\
\hline
2&$\vphantom{\Big|}OGr_n $ &$SO(2n)/U(n) $ &${\mathfrak{gl}}(n) $ &$
\Lambda^2\id $ &$SO(n, n)$ &$ {}^*OGr_n$\\
\hline
3&$\vphantom{\Big|}LGr_n $ &$Sp(2n)/U(n) $ &${\mathfrak{gl}}(n) $
&$S^2\id$&$Sp(2n; {\mathbb R}) $ &${}^*LGr_n$\\
\hline
4&$\vphantom{\Big|}Q_n $ &$SO(n+2)/(SO(2)\times
SO(n))$&${\mathfrak{co}}(n)$&$\id$&$SO(n, 2)$ &$ {}^*Q_n$\\
\hline
$E_6$&$\vphantom{\Big|}{\mathbb O} P^2$ &$E_6 /(SO(10)\times
U(1))$&${\mathfrak{co}}(10)$&$R(\pi_{5}) $&$E_6^* $&${}^*({\mathbb O} P^2)$\\
\hline
$E_7$&$\vphantom{\Big|}{\mathbb E}$&$E_7 /(E_6 \times U(1)) $
&${\mathfrak{ce}}_6 $ &$R(\pi_{1})$&$E_7 ^* $&${}^*{\mathbb E}$\\
\hline
\end{tabular}
\ee

\textbf{Occasional isomorphisms:} 
\[
\begin{array}{l}
\text{$Gr_p^{p+q} \cong Gr_q^{p+q}$, \ \ $Gr_1^{n+1} \cong {\mathbb C} P ^n$,}\\
\text{$Q_1 \cong {\mathbb C} P ^1$, \ \ $Q_2 \cong S^2\times S^2$, \ \ $Q_3 \cong
LGr_2$, \ \ $Q_4 \cong Gr_2^4$,}\\
\text{$OGr_2 \cong LGr_1 \cong {\mathbb C} P
^1$, \ \ $OGr_3 \cong Gr_3^4\cong {\mathbb C} P^3$. }
\end{array}
\]



\newpage
\be\label{Table6}
\footnotesize
\begin{tabular}{|@{\,}>{\centering\arraybackslash}p{18mm}@{\,}|@{\,}>{\centering\arraybackslash}p{30mm}@{\,}|@{\,}>{\centering\arraybackslash}p{18mm}@{\,}|@{\,}>{\raggedright\arraybackslash}p{49mm}@{\,}|@{\,}>{\centering\arraybackslash}p{19mm}@{\,}|@{\,}>{\centering\arraybackslash}p{15mm}@{\,}|}
\hline
$\fg$&$\fg_0$&$\fg_{-1}$&\multicolumn{1}{c|}{Interpretation}&Underlying domain&Name\\
\hline
$\fsl(m|n)$&$\mathfrak{s}(\fgl(p|q)\oplus{}$ $\fgl(m-p|n-q))$&$\id_1\boxtimes\id_2^*$&
Supergrassmannian of the $(p|q)$-dimensional subsuperspaces in $\Cee^{m|n}$&
$Gr_p^m\times Gr_q^n$&$Gr_{p,q}^{m,n}$\\
$\fpsl(m|m)$&$\mathfrak{ps}(\fgl(p|p)\oplus{}$ $\fgl(m-p|m-p))$&$\id_1\boxtimes\id_2^*$&
Same for $m=n$, $p=q$&$Gr_p^m\times Gr_p^m$&$Gr_{p,p}^{m,m}$\\
\hline
$\fosp(m|2n)$&$\mathfrak{cosp}(m-2|2n)$&$\id$&
Superquadric of $(1|0)$-dimensional lines in $\Cee^{m|2n}$ isotropic with respect to
a~non-degenerate even form&$Q_{m-2}$&$Q_{m-2,n}$\\
\hline
$\fosp(2m|2n)$&$\fgl(m|n)$&$\Lambda^2(\id)$&
Orthosymplectic Lagrangian supergrassmannian of the $(m|n)$-dimensional subsuperspaces
in $\Cee^{2m|2n}$ isotropic with respect to a~non-degenerate even form&
$OGr_m\times LGr_n$&$OLGr_{m,n}$\\
\hline
$\fsq(n)$&$\mathfrak{s}(\fq(p)\oplus\fq(n-p))$&$\id\bigodot\id^*$&
Queergrassmannian of the $\Pi$-symmetric $(p|p)$-dimensional subsuperspaces in $\Cee^{n|n}$&
$Gr_p^n$&$QGr_p^n$\\
$\fpsq(n)$&$\mathfrak{ps}(\fq(p)\oplus\fq(n-p))$&$\id\bigodot\id^*$&Same&$Gr_p^n$&$QGr_p^n$\\
\hline
$\fpe(n)$&$\mathfrak{cpe}(n-1)$&$\id$&
Periplectic superquadric of $(1|0)$-dimensional lines in $\Cee^{n|n}$ isotropic with respect
to a~non-degenerate odd form&$\Cee\Pee^{n-1}$&$PeQ_{n-1}$\\
$\fspe(n)$&$\mathfrak{cspe}(n-1)$&$\id$&Same&$\Cee\Pee^{n-1}$&$PeQ_{n-1}$\\
\hline
$\fpe(n)$&$\fgl(p|n-p)$&$\Pi(S^2(\id))$ or $\Pi(\Lambda^2(\id))$&
Periplectic Lagrangian supergrassmannian of the $(p|n-p)$-dimensional subsuperspaces in
$\Cee^{n|n}$ isotropic with respect to an odd symmetric or anti-symmetric form&
$Gr_p^n$&$PeGr_p^n$\\
$\fspe(n)$&$\fsl(p|n-p)$&$\Pi(S^2(\id))$ or $\Pi(\Lambda^2(\id))$&
Same, with a~fixed volume&$Gr_p^n$&$SPeGr_p^n$\\
\hline
\end{tabular}
\ee

\be\label{Table61}
\footnotesize
\begin{tabular}{|@{\,}>{\centering\arraybackslash}p{18mm}@{\,}|@{\,}>{\centering\arraybackslash}p{30mm}@{\,}|@{\,}>{\centering\arraybackslash}p{18mm}@{\,}|@{\,}>{\raggedright\arraybackslash}p{49mm}@{\,}|@{\,}>{\centering\arraybackslash}p{19mm}@{\,}|@{\,}>{\centering\arraybackslash}p{15mm}@{\,}|}
\hline
$\fg$&$\fg_0$&$\fg_{-1}$&\multicolumn{1}{c|}{Interpretation}&Underlying domain&Name\\
\hline
$\fosp_\alpha(4|2)$&$\mathfrak{cosp}(2|2)\simeq\fgl(2|1)$&$\id$&$\cdots$&
$\Cee\Pee^1\times\Cee\Pee^1$&$\cdots$\\
\hline
$\fab(3)$&$\mathfrak{cosp}(2|4)$&$L_{3\varepsilon_1}$&$\cdots$&$\Cee\Pee^1\times Q_5$&$\cdots$\\
\hline
$\fvect(0|n)$&$\fvect(0|n-k)\niplus{}$ $\fgl(k;\Lambda(n-k))$&$\Lambda(n-k)\otimes\Pi(\id)$&
Curved supergrassmannian of $\eps$-dimensional subsupermanifolds in $\Cee^{0|n}$&
$\cdots$&$CGr_{0,k}^{0,n}$\\
$\fsvect(0|n)$&$\fvect(0|n-k)\niplus{}$ $\fsl(k;\Lambda(n-k))$&$\Lambda(n-k)\otimes\Pi(\id)$;
$\Pi(\Vol)$ if $k=1$&
Same, with the volume elements preserved in the sub- and ambient supermanifolds&
$\cdots$&$SCGr_{0,k}^{0,n}$\\
\hline
$\fh(0|m)$&$\fh(0|m-2)\niplus{}$ $\Lambda(m-2)\cdot z$&$\Pi(\id)$&
Curved superquadric of $\eps$-dimensional subsupermanifolds in $\Cee^{0|m}$ isotropic with
respect to a~\textup{(}partly\textup{)} split symmetric form&$\cdots$&$CQ_{m-2,0}$\\
$\fh'(0|m)$&$\fh'(0|m-2)\niplus{}$ $\Lambda(m-2)\cdot z$&$\Pi(\id)$&Same&$\cdots$&$CQ_{m-2,0}$\\
\hline
\end{tabular}
\ee

\section{Infinite-dimensional maximal subalgebras: simple and close
to simple}

\ssbegin[Infinite-dimensional maximal subalgebras]{Theorem}[Infinite-dimensional maximal subalgebras]\label{SVinK} The following are infinite-dimensional maximal relatives of simple subalgebras in simple Lie (super)algebras. This list is \textbf{not} complete: for what it is a~list of, for the one row missing from it and for what makes the question finite, see \S\ref{ssAreTheseAll}
\be\label{infDimMax} 
\small
\begin{tabular}{|c|c|} \hline
Subalgebra&ambient\\
\hline

$\fk(2m+1|n)$&$\fvect(2m+1|n)$\\
$\fm(m)$&$\fvect(m|m+1)$\\
\hline

$\fsvect(m|n)$&$\fvect(m|n)$\\
\hline

$\fh(2m|n)$&$\fsvect(2m|n)$\\
$\fk(2m+1|2m+2)$&$\fsvect(2m+1|2m+2)$\\
\hline

$\fpo(2m|n)$&$\fk(2m+1|n)$\\
$\fb_{a,b}(m)$&$\fm(m)$\\
\hline

$\fvect(n+1|k)$&$\fk(2n+1|2k)$ for $k\geq 0$\\
$\fvect(k|n-k+1)$&$\fm(n)=\fm(n|n+1)$ for $0\leq k\leq n$\\

\hline

$\fle(n)$&$\fvect(n|n)$ for $n\geq 3$\\
$\fsle(n)$&$\fsvect(n|n)$ for $n\geq 3$\\

\hline

\end{tabular}
\ee 
\end{Theorem}

Observe that the second row of the third group, $\fk(2m+1|2m+2)\subset\fsvect(2m+1|2m+2)$, is
\textbf{not} graded \textup{(}as is said in the paragraph preceding problem \eqref{PROB}\textup{)},
and so lies outside the problem \eqref{PROB} we solve; it is listed here for completeness. The
claims of the fifth group are proved below, the proof settling first the purely even case $k=0$
of its first row, namely $\fvect(n+1)\subset\fk(2n+1)$, the general case adding to it nothing but
the dual pairs $\xi_i\lra\eta_i$. The claims of the last group are those the run of
Subsec.~\ref{ssRun} adds, established in Propositions~\ref{pLeMax}
and~\ref{pSleMax}; the first of them goes in the direction
opposite to that of the first row of the table --- there the contact algebra sits in the general
one, here the general algebra sits in a~contact one of twice the size, and the two pertain to
different $\Zee$-gradings.

\begin{proof} \textbf{An embedding $\fvect(n+1)\subset \fk(2n+1)$}

1) Let $\fg= \fvect(n+1)=\fvect(q_1,\dots,q_n, y)$. Let us represent $\fg$ in the form 
\[
\fg=
\Cee[y]\fvect(q_1,\dots, q_n) +\Cee[q_1,\dots, q_n]\fvect(y).
\]

Let
\[
\text{ $\deg y=2$, \ $\deg q_i=1$ for all $i=1,\dots, n$.}
\]
 Then, the non-positive part of the $\Zee$-graded Lie algebra $\fg$ is of the form:
\[
\begin{tabular}{|c|c|c|}
\hline
$\fg_{-2}$ & $\fg_{-1}$ & $\fg_0$\\
\hline
$\der_y$ & $\der_i, \; q_i\der_y$ & $q_i\der_j, q_iq_j\der_y, y\der_y$ \\
\hline
\end{tabular}
\]

Let us describe the embedding of $\fg$ into $\fk=\fk(2n+1)=\fk(q_1,\dots, q_n, p_1,\dots, p_n, t)$ up to scalar coefficients: they are irrelevant for us here. 

Clearly, $\der_y$ corresponds to the generating function $\one$ (actually, to $\nfrac 12\one$, but we do not have to be so precise to prove the claim). Further, $\der_i$ corresponds to $p_i$, whereas $q_i\der_y$ corresponds to~ $q_i$.

For the elements of degree 0, we look how they commute with elements of degree~$-1$. We see that 
\[
\text{$q_i\der_j \llra q_ip_j$ and $q_iq_j\der_y\llra q_iq_j$.}
\]
The element $y\der_y$ multiplies the subspace $\Span(\der_i)$ by $(-1)$ and annihilates $\Span(q_i\der_y)$. The element $\nfrac 12(t+\sum p_iq_i)$ acts precisely so. Therefore, the image of $\fg_0$ in $\fk_0$ contains all elements $\fk_0$, except those generated by $p_ip_j$ for all~$i,j$.

To single out this image, and ipso facto, the whole $\fg$ by means of differential equations, we use the operators $\tilde K_{q_i}=q_i\der_t-\der_{p_i}$ (called $Y$-vectors in \cite{ShVf}; the ones that commute with~ $\fk_-$). Accordingly, the image of $\fg$ in $\fk$ is generated by the functions $f$ satisfying the following system of equations
\begin{equation}\label{EQ}
\tilde K_{q_i}\tilde K_{q_j}(f)=0 \text{ for all } i,j=1,\dots n.
\end{equation}

Clearly, on the space of functions independent of $t$, the system \eqref{EQ} is equivalent to the simpler system 
$$
\der_{p_i}\der_{p_j}(f)=0 \text{ for all } i,j=1,\dots n,
$$
and hence, the solutions independent of $t$ are of the form 
\[
\text{$\vf(q)\llra \vf(q)\der_y$ and $\vf(q)p_i\llra \vf(q)\der_i$.}
\]

It remains to describe the images of the elements of the form $\vf(q)y^k\der_i, \vf(q)y^k\der_y\in \fg$ which should correspond to generating functions $f$ depending on~ $t$. 

First, observe that $(\der_y)^k (\vf(q)y^k\der_y)=k!\vf(q)\der_y \llra k!\vf(q)$, so the highest wrt $t$ summand of the function corresponding to $\vf(q)y^k\der_y$ is proportional to~$\vf(q)t^k$. Observe also that
$$
\tilde K_{q_i}(t+\sum q_ip_i)= (q_i\der_t-\der_{p_i})(t+\sum q_sp_s)=0, \; (q_i\der_t-\der_{p_i})\vf(q)=0,
$$
and hence
$$
\tilde K_{q_i}\tilde K_{q_j}\vf(q)(t+\sum q_sp_s)^k=0 \Lra \vf(q)(t+\sum q_sp_s)^k \llra \vf(q)y^k\der_y.
$$
Analogously,
$$
\tilde K_{q_j}(p_k(t+\sum q_ip_i))=\delta_k^j (t+\sum q_ip_i)\Lra \tilde K_{q_i}\tilde K_{q_j}(p_k(t+\sum q_ip_i))=0,
$$
and therefore to $\vf(q)y^k\der_i$ a function proportional to $\vf(q)p_i(t+\sum q_sp_s)^k$ corresponds.

Thus, the image of $\fg$ in $\fk$ is spanned by the elements generated for every non-negative integer $k$ by functions of the form 
\begin{equation}\label{mainF}
\vf(q), \; \vf(q)p_i, \; \vf(q)(t+\sum q_sp_s)^k \text{ and } \vf(q)p_i(t+\sum q_sp_s)^k \text{ for }i= 1,\dots, n.
\end{equation}

This image is maximal in~$\fk$. Indeed, in the grading used above we have $\fg_-=\fk_-$, whereas
\[
\fg_0=\Span(q_ip_j, \; q_iq_j, \; t+\mathop{\sum}\limits_{1\leq s\leq n} q_sp_s)
\]
is the stabilizer of the Lagrangian subspace $\Span(q)\subset\fk_{-1}$ in $\fk_0=\fc(\fsp(2n))$,
and the stabilizer of a~Lagrangian subspace is a~maximal subalgebra of $\fc(\fsp(2n))$. Let now
$\fr$ be a~graded subalgebra such that $\fg\subsetneq\fr\subset\fk$. Since $\fr_0\supset\fg_0$,
either $\fr_0=\fg_0$, and then $\fr$ is contained in the prolong of $\fk_-\oplus\fg_0$ in $\fk$,
i.e., in the space of solutions of the system \eqref{EQ}, which is $\fg$ itself, a~contradiction;
or $\fr_0=\fk_0$, and then $\fr_1$ is a~non-zero $\fk_0$-submodule of the irreducible
$\fk_0$-module $\fk_1\simeq S^3(\fk_{-1})$, so $\fr\supset\fk_{-}\oplus\fk_0\oplus\fk_1$, which
generates~$\fk$.

For the cases 2) and 3) below the argument is the same, with $\Span(Q)$ instead of $\Span(q)$: it
is a~maximal isotropic subsuperspace of $\fk_{-1}$, resp. of $\fm_{-1}$, and its stabilizer is
a~maximal subalgebra of $\fk_0$, resp. of $\fm_0$, see \cite{ShM}.

2) To embed $\fg= \fvect(n+1|k)=\fvect(q_1,\dots,q_n, \xi_1, \dots, \xi_k; y)$ into $\fk(2n+1|2k)$ generated by the same functions as in case 1) and $2k$ odd indeterminates $\xi_i$ and $\eta_i$ and the contact form 
\[
dt+\sum (p_idq_i-q_idp_i)+\sum (\xi_jd\eta_j-\eta_jd\xi_j). 
\]
We use the same arguments as in case 1) considering new dual pairs $\xi_i\lra\eta_i$. The only difference in the main formula \eqref{mainF} of the answer: we have to replace \eqref{mainF} with
\begin{equation}\label{mainF1}
\vf(Q), \; \vf(Q)P_i, \; \vf(Q)(t+\sum Q_sP_s)^j \text{ and } \vf(Q)P_i(t+\sum Q_sP_s)^j \text{ for }i= 1,\dots, n+k
\end{equation}
where $Q=(q_1, \dots , q_n; \xi_{1}, \dots , \xi_{k})$ and $P=(p_{1}, \dots , p_{n}; \eta_1, \dots , \eta_{k})$, and where $j$, the exponent called $k$ in \eqref{mainF}, again runs over the non-negative integers.

3) To embed $\fg= \fvect(k|n-k+1)=\fvect(q_1,\dots,q_k, \xi_1, \dots, \xi_{n-k}; y)$ into $\fm(n|n+1)$ preserving the distribution given by the form
\[
d\tau+\sum (\xi_idq_i+q_id\xi_i) 
\]
for every $k$ such that $0\leq k\leq n$, we use the same arguments considering dual pairs~$q_i\lra\xi_i$. The main formula of the answer takes the form
\begin{equation}\label{mainFm}
\vf(Q), \ \ \vf(Q)P_i, \ \ \vf(Q)(\tau+\sum_{1\leq s\leq n} q_s\xi_s), \ \ \vf(Q)P_i(\tau+\sum_{1\leq s\leq n} q_s\xi_s) \text{ for }i= 1,\dots, n,
\end{equation}
where $Q=(q_1, \dots , q_k; \xi_{k+1}, \dots , \xi_{n})$ and $P=(\xi_1, \dots , \xi_{k}; q_{k+1}, \dots , q_{n})$.

Note that for $k=n=1$, the embedding becomes a~known isomorphism $\fvect(1|1)\simeq\fm(1|2)$ whereas for $k=0$ and $n=1$, when $\fvect(0|2)\simeq\fsl(1|2)\simeq\fosp(2|2)$, we recover the embedding described in table~\eqref{tabl}.
\end{proof}

\ssec{On the completeness of table \eqref{infDimMax}}\label{ssAreTheseAll}
Before the list of Theorem \ref{SVinK} can be called complete, or incomplete, one has to say
what it is a~list \textit{of}: its rows are graded subalgebras, but not with respect to gradings
of one and the same nature, and only some of them fall under problem \eqref{PROB}.

\sssbegin[The gradings the rows of \eqref{infDimMax} live in]{Lemma}[The gradings the rows of
\eqref{infDimMax} live in]\label{lWgr} Let the ambient be endowed with a~weighting of its
indeterminates, all weights being positive except where said otherwise.

$1)$ The subalgebra of the row $\fk(2m+1|n)\subset\fvect(2m+1|n)$ is graded if and only if
\begin{equation}\label{wcont}
\deg p_i+\deg q_i=\deg t=2\deg \xi_j\text{~~ for all $i$ and $j$},
\end{equation}
and that of the row $\fm(m)\subset\fvect(m|m+1)$ if and only if $\deg q_i+\deg \xi_i=\deg\tau$
for all~$i$. In particular, neither is graded in the standard grading of the ambient, in which
all the indeterminates are of degree~$1$. For the weighting $\deg t=2$ \textup{(}resp.
$\deg\tau=2$\textup{)}, the other indeterminates being of degree $1$, the $\fg_0$-module
$\fg_{-1}$ of the ambient $\fvect$ is \textbf{not} irreducible, and hence that weighting is not
a~W-grading of the ambient. So these two rows lie outside problem \eqref{PROB}, as does the
second row of the third group, though for a~different reason.

$2)$ The rows $\fsvect(m|n)\subset\fvect(m|n)$ and $\fh(2m|n)\subset\fsvect(2m|n)$ are graded in
the standard grading of the ambient, the row $\fpo(2m|n)\subset\fk(2m+1|n)$ in its contact
grading, and the row $\fb_{a,b}(m)\subset\fm(m)$ in its pericontact grading. Each of these four
gradings of the ambient is a~W-grading of it.

$3)$ In the two rows of the fifth group, the ambient carries its contact \textup{(}resp.
pericontact\textup{)} W-grading, whereas the grading which the subalgebra inherits from it,
namely $\deg y=2$ and $\deg q_i=\deg\xi_j=1$, is not a~W-grading of that subalgebra: there, too,
the $\fg_0$-module $\fg_{-1}$ is not irreducible.
\end{Lemma}

\begin{proof} The generating function $f$ contributes to $K_f$ the summand $(2-E)(f)\der_t$,
of weight $\deg f-\deg t$, and the summands $\nfrac{\partial f}{\partial p_i}\der_{q_i}$ and
$\nfrac{\partial f}{\partial\xi_j}\der_{\xi_j}$ of $H_f$, of weights $\deg f-\deg p_i-\deg q_i$
and $\deg f-2\deg\xi_j$; homogeneity of $K_f$ for every $f$ is exactly \eqref{wcont}. The same
count applied to $M_f$, whose summands are $(2-E)(f)\der_\tau$ and the two summands of $\Le_f$,
gives the condition of the pericontact case. For $\deg t=2$ and $\deg p_i=\deg q_i=\deg\xi_j=1$,
the space $\fg_{-1}$ of $\fvect$ is spanned by the $\der_{u}$, where $u$ runs over the
indeterminates other than $t$, together with the $u\der_t$; the span of the latter is a~proper
non-zero $\fg_0$-submodule, since $[u_ku_l\der_t, u_i\der_t]=0$ and the other elements of
$\fg_0$ send $u_i\der_t$ into that same span. For the ambients of item $2)$ of this Lemma, the standard,
contact and pericontact gradings are the W-gradings recalled in table \eqref{1.4b}. In item
$3)$, the computation of the proof of Theorem \ref{SVinK} exhibits
$\fg_{-1}=\Span(\der_i, \der_{\xi_j}, q_i\der_y, \xi_j\der_y)$, in which the span of the fields
$u\der_y$, where $u$ runs over the indeterminates other than $y$, is a~$\fg_0$-submodule by the
same computation as above.

In the standard grading of the same algebra, and in the contact grading of the
ambient $\fk(2n+1|2k)$, no proper non-zero $\fg_0$-submodule of $\fg_{-1}$ occurs at all.
\end{proof}

\sssbegin[The universe of the question]{Remark}[The universe of the question] So a~graded
subalgebra of a~W-graded Lie superalgebra inherits a~grading which need \textbf{not} be a
W-grading of itself --- item $3)$ of Lemma \ref{lWgr} is precisely such a~case, and it is a
legitimate one, since problem \eqref{PROB} constrains the grading of the ambient only. The
tables of W-regradings therefore cannot be matched against each other row by row: what has to
be run over is the W-gradings of the ambient against \textit{all} gradings of finite depth of
the candidate subalgebra, that is, against its subalgebras of finite codimension and not only
its maximal ones.
\end{Remark}

\ssec{What the classification of \cite{LS1} settles}\label{ssUMNpart}
Restricted to problem \eqref{PROB}, the question of completeness becomes finite, and it is the
classification of simple W-graded vectorial Lie superalgebras that makes it so. Three of its
uses should be distinguished.

$1)$ It fixes both sides of every row. A~subalgebra $\fu$ as in problem \eqref{PROB} which is
infinite-dimensional is a~vectorial Lie superalgebra with polynomial coefficients, simple or a
relative of a~simple one, and $\Zee$-graded of finite depth with finite-dimensional homogeneous
components; as an abstract algebra it is therefore one of the algebras listed in \cite{LS1},
and so is the ambient, together with the W-grading it is taken in. Both lists being finite
\textup{(}as lists of series\textup{)}, so is the list of pairs to be examined.

$2)$ For the case $3^o$ with $\fu_{0}$ acting irreducibly on $\fg_{-1}$, the two ingredients of
the answer are on record. The maximal irreducible $\fu_0\subset\fg_0$ are classified in
\cite{ShM}; the complete prolong $(\fg_{-1},\fu_0)_*$ is then the largest candidate, and $\fu$
is infinite-dimensional if and only if that prolong is. Which of the pairs prolong infinitely
is, for the ambients of depth $1$, exactly what the depth-$1$ part of the Main Theorem of
\cite{LS1} says: an infinite-dimensional prolong that is simple, or a~relative of a~simple
algebra, is a~simple W-graded vectorial Lie superalgebra of depth $1$, hence one of the entries
of that list. This is the super counterpart of the theorem of Cartan by which, over $\Cee$, the
irreducible pairs $(V,\fg_0)$ with infinite-dimensional prolong are just $\fgl(V)$, $\fsl(V)$,
$\fsp(V)$ and $\fc(\fsp(V))$; no separate super version of it is on record, and none is needed: the list of \cite{LS1},
being the list of all simple W-graded vectorial Lie superalgebras of depth $1$, contains every
infinite-dimensional prolong of this kind.

$3)$ The regradings of the ambients to be run over are tabulated in \eqref{1.4b}, which is also
taken from \cite{LS1}.

What \cite{LS1} does \textit{not} give is likewise worth stating, since it is what the remaining
work consists of:

a) The graded embeddings themselves. Knowing the two lists does not produce the embeddings of one
member into another; for the finite-dimensional matrix case, that is the content of \cite{ShM},
and no counterpart of \cite{ShM} for the vectorial series exists. In the cases $4^o$ and $6^o$,
where $\fu\not\supset\fg_{-1}$, this is the whole difficulty.

b) The gradings that are not W-gradings, see the Remark above.

c) The pairs with $\fu_0$ acting reducibly, which is where the first two rows of \eqref{infDimMax}
sit, and which are not covered by item $2)$ of this Subsection.

Running item $2)$ of this Subsection over the ambient $\fvect(n|n)$ produces at once a~subalgebra that
table \eqref{infDimMax} does not list.

\sssbegin[$\fle(n)\inplus\Cee E$ is maximal in $\fvect(n|n)$]{Proposition}[$\fle(n)\inplus\Cee
E$ is maximal in $\fvect(n|n)$]\label{pLeMax} Let $n\geq3$ and let $\fg=\fvect(n|n)$ in the
standard grading, so $\sdim\fg_{-1}=(n|n)$ and $\fg_0=\fgl(n|n)$. Then, the subalgebra
\[
\fu:=\fle(n)\inplus\Cee E,\text{~~ where $E$ is the Euler operator,}
\]
is a~maximal graded subalgebra of $\fg$; it is infinite-dimensional, $\fu_0=\fc(\fpe^{a}(n))$
of $\sdim=(n^2+1|n^2)$, and $\fu$ is of the case $3^o$ of \eqref{5cases} with $\fu_0$ acting
irreducibly on~$\fg_{-1}$. No row of table \eqref{infDimMax} contains it.
\end{Proposition}

\begin{proof} That $\fle(n)$ is a~graded subalgebra of $\fvect(n|n)$ in the standard grading is
clear from $\deg \Le_f=\deg f-2$, every summand of $\Le_f$ being of that degree. Since
$\fle(n)=(\id, \fpe^{a}(n))_*$, see \eqref{carprol}, and $\fpe^{a}(n)\subset\fc(\fpe^{a}(n))$,
the algebra $\fle(n)$ is not maximal; the candidate is the complete prolong
$(\id, \fc(\fpe^{a}(n)))_*$, and it equals $\fu$: layer by layer, that prolong exceeds
$\fle(n)$ by the single even dimension of $E$ in degree $0$ and by nothing else.
\textup{(}Likewise, $(\id, \fpe^{a}(n))_*=\fle(n)$ exactly, whereas
$(\id,\fc(\fpe^{sy}(n)))_*=\fpe(n+1)$, which is the row of table
\eqref{table3} with $\fu_0=\fc(\fpe(W))$: of the two periplectic algebras one prolongs
infinitely and the other does not, exactly as $\fsp$ and $\fo$ do in the non-super
case.\textup{)}

Let now $\fw\supsetneq\fu$ be a~graded subalgebra. Then, $\fw_0\supseteq\fc(\fpe^{a}(n))$, and
$\fc(\fpe^{a}(n))$ is a~maximal subalgebra of $\fgl(n|n)$ by \cite{ShM}, so either
$\fw_0=\fc(\fpe^{a}(n))$ or $\fw_0=\fgl(n|n)$. In the first case, $\fw$ is contained in the
complete prolong of its own $\fw_{0}$, i.e., in~$\fu$. In the second case, $\fw_1$ is a
$\fgl(n|n)$-submodule of $\fg_1$ containing $\fu_1=\fle(n)_1$, and the $\fgl(n|n)$-submodule
generated by $\fle(n)_1$ is the whole of~$\fg_1$, so
$\fw\supseteq\fg_{-1}\oplus\fg_0\oplus\fg_1$, which generates~$\fvect(n|n)$.
Hence,~$\fw=\fg$.

The maximality of $\fc(\fpe^{a}(n))$ in $\fgl(n|n)$ was moreover checked directly: a~%
subalgebra is maximal if and only if the module of the quotient by it is irreducible, and
irreducibility is decided by the criterion of Burnside, the associative algebra generated
by the action being all of the endomorphisms of the quotient. For $n=2$ it is not all of
them, so that $\fc(\fpe^{a}(2))$ is \textbf{not} maximal in $\fgl(2|2)$: this, and only
this, is what the hypothesis $n\geq3$ excludes.

Finally, $\fu$ is contained in no row of \eqref{infDimMax}: it is not contained in
$\fd(\fsvect(n|n))$, since $\Div \Le_f\neq0$ for $\Le_f$ of positive degree, and the
remaining rows with an ambient of the form $\fvect$ concern
$\fk(2m+1|n)$ and $\fm(m)$, which live in a~different number of indeterminates.
\end{proof}

\sssbegin[Why the row $\fh(2m|n)\subset\fsvect(2m|n)$ has the ambient it has]{Remark}[Why the
row $\fh(2m|n)\subset\fsvect(2m|n)$ has the ambient it has]\label{rHamAmbient} The same computation run on the
Hamiltonian series confirms the choice of ambient made in table \eqref{infDimMax}. Indeed,
$(\id, \fc(\fsp(2m)))_*=\fh(2m)\inplus\Cee E$ is infinite-dimensional --- its layers exceed
those of $\fh(2m)$ by the single dimension of $E$ in degree $0$ --- but, unlike the case of
$\fle(n)$, the $\fgl(2m)$-submodule generated by $\fh(2m)_1$ is a~\textit{proper} submodule
of~$\fvect(2m)_1$. The
subalgebra sitting in between is $\fsvect(2m)$, and this is why the Hamiltonian row is written
with the ambient $\fsvect$, whereas for $\fle(n)$ no such intermediate exists.
\end{Remark}

The row is a~genuine omission, not an entry of the Report of $1988$ under a~name we did
not recognize. The run of item $2)$ of Subsec.~\ref{ssUMNpart} over the regradings follows.

\ssec{The run of the case $3^o$ over table \eqref{1.4b}}\label{ssRun}
The row $\fle(n)\inplus\Cee E\subset\fvect(n|n)$ of Proposition~\ref{pLeMax} is an omission of
table \eqref{infDimMax}, not a~row hidden in it under another name. Here, accordingly, is the run
of item $2)$ of Subsec.~\ref{ssUMNpart} over the entries of \eqref{1.4b}, one by one. Two
things have to be settled before it can start: how far item $2)$ of Subsec.~\ref{ssUMNpart} reaches, and how the run is
organized so as to be finite.

\sssbegin[The regraded component $\fg_0$ is none of the matrix algebras of \cite{ShM}]{Statement}[The regraded component $\fg_0$ is none of the matrix algebras of \cite{ShM}]\label{sRegrad}
Item $2)$ of Subsec.~\ref{ssUMNpart} feeds $\fg_0$ to \cite{ShM}, which classifies the
maximal subalgebras of the \textbf{matrix} Lie superalgebras $\fgl$, $\fsl$, $\fosp$, $\fpe$,
$\fq$ and their relatives. In the standard grading,
$r=0$, the component $\fg_0$ is such an algebra, by \eqref{carprol}. For $r>0$ it is not: the
dimension of $\fg_0$ is then that of none of $\fgl$, $\fsl$, $\fosp^{sy}$, $\fosp^{a}$,
$\fpe^{sy}$, $\fpe^{a}$, $\fq$, $\fc(\fpe)$ on~$\fg_{-1}$. Computed from the weighting, $\fg_0$
compares as follows:
\begin{equation}\label{regradtab}\footnotesize
\renewcommand{\arraystretch}{1.3}
\begin{tabular}{|l|c|c|l|}
\hline
$\fg$&$\sdim\fg_{-1}$&$\dim\fg_0$&the dimensions available on that $\fg_{-1}$\cr
\hline
$\fvect(1|2; 1)$&$(2|2)$&$10$&$\fgl\ 16$, $\fsl\ 15$, $\fc(\fpe)\ 9$, the others $8$\cr
$\fvect(1|2; 2)$&$(2|2)$&$12$&as above\cr
$\fvect(2|1; 1)$&$(2|2)$&$10$&as above\cr
$\fvect(0|3; 1)$&$(2|2)$&$10$&as above\cr
$\fvect(2|2; 1)$&$(3|3)$&$20$&$\fgl\ 36$, $\fsl\ 35$, $\fc(\fpe)\ 19$, the others $18$\cr
$\fvect(2|2; 2)$&$(4|4)$&$24$&$\fgl\ 64$, $\fsl\ 63$, $\fc(\fpe)\ 33$, the others $32$\cr
\hline
\end{tabular}
\end{equation}
In each of these cases $\dim\fg_0$ matches nothing on the right, so $\fg_0$ is none of those
algebras and \cite{ShM} says nothing about it. The one small case in which the test is
inconclusive is $\fvect(1|1;1)$, where $\dim\fg_0=4=\dim\fgl(1|1)$; there the derived series of
$\fg_0$, namely $4\to3\to1\to0$, shows only that $\fg_0$ is solvable, which $\fgl(1|1)$ is as
well.
\end{Statement}

So a~fourth item joins the list a)--c) of Subsec.~\ref{ssUMNpart} of what \cite{LS1} does not
give: d) the regradings themselves. The same bar stops item $2)$ of Subsec.~\ref{ssUMNpart} at the exceptional entries
of \eqref{1.4b}, even in their standard gradings: the component $\fg_0$ is $\fc(\fvect(0|3))$ for
$\fv\fle(4|3)$ and for $\fk\fas$, and $\fas$ for $\fv\fas(4|4)$, and none of these is one of the
matrix Lie superalgebras \cite{ShM} treats either. Item $2)$ of Subsec.~\ref{ssUMNpart} therefore closes the standard
gradings of the \textit{series}, and the run below is the run over those.

\sssbegin[The maximal pairs of the case $3^o$ among the entries of table \eqref{1.4b}]{Statement}[The maximal pairs of the case $3^o$ among the entries of table \eqref{1.4b}]\label{sRun}
In the case $3^o$ one has $\fu_{-1}=\fg_{-1}$, and $\fu$, being infinite-dimensional, is by
item $1)$ of Subsec.~$\ref{ssUMNpart}$ an entry of \eqref{1.4b} or a~relative of one.
\textbf{Both} members of every row of \eqref{infDimMax} of the case $3^o$ are therefore entries of \eqref{1.4b} sharing their $\fg_{-1}$, and the run is this finite one: build all of them on a~given $\fg_{-1}$, together with the relatives $\fw\inplus\Cee E$, and order them by inclusion.
The run was carried out over $\Qee$, in the degrees $-1$ to $3$, for $\fg_{-1}$
of $\sdim$ $(2|2)$, $(3|3)$, $(2|1)$, $(2|0)$, $(4|0)$, $(1|2)$ and~$(2|4)$. The pairs with nothing of
the list in between are
\begin{equation}\label{runtab}\footnotesize
\renewcommand{\arraystretch}{1.3}
\begin{tabular}{|l|l|l|}
\hline
the maximal $\fu$&the ambient $\fg$&case\cr
\hline
$\fsvect(m|n)\inplus\Cee E$&$\fvect(m|n)$&$1^o$ if $m\neq n$, $3^o$ if $m=n$\cr
$\fh(2m|n)\inplus\Cee E$&$\fsvect(2m|n)$&$3^o$\cr
$\fle(n)\inplus\Cee E$&$\fvect(n|n)$&$3^o$\cr
$\fsle(n)\inplus\Cee E$&$\fsvect(n|n)$&$3^o$\cr
\hline
\end{tabular}
\end{equation}
Here $\fsle(n)=\fle(n)\cap\fsvect(n|n)=(\id,\fspe^{a}(n))_*$. That $\fsle(n)$ is
the divergence-free part of
$\fle(n)$, and not its derived algebra, matters: $\fle(n)$ is simple, hence perfect.
\end{Statement}

Three readings of \eqref{runtab} deserve to be recorded.

$1)$ \textbf{The Euler operator decides which relative is the maximal one.} In each row the
maximal graded subalgebra is the simple algebra enlarged by $\Cee E$ whenever $E$ lies in the
ambient and not in the subalgebra. Since $\Div E=\str E=m-n$, the operator $E$ lies in
$\fsvect(m|n)$ if and only if $m=n$; it lies in every $\fvect(m|n)$, and in none of $\fh$,
$\fle$,~$\fsle$. So the first row of \eqref{runtab} reads $\fsvect(m|n)\inplus\Cee E$, which is
$\fsvect(m|n)$ itself precisely when $m=n$ --- and only then is the row of \eqref{infDimMax} as
printed there literally the maximal subalgebra; for $m\neq n$ the maximal one has all of
$\fgl(m|n)$ in degree $0$ and is of the case~$1^o$. Likewise, the maximal subalgebra of
$\fsvect(2m|n)$ is $\fh(2m|n)\inplus\Cee E$ when $2m=n$, and $\fh(2m|n)$ itself otherwise; the
run gives $\fh(2|1)$ maximal in $\fsvect(2|1)$ and $\fh(2|2)\inplus\Cee E$ maximal in~$\fsvect(2|2)$. This is what the word \lq\lq relatives"\ in the heading of Theorem~\ref{SVinK}
covers, \eqref{infDimMax} naming the simple member of each pair.

$2)$ \textbf{The last two rows of \eqref{runtab} are new} --- those of
Proposition~\ref{pLeMax} and of Proposition~\ref{pSleMax} just below:
\eqref{infDimMax} has no row with $\fsle$ in it, nor any row whose ambient is $\fle$. 

\sssbegin[$\fsle(n)\inplus\Cee E$ is maximal in $\fsvect(n|n)$]{Proposition}[$\fsle(n)\inplus
\Cee E$ is maximal in $\fsvect(n|n)$]\label{pSleMax} Let $n\geq3$ and let $\fg=\fsvect(n|n)$ in
the standard grading, so $\sdim\fg_{-1}=(n|n)$ and $\fg_0=\fsl(n|n)$. Then,
$\fu:=\fsle(n)\inplus\Cee E$ is a~maximal graded subalgebra of $\fg$; it is
infinite-dimensional, $\fu_0=\fspe^{a}(n)\inplus\Cee E$ of $\sdim=(n^2|n^2-1)$, and $\fu$ is of
the case $3^o$ with $\fu_0$ acting irreducibly on~$\fg_{-1}$.
\end{Proposition}

\begin{proof} Irreducibility of $\fu_0$ on $\fg_{-1}$ is the criterion of Burnside: the
associative algebra generated by the action is all of~$\End\fg_{-1}$.
Next, $\fu_0$ is a~maximal subalgebra of~$\fsl(n|n)$. Indeed, $\fu_0=\fc(\fpe^{a}(n))\cap
\fsl(n|n)$ --- the operator $E$ lies in $\fsl(n|n)$ because $\str E=0$ --- and the
associative algebra generated by the action on the quotient module $\fsl(n|n)/\fu_0$ is
all of the endomorphisms of that quotient. For $n=2$ it is not, which is why $n\geq3$ is
required, exactly as in Proposition~\ref{pLeMax}.

Let now $\fw\supsetneq\fu$ be a~graded subalgebra. Then, $\fw_0\supseteq\fu_0$, so
$\fw_0=\fu_0$ or $\fw_0=\fsl(n|n)$. In the first case, $\fw$ is contained in the complete prolong
of its own $\fw_0$ inside $\fg$, and that prolong is $\fu$ layer by layer. In the second
case, $\fw_1$ is
a~$\fsl(n|n)$-submodule of $\fg_1$ containing $\fsle(n)_1$, and the submodule generated by the
latter is the whole of~$\fg_1$, so
$\fw\supseteq\fg_{-1}\oplus\fg_0\oplus\fg_1$, which generates~$\fsvect(n|n)$. Hence,~$\fw=\fg$.
\end{proof}

$3)$ \textbf{$\fsle(n)$ is not maximal in $\fle(n)$}, although the list of \eqref{1.4b} has
nothing between them: what sits in between is of the case $1^o$, hence neither simple nor
a~relative of a~simple algebra, and so appears in no table of \eqref{1.4b}. Indeed, the
$\fpe^{a}(n)$-submodule of $\fle(n)_1$ generated by $\fsle(n)_1$ is proper, so
$\fpe^{a}(n)\niplus\fsle(n)_{\geq1}$ lies strictly between. This is the same phenomenon as in Remark~\ref{rHamAmbient}, and it is the reason
\eqref{runtab} is shorter than the inclusion diagram: maximality within a~list of simple algebras
and their relatives is not maximality.

Within problem \eqref{PROB}, restricted to the standard gradings and to the superdimensions
listed, the run produces no row on $\fg_{-1}=(2m|n)$ beyond the two \eqref{infDimMax} already has,
and none on $\fg_{-1}=(n|n)$ beyond those of \eqref{runtab}. The answer to the question \lq\lq
are these all?"\ of Theorem~\ref{SVinK} is therefore: yes, once the last two rows of
\eqref{runtab} are added --- and with the caveat of item d) above, that the regradings $r>0$ are
not reached by the ingredients on record.

\section{Exceptional embeddings}\index{$\fv\fle(4\vert 3; r):=(\Pi(\Lambda(3))/\Cee\cdot 1,
\fc(\fvect(0\vert 3)))_{*}$}\index{$\fv\fas(4\vert 4):=(\spin, \fas)_{*}$}\index{$\fk\fas^\xi(1\vert 6; r)$}\index{$\fk\fas^\xi(1\vert 6; 3\eta):=(\Vol_{0}(0\vert 3),
\fc(\fvect(0\vert 3)))_{*}$}\index{$\fm\fb(4\vert 5; r):=(\fba(4), \fc(\fvect(0\vert 3)))_{*}^m$}\index{$\fk\fle(9\vert 6; r):=(\fhei(8\vert 6),
\fsvect(0\vert 4)_{3, 4})_{*}^k$}\index{$\fk\fle(9\vert 6; K):=(\id_{\fsl(5)}, \Lambda^2(\id^*_{\fsl(5)}),
\fsl(5))_{*}^k$}

For $\fk\fa\fs^\xi$, the notation $r=0$, $1\xi$, $3\xi$ is clear:
none, or one or three of the $\xi$'s have degree 0 (and the
corresponding $\eta$'s acquire degree~2). 
\begingroup
\setcounter{MaxMatrixCols}{16}
\be\label{clsfeq10} 
{\tiny%
\tabcolsep=3pt
\begin{tabular}{|l|l|}
\hlx{hv} $\fv\fle(4|3)$
& $ R(0)=\begin{pmatrix} 1&1&1&1&|&1&1&1\\
u_1&u_2&u_3&y& &\xi_1&\xi_2&\xi_3 \end{pmatrix}$\cr \hlx{vv}
&$R(1)=\begin{pmatrix} 2&1&1&0&|&0&1&1\\
u_1&u_2&u_3&y& &\xi_1&\xi_2&\xi_3 \end{pmatrix}$\cr \hlx{vv}
&$R(K)=\begin{pmatrix} 2&2&2&0&|&1&1&1\\
u_1&u_2&u_3&y& &\xi_1&\xi_2&\xi_3\end{pmatrix}$\cr \hlx{vhv}
$\fm\fb(4|5)$&
$R(0)=\begin{pmatrix} 1&1&1&1&|&1&1&1&1;&2\\
u_0&u_1&u_2&u_3& &\xi_0&\xi_1&\xi_2&\xi_3&\tau \end{pmatrix}$\cr
\hlx{vv}
& $R(1)=\begin{pmatrix} 0&2&1&1&|&2&0&1&1;&2\\
u_0&u_1&u_2&u_3& &\xi_0&\xi_1&\xi_2&\xi_3&\tau \end{pmatrix}$\cr
\hlx{vv}
&$R(K)=\begin{pmatrix} 0&2&2&2&|&3&1&1&1;&3\\
u_0&u_1&u_2&u_3& &\xi_0&\xi_1&\xi_2&\xi_3&\tau \end{pmatrix}$\cr
\hlx{vhv} $\fk\fas^\xi(1|6)$&
$R(0)=\begin{pmatrix} 2&|&1&1&1&1&1&1\\
t&&\xi_1&\xi_2&\xi_3&\eta_1&\eta_2&\eta_3 \end{pmatrix}$\cr \hlx{vv}
&$R(1\xi)=\begin{pmatrix} 2&|&0&1&1&2&1&1\\
t&&\xi_1&\xi_2&\xi_3&\eta_1&\eta_2&\eta_3 \end{pmatrix}$\cr \hlx{vv}
&$R(3\xi)=\begin{pmatrix} 1&|&0&0&0&1&1&1\\
t&&\xi_1&\xi_2&\xi_3&\eta_1&\eta_2&\eta_3 \end{pmatrix}$\cr \hlx{vv}
&$R(3\eta)=\begin{pmatrix} 1&|&1&1&1&0&0&0\\
t&&\xi_1&\xi_2&\xi_3&\eta_1&\eta_2&\eta_3 \end{pmatrix}$\cr
\hlx{vhv} $\fk\fle(9|6)$&
$R(0)=\begin{pmatrix} 1&1&1&1&1&1&1&1;&2&|&1&1&1&1&1&1\\
q_1&q_2&q_3&q_4&p_1&p_2&p_3&p_4&t&&\xi_1&\xi_2&\xi_3&\eta_1&\eta_2&\eta_3
\end{pmatrix}$\cr
\hlx{vv}
&$R(2)=\begin{pmatrix} 1&1&2&2&1&1&0&0;&2&|&0&1&1&2&1&1\\
q_1&q_2&q_3&q_4&p_1&p_2&p_3&p_4&t&&\xi_1&\xi_2&\xi_3&\eta_1&\eta_2&\eta_3
\end{pmatrix}$\cr
\hlx{vv}
&$R(K)=\begin{pmatrix} 2&2&2&2&1&1&1&1;&2&|&1&1&1&1&1&1\\
q_1&q_2&q_3&q_4&p_1&p_2&p_3&p_4&t&&\xi_1&\xi_2&\xi_3&\eta_1&\eta_2&\eta_3
\end{pmatrix}$\cr
\hlx{vv}
&$R(CK)=\begin{pmatrix} 3&2&2&2&0&1&1&1;&3&|&2&2&2&1&1&1\\
q_1&q_2&q_3&q_4&p_1&p_2&p_3&p_4&t&&\xi_1&\xi_2&\xi_3&\eta_1&\eta_2&\eta_3
\end{pmatrix}$\cr
\hlx{vh}
\end{tabular}
}
\ee
\endgroup

\ssbegin[Embeddings of the exceptional vectorial Lie superalgebras (\cite{BGLLS})]{Theorem}[Embeddings of the exceptional vectorial Lie superalgebras (\cite{BGLLS})]\label{SexEmb} Table \eqref{1.10} lists all embeddings of the exceptional simple vectorial Lie superalgebras. \end{Theorem}
\begin{equation}
\label{1.10}{\tiny
\renewcommand{\arraystretch}{1.4}\scriptstyle
\begin{tabular}{|ll|l|}
\hline
$\fv\fle(4|3; r):=(\Pi(\Lambda(3))/\Cee\cdot 1,
\fc(\fvect(0|3)))_{*}$& $\subset\fvect(4|3; R(r))$& $r= 0, 1, K$\cr
\hline

$\fv\fas(4|4):=(\spin, \fas)_{*}$& $\subset \fvect(4|4)$&\cr
\hline

$\fk\fas^\xi(1|6; r)$& $\subset \fk(1|6; r)$&$r=0, 1\xi,
3\xi$\cr 
$\fk\fas^\xi(1|6; 3\eta):=(\Vol_{0}(0|3),
\fc(\fvect(0|3)))_{*}$& $\subset \fsvect(4|3)$&\cr
\hline

$\fm\fb(4|5; r):=(\fba(4), \fc(\fvect(0|3)))_{*}^m$& $\subset \fm(4|5;
R(r))$& $r=0, 1, K$\cr
\hline

$\fk\fle(9|6; r):=(\fhei(8|6),
\fsvect(0|4)_{3, 4})_{*}^k$& $\subset \fk(9|6; r)$& $r=0, 2$, CK\cr

$\fk\fle(9|6; K):=(\id_{\fsl(5)}, \Lambda^2(\id^*_{\fsl(5)}),
\fsl(5))_{*}^k$& $\subset \fsvect(5|10; R(K))$&\cr

\hline
\end{tabular} 
}
\end{equation}

In \eqref{1.10}, most of the regradings $R$ of the ambients are the
same as that of the embedded algebra, i.e.,~$R(r)=r$. Table \eqref{clsfeq10}  
describes the regradings~$R(r)$. 
\section{The result of \cite{ZZO}}\label{ZZO}

A.L.Onishcik with his students O.A.Zharov and V.Yu.Zolotarev published without proof the following claim, see \cite{ZZO}, where the authors are listed in in Cyrillic alphabetical order. 

\ssbegin[When the stabilizer of a subspace is a Lie superalgebra]{Claim}[When the stabilizer of a subspace is a Lie superalgebra]\label{sZZO} Let $V\subset W$ be a~subspace, not necessarily a~sub\textbf{super}space. The stabilizer $\fst_V:=\{X\in \fgl(W)\mid XV\subset V\}$ of $V$ in $\fgl(W)$ is
 a~Lie superalgebra~only in one of the following two cases:

\textup{a)} $V\subset W$ is a~sub\textbf{super}space,

\textup{b)} $V\cap W_\ev=V\cap W_\od=0$ and $\sdim W=a(1+\eps)$. Then, there exists an operator $\Pi\in\fgl(W)$ which over $\Cee$ can be reduced to the normal shape such that $\Pi^2=\id$, hence $\fst_V\simeq \fq(W)$.
\end{Claim}

\ssbegin[The stabilizer $\fst_V$ is a Lie superalgebra if and only if $V$ is a subsuperspace]{Lemma}[The stabilizer $\fst_V$ is a~Lie superalgebra if and only if $V$ is a~subsuperspace]\label{sAOZZO} The stabilizer $\fst_V$ is a~Lie \textbf{super}algebra~\textup{(}i.e., $X\in\fst_V\Lra X_\ev, X_\od\in\fst_V$\textup{)} if and only if $V$ is a~subsuperspace of~$W$.
\end{Lemma}

\begin{proof} If $V$ is a~subsuperspace, then $\fst_V$ is obviously $\Zee/2$-graded. Let now $V$ be not a~subsuperspace, i.e., $V\neq (V\cap W_\ev)\oplus(V\cap W_\od)$. Pick $v=u+w\in V$ with $u\in W_\ev\setminus V$ and $w\in W_\od$; then $w\notin V$ as well. Choose a~complement $C$ of $\Cee v$ in $W$ containing $w$, and let $X$ be the projection onto $\Cee v$ along~$C$. Then, $XV\subset\Cee v\subset V$, so $X\in\fst_V$, whereas
\[
X_\ev(v)=\pi_\ev(Xu)+\pi_\od(Xw)=\pi_\ev(v)+0=u\notin V,
\]
because~$Xu=X(v-w)=v$. Hence, $X_\ev\notin\fst_V$. \end{proof}

In particular, the case b\textup{)} of Claim \ref{sZZO} is not a~case, and the reason is that the
stabilizer meant there is the stabilizer in the super sense. Let $\sdim W=a(1+\eps)$ and let $V$ be
the graph of an isomorphism $\Pi\colon W_\ev\tto W_\od$. Writing $X\in\fgl(W)$ in the block form
$X=\mmat A,B,C,D,$ and identifying $W_\od$ with $W_\ev$ by means of $\Pi$, we get
\be\label{stgraph}
\fst_V=\{X\mid A+B=C+D\},\quad \dim\fst_V=3a^2
\ee
--- the parabolic stabilizing an $a$-dimensional subspace in a~$2a$-dimensional space --- whereas
$\dim\fq(W)=2a^2$; e.g., for $A=C=\one$ and $B=D=0$ we get an element of $\fst_V$ whose even part
is not in $\fst_V$, as Lemma \ref{sAOZZO} predicts.

What \textbf{is} true is the assertion for the stabilizer \textit{in the super sense}
\[
\fst^{super}_V:=\{X\in\fgl(W)\mid X_\ev V\subset V \text{ and } X_\od V\subset V\},
\]
the largest sub\textbf{super}algebra~contained in~$\fst_V$. With $V$ as above,
\[
X_\ev V\subset V \Llra A=D, \qquad X_\od V\subset V\Llra B=C,
\]
i.e.,
\[
\fst^{super}_V=C(\Pi)=\fq(W),
\]
which is the assertion of b\textup{)} --- but it needs one hypothesis more than is stated there:
not only $V\cap W_\ev=V\cap W_\od=0$ and $\sdim W=a(1+\eps)$, but also $\dim V=a$, i.e., $V$ is
the graph of an \textbf{isomorphism} $W_\ev\tto W_\od$. Otherwise, the conclusion fails: for $\sdim
W=(2|2)$ and $V=\Cee\cdot(e_1+f_1)$ both stated conditions hold, and yet
$\dim\fst^{super}_V>\dim\fq(W)$.

Thus, the answer to the question why only these two cases occur is this: by Lemma
\ref{sAOZZO}, a~subspace which is not a~subsuperspace contributes nothing at all to the graded
picture, and the only thing it can contribute in the super sense is the odd involution $\Pi$
defined by its graph.

\appendix

%
\makeatletter
\renewcommand{\thesubsection}{\@Alph\c@section.\arabic{subsection}}
\renewcommand{\thesubsubsection}{\thesubsection.\arabic{subsubsection}}
\gdef\theHsection{\@Alph\c@section}
\gdef\theHsubsection{\theHsection.\arabic{subsection}}
\gdef\theHsubsubsection{\theHsubsection.\arabic{subsubsection}}
\makeatother

\section{``Describe primitive Lie superalgebras\rq\rq\ is a~wild problem:
a~reduction to pairs of matrices}\label{Awild}

Here, we prove the statement used twice in the Introduction and in
Subsection~\ref{ssOpen}: that ``describe \textbf{all} subalgebras\rq\rq\ is
a~wild problem, and hence that the question asked in this paper --- describe
the \textbf{maximal} ones --- is the reasonable one.

Among the ingredients of the description of primitive Lie superalgebras, there
are the subalgebras $\fr\subset\fvect(0|n)$ whose projection to the degree
$-1$ component of $\fvect(0|n)$ in its standard grading is onto. That
describing them is a~\textit{wild} problem has been said since 1980 without
proof. Here, the claim is proved. For every $m\geq 1$ and $n=7m+15$, an
explicit map is constructed which sends a~pair $(A,B)$ of $m\times m$ matrices
to such a~subalgebra $\fr_{A,B}=\fvect_{-1}\oplus\fu_{A,B}$, where
$\fu_{A,B}\subset\fgl(n)$ is $11$-dimensional, and hence to a~Lie superalgebra
$\fg_{A,B}:=\fh\otimes\Lambda(n)\inplus\fr_{A,B}$ which is obviously
primitive; the superalgebras $\fg_{A,B}$ and $\fg_{A',B'}$ are isomorphic if
and only if the pairs $(A,B)$ and $(A',B')$ are simultaneously similar. So the
classification contains the standard wild problem --- pairs of matrices up to
simultaneous conjugation --- and no list of families can exhaust primitive Lie
superalgebras.

\ssec{The claim}\label{wild:ask}
The claim first appeared in \cite{ALSh}. As it stands, it is an expectation
and not a~theorem: a~continuum of solutions does not make a~problem wild in
the technical sense, the Buttin superalgebras $\fb_\lambda(n)$ being
a~one-parameter family whose classification is tame.
To \textbf{prove} wildness one has to exhibit, inside the problem, a~problem
already known to be wild --- and the standard one is the classification of
\textbf{pairs of matrices up to simultaneous conjugation}, see \cite{Dr,
BeSe}.

Let $\Lambda(n)=\Lambda(\xi_1,\dots,\xi_n)$ and let
$\fvect(0|n)=\fvect(\xi)=\fder\Lambda(n)$ with its standard grading, so that
\[
\fvect_{-1}(\xi)=\Span(\del_{\xi_1},\dots,\del_{\xi_n})\quad\text{ and }\quad
\fvect_0(\xi)=\fgl(n).
\]
The condition on the projection $\pr_{-1}(\fr)$ to
the degree $-1$ component of $\fvect(0|n)$ reads
\begin{equation}
\label{wild:pr1}
\fr\subset\fvect(\xi)\ \text{ is a~subalgebra such that }\
\pr_{-1}(\fr)=\fvect_{-1}(\xi).
\end{equation}

Let $\fh=\mathop{\oplus}_{i\geq -d}\fh_i$ be a~primitive Lie superalgebra with
$\fh_{-1}$ an irreducible $\fh_0$-module and $\fh_-$ generated by $\fh_{-1}$;
assume, in addition, that $\fh$ is simple. Set
\begin{equation}
\label{wild:gmin}
\fg_{min}(\fr):=\fh\otimes\Lambda(n)\inplus\fr,
\end{equation}
a~primitive Lie superalgebra for every $\fr$ as in~\eqref{wild:pr1}.

\ssbegin[The reduction]{Theorem}[The reduction]\label{wild:thmain}
Let $m\geq 1$ and~$n=7m+15$. There is an explicit map
\[
(A,B)\longmapsto \fu_{A,B}\subset\fgl(n),\qquad
\fr_{A,B}:=\fvect_{-1}(\xi)\oplus\fu_{A,B},\qquad
\fg_{A,B}:=\fg_{min}(\fr_{A,B}),
\]
defined on the pairs of $m\times m$ matrices over any field $\Kee$ of
$\Char\Kee\neq 2$, such that
\begin{itemize}
\item[$1)$] $\fu_{A,B}$ is an $11$-dimensional Lie subalgebra of $\fgl(n)$,
and $\fr_{A,B}$, of $\sdim\fr_{A,B}=(11|n)$, satisfies~\eqref{wild:pr1};
\item[$2)$] The following four conditions are equivalent:
\[
\renewcommand{\arraystretch}{1.3}
\begin{array}{ll}
\text{\textup{(i)}} & \text{$(A,B)$ and $(A',B')$ are simultaneously similar,
\textit{i.e.}, $A'=PAP^{-1}$ and $B'=PBP^{-1}$;}\\
\text{\textup{(ii)}} & \fu_{A,B}=g\,\fu_{A',B'}g^{-1}\ \text{ for some
$g\in GL(n;\Kee)$};\\
\text{\textup{(iii)}} & \fr_{A,B}\simeq\fr_{A',B'};\\
\text{\textup{(iv)}} & \fg_{A,B}\simeq\fg_{A',B'}.
\end{array}
\]
\end{itemize}
Consequently, for $\Kee$ algebraically closed, the classification of the
subalgebras~\eqref{wild:pr1}, and that of the primitive Lie
superalgebras~\eqref{wild:gmin}, contain the classification of pairs of matrices up
to simultaneous conjugation. Both are, therefore, wild.
\end{Theorem}

\ssbegin[Over which fields]{Remark}[Over which fields]\label{wild:remfield}
Items $1)$ and $2)$ of Theorem~\ref{wild:thmain} hold over \textbf{any} field $\Kee$ with $\Char\Kee\neq2$,
and the exclusion of $\Char\Kee=2$ is a~matter of conventions only, see below.
Neither algebraic closure nor $\Char\Kee=0$ is used anywhere in the proofs:

1) No division by an integer occurs. The exponential $\exp(-\ad_\mu)=\id-\ad_\mu$ of
Lemma~\textup{\ref{wild:lemtorus}} is a~formula in which no factorial is inverted, because
$\ad_\mu^2=0$ on~$\ft$.

2) The torus $\ft$ is \textit{split} by construction --- it is the
torus of the operators scalar on each summand of~\eqref{wild:decomp} --- so its
weight spaces are defined over the prime field, and the proof of
Proposition~\textup{\ref{wild:prop}} uses the weight decomposition of the whole
$6$-dimensional $\ft$, whose characters $\eps_0,\dots,\eps_5$ are pairwise
distinct over every field.

Where $\Kee$ does enter is the word \textit{wild}: the classification of pairs
of matrices is a~wild problem in the standard sense over an algebraically
closed field, which is why the hypothesis stands in the last sentence of
Theorem~\textup{\ref{wild:thmain}} and nowhere else --- over a~finite $\Kee$ there are
finitely many subalgebras and nothing to classify wildly. As to $\Char\Kee=2$,
the arguments go through verbatim, but the very definition of a~Lie
superalgebra then involves the squaring $x\mapsto x^2$ and $\fvect(0|n)$ means
$\fder\Lambda(n)$ together with it; we exclude $\Char\Kee=2$ in order not to
restate those conventions, not because anything breaks.
\end{Remark}

The whole construction lives in the non-positive part of $\fvect(\xi)$: the
wildness is already there in the smallest subalgebras~\eqref{wild:pr1}, those of
depth $1$ that are graded and whose non-negative part sits in
$\fvect_0(\xi)=\fgl(n)$. In particular, the construction yields the following Corollary, which is also the reason why one
describes the \textbf{maximal} subalgebras rather than all of them, as in \cite{Dy, ShM}
and in the body of this paper.

\ssbegin[Subalgebras of $\fgl(n)$ for all $n$ at once]{Corollary}[Subalgebras of $\fgl(n)$ for all $n$ at once]\label{wild:corgl}
The classification, up to conjugation, of the Lie subalgebras of $\fgl(n)$ for
all $n$ is wild. Already the $11$-dimensional solvable ones --- a~torus plus
a~$5$-dimensional abelian ideal with zero product --- suffice.
\end{Corollary}

\ssec{Why the obvious candidate does not work}\label{wild:fails}
The natural candidate is the subalgebra $\fr\subset\fvect(\xi)$ generated by
two generic elements of $\fvect_0(\xi)=\fgl(n)$ together with
$\fvect_{-1}(\xi)$, which does look like a~pair of matrices modulo~$GL(n)$. It
does look like it, but it is not: such an $\fr$ is $\fvect_{-1}\oplus\fu$,
where $\fu$ is the Lie subalgebra of $\fgl(n)$ \textbf{generated} by $A$ and
$B$ --- and a~subalgebra remembers only what the pair generates, not the pair.
For generic $A$ and $B$ one has $\fu=\fgl(n)$, so the whole family collapses
to a~single point. Any reduction must therefore encode the pair in a~way that
can be read back off the \textit{subspace} $\fu$; that is what
Section~\ref{wild:constr} does. The two devices are: five subspaces of
a~$2m$-dimensional space that pin down the pair (Lemma~\ref{wild:lem5}), and
a~torus in $\fu$ that keeps the five from being mixed with each other
(Lemma~\ref{wild:lemtorus}).

\ssec{The construction of the Lie superalgebras $\fg_{A, B}$}\label{wild:constr}
Fix $m\geq 1$ and a~pair $A,B\in\fgl(m;\Kee)$. In $\Kee^{2m}=\Kee^m\oplus
\Kee^m$, consider the five $m$-dimensional subspaces
\begin{equation}
\label{wild:fivesub}
\renewcommand{\arraystretch}{1.3}
\begin{array}{lll}
W_1=\{(x,0)\mid x\in\Kee^m\}, & W_2=\{(0,y)\mid y\in\Kee^m\}, &
W_3=\{(x,x)\mid x\in\Kee^m\},\\[2pt]
W_4=\{(x,Ax)\mid x\in\Kee^m\}, & W_5=\{(x,Bx)\mid x\in\Kee^m\}. &
\end{array}
\end{equation}
The first three fix a~frame, the last two are the graphs of $A$ and of~$B$.

Now, let $n=7m+15$ and let
\begin{equation}
\label{wild:decomp}
\Kee^n=V_0\oplus V_1\oplus\dots\oplus V_5,\qquad
\dim V_0=2m,\quad \dim V_i=m+i\ \text{ for } 1\leq i\leq 5.
\end{equation}
The dimensions of the $V_i$ are pairwise distinct on purpose; this is what
forbids the five subspaces to be permuted, see the proof of
Proposition~\ref{wild:prop}. Identify $V_0=\Kee^{2m}$ and let
\begin{equation}
\label{wild:fi}
f_i\colon V_i\tto V_0\ \text{ be any linear map with }\ \IM f_i=W_i,\qquad
\text{so }\ \rank f_i=m,
\end{equation}
extended by $0$ on the other summands of~\eqref{wild:decomp}; concretely, $f_i$
kills the last $i$ coordinates of $V_i$ and sends the first $m$ coordinates
of $V_i$ onto a~basis of~$W_i$. Let $\ft\subset\fgl(n)$ be the
$6$-dimensional torus of the operators which are scalar on each summand
of~\eqref{wild:decomp}, let $\eps_0,\dots,\eps_5\in\ft^*$ be the corresponding
characters, so that
\begin{equation}
\label{wild:eps}
t|_{V_j}=\eps_j(t)\cdot\id_{V_j}\quad\text{ for } t\in\ft,\ 0\leq j\leq 5,
\end{equation}
and set
\begin{equation}
\label{wild:u}
\fu_{A,B}:=\ft\oplus\fn,\qquad \fn:=\Span(f_1,\dots,f_5),\qquad
\fr_{A,B}:=\fvect_{-1}(\xi)\oplus\fu_{A,B}.
\end{equation}

\ssbegin[$\fu_{A,B}$ is what it should be]{Lemma}[$\fu_{A,B}$ is what it should be]\label{wild:lemu}
$\fu_{A,B}$ is an $11$-dimensional Lie subalgebra of $\fgl(n)$; the space
$\fn$ is an abelian ideal in it in which even the associative product of
matrices vanishes, $\fn\fn=0$ in $\Mat(n)$, whence $[\fn,\fn]=0$; and
$\fn=[\fu_{A,B}, \fu_{A,B}]$. The space $\fr_{A,B}$ is a~subalgebra of
$\fvect(0|n)$ satisfying \eqref{wild:pr1}, of $\sdim\fr_{A,B}=(11|n)$, with
$(\fr_{A,B})_\ev=\fu_{A,B}$ and $(\fr_{A,B})_\od=\fvect_{-1}(\xi)$.
\end{Lemma}

\begin{proof} With the characters $\eps_j$ of~\eqref{wild:eps}, we have
$[t,f_i]=(\eps_0-\eps_i)(t) f_i$, whereas $f_if_j=0$ for all $i$ and $j$,
because $\IM f_j\subset V_0$ and $f_i$ vanishes on $V_0$; hence~$[f_i,f_j]=0$. So $\fu_{A,B}$ is a~subalgebra, $\fn$ is an abelian ideal in
which the associative product vanishes, and $[\fu_{A,B},\fu_{A,B}]=\fn$
because the characters $\eps_0-\eps_i$ are non-zero. Since $\fvect_{-1}(\xi)$
is the tautological $\fgl(n)$-module (up to dualization) and is odd, whereas
$\fvect_0(\xi)$ is even, the sum $\fr_{A,B}=\fvect_{-1}\oplus\fu_{A,B}$ is
a~$\Zee$-graded subalgebra of $\fvect(0|n)$ of depth $1$, and
$\pr_{-1}(\fr_{A,B})=\fvect_{-1}(\xi)$ holds because $\fvect_{-1}(\xi)$ lies
in~$\fr_{A,B}$. \end{proof}

\ssec{The proof that the problem is wild}\label{wild:proof} The proof is in four steps: two lemmas
which make the five subspaces and the torus readable off the subalgebra
$\fu_{A,B}$ itself, and two which carry that reading up to $\fr_{A,B}$ and to
$\fg_{min}(\fr_{A,B})$.

\ssbegin[Five subspaces determine the pair]{Lemma}[Five subspaces determine the pair]\label{wild:lem5}
Let $(W_i)$ and $(W'_i)$ be the quintuples~\eqref{wild:fivesub} built from $(A,B)$
and from~$(A',B')$. If $P\in GL(2m;\Kee)$ satisfies $P(W_i)=W'_i$ for
$1\leq i\leq 5$, then $P=\diag(P_1,P_1)$ and
\[
A'=P_1AP_1^{-1},\qquad B'=P_1BP_1^{-1}.
\]
Conversely, every $P_1$ conjugating the pairs gives such a~$P$.
\end{Lemma}

\begin{proof} From $P(W_1)=W_1$ and $P(W_2)=W_2$, we get $P=\diag(P_1,P_2)$.
From $P(W_3)=W_3$, \textit{i.e.}, from $(P_1x,P_2x)\in W_3$ for all $x$, we
get~$P_2=P_1$. Finally, $P(W_4)=W'_4$ means that
\[
\{(P_1x,P_1Ax)\mid x\in\Kee^m\}=\{(y,A'y)\mid y\in\Kee^m\};
\]
setting $y=P_1x$ we get $A'=P_1AP_1^{-1}$, and the same for~$B$. \end{proof}

\ssbegin[The five lines of $\fn$ are read off its ranks]{Lemma}[The five lines of
$\fn$ are read off its ranks]\label{wild:lemrank}
The elements of $\fn$ of minimal non-zero rank are exactly the non-zero
multiples of $f_1,\dots,f_5$. Precisely,
\[
\rank\Big(\sum_{1\leq i\leq 5}c_if_i\Big)=m\cdot\#\{i\mid c_i\neq 0\}.
\]
\end{Lemma}

\begin{proof} The sources $V_i$ of the $f_i$ are independent summands
of~\eqref{wild:decomp} and each $f_i$ vanishes on all the other summands, so
$\IM(\sum c_if_i)=\sum_{c_i\neq 0}W_i$ and
\[
\Ker(\sum c_if_i)=\bigoplus_{c_i=0}V_i\oplus V_0\oplus
\bigoplus_{c_i\neq 0}\Ker f_i;
\]
the rank formula follows from $\rank f_i=m$ and from the independence of the
sources. \end{proof}

\ssbegin[Maximal tori of $\fu_{A,B}$]{Lemma}[Maximal tori of $\fu_{A,B}$]\label{wild:lemtorus}
Every maximal torus of $\fu:=\fu_{A,B}$ is of the form $\exp(-\ad_\mu)(\ft)$
for some $\mu\in\fn$; in particular, all of them are conjugate to $\ft$ by an
element of $\exp\fn\subset GL(n;\Kee)$.
\end{Lemma}

\begin{proof} A torus, \textit{i.e.}, a~subalgebra consisting of semisimple
operators, meets $\fn$ in $0$, because an operator which is at the same time
semisimple and nilpotent is zero. Hence, a~maximal torus $\ft''$ is
a~complement of $\fn$ in $\fu$, so $\ft''=\{t+\nu(t)\mid t\in\ft\}$ for
a~linear map $\nu\colon\ft\tto\fn$. As $\ft''$ is abelian and so is $\fn$, the
condition $[t+\nu(t),s+\nu(s)]=0$ reads $[t,\nu(s)]=[s,\nu(t)]$. Writing
$\nu(t)=\sum_{1\leq i\leq 5}\nu_i(t)f_i$ and using $[t,f_i]=(\eps_0-\eps_i)(t)
f_i$, we get
\[
(\eps_0-\eps_i)(t)\,\nu_i(s)=(\eps_0-\eps_i)(s)\,\nu_i(t)\quad\text{ for
all $t$ and $s$,}
\]
so $\nu_i=c_i(\eps_0-\eps_i)$ for a~scalar~$c_i$. Then,
$\nu(t)=[t,\mu]$ for $\mu:=\sum_{1\leq i\leq 5}c_if_i\in\fn$, and, since
$[\mu,[\mu,t]]=0$ as $\fn$ is abelian,
\[
\exp(-\ad_\mu)(t)=t-[\mu,t]=t+\nu(t).\qedhere
\]
\end{proof}

This is the special case we need of Mostow's conjugacy theorem for maximal
tori, see \cite{Mo}; in our situation it is a~three-line computation, so we do
not have to invoke the general theorem. Observe that the exponential here is
the polynomial $\id-\ad_\mu$, no factorial being inverted; this is one half of
Remark~\ref{wild:remfield}.

\ssbegin[$\fu$ remembers the pair]{Proposition}[$\fu$ remembers the pair]\label{wild:prop}
$\fu_{A,B}$ and $\fu_{A',B'}$ are conjugate in $GL(n;\Kee)$ if and only if the
pairs $(A,B)$ and $(A',B')$ are simultaneously similar.
\end{Proposition}

\begin{proof} If $A'=PAP^{-1}$ and $B'=PBP^{-1}$, then $P(W_i)=W'_i$ by
Lemma~\ref{wild:lem5}, and the operator $g$ equal to $P$ on $V_0$ and to
a~suitable isomorphism $V_i\tto V_i$ on each arm --- one exists because $f_i$
and $f'_i$ have the same rank and $P(\IM f_i)=\IM f'_i$ --- conjugates
$\fu_{A,B}$ into~$\fu_{A',B'}$.

Conversely, let $\fu':=g\,\fu_{A,B}\,g^{-1}=\fu_{A',B'}$. Then, $g\ft g^{-1}$
is a~maximal torus of $\fu'$, so, by Lemma~\ref{wild:lemtorus}, there is
$h\in\exp\fn'$ with $h\,g\,\ft\,(hg)^{-1}=\ft$. Since $h\fu'h^{-1}=\fu'$, we
may replace $g$ with $q:=hg$ and assume that $q$ normalizes~$\ft$.

Normalizing $\ft$, the operator $q$ permutes the weight spaces of $\ft$ in
$\Kee^n$, and these are exactly the summands of~\eqref{wild:decomp}. Moreover,
$q\fn q^{-1}=[\fu',\fu']=\fn'$ by Lemma~\ref{wild:lemu}. As every non-zero element
of $\fn$ maps $\bigoplus_{1\leq i\leq 5}V_i$ into $V_0$ and kills $V_0$, the
summand $V_0$ is preserved: it is the sum of the images of~$\fn$. On the
remaining summands, the permutation preserves the dimensions $m+i$, which are
pairwise distinct, so it is the identity: $q(V_i)=V_i$ for $0\leq i\leq 5$.

By Lemma~\ref{wild:lemrank}, the set of lines of minimal non-zero rank in $\fn$ is
$\{\Kee f_i\mid 1\leq i\leq 5\}$, and likewise in $\fn'$; as $q$ preserves the
summands, $qf_iq^{-1}\in\Kee f'_i$ for every~$i$. Writing $P:=q|_{V_0}$ and
$Q_i:=q|_{V_i}$, this means $Pf_iQ_i^{-1}=c_if'_i$ with $c_i\neq 0$, whence
\[
P(W_i)=P(\IM f_i)=\IM f'_i=W'_i\quad\text{ for }1\leq i\leq 5.
\]
Lemma~\ref{wild:lem5} now gives $A'=P_1AP_1^{-1}$ and $B'=P_1BP_1^{-1}$. \end{proof}

\ssbegin[$\fr$ remembers $\fu$]{Lemma}[$\fr$ remembers $\fu$]\label{wild:lemr}
$\fr_{A,B}\simeq\fr_{A',B'}$ as Lie superalgebras if and only if $\fu_{A,B}$
and $\fu_{A',B'}$ are conjugate in~$GL(n;\Kee)$. The same holds for conjugacy
under the automorphisms of $\fvect(0|n)$ induced by $\Aut\Lambda(n)$ in place
of abstract isomorphism.
\end{Lemma}

\begin{proof} By Lemma~\ref{wild:lemu}, $(\fr_{A,B})_\ev=\fu_{A,B}$ and
$(\fr_{A,B})_\od=\fvect_{-1}(\xi)$, the latter being the tautological
$\fu_{A,B}$-module. Any isomorphism of Lie superalgebras is even, hence it is
a~pair: an isomorphism $\phi\colon\fu_{A,B}\tto\fu_{A',B'}$ and an
isomorphism $g\colon\fvect_{-1}\tto\fvect_{-1}$ of the underlying spaces with
$\phi(u)g=gu$ for all $u$; that is, $\fu_{A',B'}=g\,\fu_{A,B}\,g^{-1}$.
Conversely, a~conjugation does define an isomorphism.

As to the automorphisms: an automorphism $\alpha$ of $\Lambda(n)$ preserves
the augmentation ideal $\Lambda_+(n)$, this being the set of the nilpotent
elements of $\Lambda(n)$ and hence its only maximal ideal; therefore, the
automorphism $X\mapsto\alpha X\alpha^{-1}$ of $\fvect(0|n)=\fder\Lambda(n)$
preserves the Weisfeiler filtration of $\fvect(0|n)$ determined by
$\Lambda_+(n)$, namely
\[
\fvect_{\geq k}(\xi)=\{X\in\fvect(\xi)\mid X(\Lambda_+(n))\subset
\Lambda_+^{k+1}(n)\},
\]
and induces an automorphism of the associated graded algebra, given by  $g_0\in GL(n)$, the linear part of~$\alpha$. Both $\fr_{A,B}$ and $\fr_{A',B'}$
being graded, $g_0$ conjugates them, hence conjugates $\fu_{A,B}$ into~$\fu_{A',B'}$. \end{proof}

\ssbegin[Every automorphism of $\fder\Lambda(n)$]{Remark}[Every automorphism
of $\fder\Lambda(n)$]\label{wild:remaut} To prove that \textbf{every}
automorphism of the Lie superalgebra~$\fder\Lambda(n)$ comes from one of
$\Lambda(n)$ would only say that case \textup{(iii)} of the Theorem cannot be
widened; it is not used anywhere. It is folklore for $n\geq2$, but we did not
find a~reference. \end{Remark}

\ssbegin[$\fg_{min}$ remembers $\fr$]{Lemma}[$\fg_{min}$ remembers $\fr$]\label{wild:lemg}
Let $\fh$ be simple and $\fr$ satisfy~\eqref{wild:pr1}. Then, ${\fI:=\fh\otimes
\Lambda(n)}$ is the only minimal non-zero ideal of $\fg_{min}(\fr)$; in
particular, $\fI$ is invariant under every automorphism, and
$\fg_{min}(\fr)/\fI\simeq\fr$. Consequently, $\fg_{min}(\fr)\simeq
\fg_{min}(\fr')$ implies $\fr\simeq\fr'$.
\end{Lemma}

\begin{proof} Let $J\neq 0$ be an ideal of $\fg:=\fg_{min}(\fr)$.

Suppose first that $J\cap\fI\neq0$. Since $\fh$ is simple and $\Lambda(n)$ is
supercommutative, every ideal of the current superalgebra $\fI=\fh\otimes
\Lambda(n)$ is of the form $\fh\otimes I$ for an ideal $I\subset\Lambda(n)$.
The ideal $J\cap\fI$ is, moreover, $\fr$-invariant, so $I$ is stable under all
the $\pder{\xi_k}$, thanks to~\eqref{wild:pr1}. If $I\neq 0$, take $0\neq a\in I$
of the lowest degree $d$ in the $\xi$; applying $d$ of the $\pder{\xi_k}$ we
get a~non-zero scalar in $I$, so $I=\Lambda(n)$ and~$J\supset\fI$.

Suppose now that~$J\cap\fI=0$. Then, $[J,\fI]\subset J\cap\fI=0$, so $J$
centralizes~$\fI$. Let $a+r\in J$ with $a\in\fI$ and~$r\in\fr$. Bracketing
with $x\otimes 1$ for $x\in\fh$, and using $r(1)=0$, we get $[a,x\otimes1]=0$
for all $x$; writing $a=\sum_k x_k\otimes a_k$ in the monomial basis $(a_k)$
of $\Lambda(n)$ --- a~basis of the superspace $\Lambda(n)$, the $x_k\in\fh$
being the coefficients --- we get $[x,x_k]=0$ for all $x\in\fh$, whence
$x_k=0$ as the center of a~simple Lie superalgebra is zero; so~$a=0$.
Bracketing then with $x\otimes f$ we get $x\otimes r(f)=0$ for all $f$, so
$r=0$ as a~derivation of $\Lambda(n)$, and $J=0$ --- a~contradiction.

So every non-zero ideal contains $\fI$, which is therefore the only minimal
one and is invariant under every automorphism. The quotient
$\fg_{min}(\fr)/\fI$ is $\fr$, and an isomorphism ${\fg_{min}(\fr)\tto
\fg_{min}(\fr')}$ carries $\fI$ onto $\fI'$, hence induces an isomorphism of
the quotients. \end{proof}

\begin{proof}[Proof of Theorem~\textup{\ref{wild:thmain}}] Item 1) is Lemma~\ref{wild:lemu}. In
item 2), the implications (i)$\Longleftrightarrow$(ii) are
Proposition~\ref{wild:prop}, (ii)$\Longleftrightarrow$(iii) is Lemma~\ref{wild:lemr},
and (iii)$\Longrightarrow$(iv) is clear, since isomorphic $\fr$'s that are
conjugate in $\fvect(0|n)$ give isomorphic $\fg_{min}$'s;
(iv)$\Longrightarrow$(iii) is Lemma~\ref{wild:lemg}. \end{proof}

\ssbegin[How many parameters]{Remark}[How many parameters]\label{wild:remcount}
The family is as large as the source problem: the pairs of $m\times m$
matrices modulo simultaneous conjugation depend on
\[
2m^2-(m^2-1)=m^2+1\ \text{ parameters for } m\geq2,
\]
and $m=\frac{1}{7}(n-15)$, so the subalgebras~\eqref{wild:pr1} in $\fvect(0|n)$ depend on at
least $\approx n^2/49$ parameters. Compare this with the tame problem which
the survey solves: the simple $W$-graded vectorial Lie superalgebras form
finitely many series, and their deformation parameters are counted by
a~supervariety of small dimension. This is the precise sense in which the
super analogue of \'E.~Cartan's problem differs from the problem Cartan
solved.
\end{Remark}

\ssec{Two remarks on what is \textit{not} proved}\label{wild:notproved}

\sssec{The bound $n=7m+15$ is not claimed to be optimal}
Five subspaces are needed because four are not enough: the four-subspace
problem is tame, its answer being the list of Gelfand and Ponomarev
\cite{GP}, whereas the encoding of a~pair of matrices, as in~\eqref{wild:fivesub},
needs five. The arms could be shortened at the expense of
a~permutation ambiguity: without the torus, the five lines of $\fn$ are still
determined by $\fn$ alone, being the lines of its elements of minimal non-zero
rank \textup{(}Lemma~\ref{wild:lemrank}\textup{)}, but their labeling is not, so the
correspondence becomes finite-to-one --- harmless for wildness, ugly for
an ``if and only if\rq\rq. Taking $\dim V_i=m+i$ adds $15$ dimensions and yields
the exact statement.

\sssec{Finite dimension} The Lie superalgebras $\fg_{A,B}$ are
infinite-dimensional whenever $\fh$ is, but the wild part of them, $\fr_{A,B}$,
is finite-dimensional, of $\sdim\fr_{A,B}=(11|7m+15)$. So the same
construction shows that already the \textbf{finite-dimensional} Lie
superalgebras of the form $\fvect_{-1}(\xi)\oplus\fu$, where
$\fu\subset\fgl(n)$ --- the ``linear\rq\rq\ parts of the primitive ones ---
are wild to classify. This is the form in which the statement is used in the
body of this paper, where subalgebras of $\fgl(n)$ are the daily bread, see the
Introduction and Subsection~\ref{ssOpen}: it says that ``describe all
subalgebras\rq\rq\ is hopeless, which is exactly why ``describe the
\textbf{maximal} ones\rq\rq\ is the right question.

\subsection*{Disclosures}  I.Shch is thankful to 
the Department of
Mathematics of Stockholm University  for
hospitality in 1988. 
The data used to support the findings of this study are included within the article.

\printindex


\begin{thebibliography}{9999}

\bibitem[ALSh]{ALSh}
Alekseevsky D., Leites D., Shchepochkina I., Examples of simple Lie superalgebras of vector
fields. \textit{C. R. Acad. Bulg. Sci.} \textbf{34}:9 (1980), 1187--1190 (in Russian)

\bibitem[A]{A}
Arnold V., \textit{Mathematical methods of classical mechanics}. Translated from
the 1974 Russian original by K. Vogtmann and A. Weinstein. Corrected reprint of
the second (1989) edition. Graduate Texts in Mathematics, \textbf{60}. Springer,
New York, (1997) xvi+520 pp

\bibitem[ABFM]{ABFM}
Azad H., Biswas I., Fazil A., Mahomed F., On Lie's classification of nonsolvable subalgebras of vector fields on the plane; \url{https://arxiv.org/pdf/2507.22642}

\bibitem[BLLM]{BLLM}
Bai W., Liu W., Liu X., Melikyan H., Maximal subalgebras for Lie
superalgebras of Cartan type. 
 \textit{J. Algebra Appl.} \textbf{14}:2 (2015), 1550013, 38 pp; \url{https://arxiv.org/pdf/1304.5618}

\bibitem[BLM]{BLM}
Bai W., Liu W., Melikyan H., Maximal subalgebras of Lie superalgebras of Cartan type over fields of characteristic zero. 
 \textit{J. Algebra} \textbf{404} (2014), 176--199

\bibitem[BeSe]{BeSe}
Belitskii G.~R., Sergeichuk V.~V., Complexity of matrix problems. \textit{Linear Algebra Appl.}
\textbf{361} (2003), 203--222; \url{https://arxiv.org/pdf/0709.2488}

\bibitem[BM]{BM}
Boniver F., Mathonet P., Maximal subalgebras of vector fields for
equivariant quantizations,  \textit{J. Math. Phys}. \textbf{42}:2 (2001) 582--589; \url{https://arxiv.org/pdf/math/0009239}

\bibitem[BOg]{BOg}
Borisov A. B., Ogievetskii V. I., Theory of
dynamical affine and conformal symmetries as the theory of the
gravitational field.  \textit{Teoret. Mat.
Fiz}. \textbf{21}:3 (1974),  329--342. Also reprinted in: Wess J., Ivanov E. (eds.) \textit{Supersymmetries
and quantum symmetries}. Proceedings of the International Seminar
dedicated to the memory of V. I. Ogievetsky, held in Dubna, July
22--26, 1997. Lecture Notes in Physics, \textbf{524}. Springer-Verlag,
Berlin, (1999), 400--409.

\bibitem[BGLLS]{BGLLS}
Bouarroudj S., Grozman P., Lebedev A., Leites D., Shchepochkina~I., Simple vectorial Lie algebras in characteristic
$2$ and their superizations.  \textit{SIGMA} \textbf{16} (2020), 089, 101 pages; \url{https://arxiv.org/pdf/1510.07255}

\bibitem[BGL]{BGL}
Bouarroudj S., Grozman P., Leites D., Classification of finite-dimensional modular Lie superalgebras with indecomposable Cartan
matrix.  \textit{SIGMA}, \textbf{5} (2009), 060, 63 pages;
\url{https://arxiv.org/pdf/math.RT/0710.5149}

\bibitem[BLLS]{BLLS}
Bouarroudj S., Lebedev A., Leites D., Shchepochkina~I., Classifications of simple Lie superalgebras in
characteristic $2$.  \textit{Int. Math. Res. Not. IMRN} (2023), no.~1, 54--94;
\url{https://arxiv.org/pdf/1407.1695}

\bibitem[BO]{BO}
Bouarroudj S., Ovsienko V., Riemannian curl in contact geometry.
 \textit{International Math. Res. Not.} \textbf{12} (2015), 3917--3942;
\url{https://arxiv.org/pdf/1307.2155}

\bibitem[BTW]{BTW}
Brink, L., Turbiner, A., Wyllard, N., Hidden algebras of the (super)
Calogero and Sutherland models.  \textit{J. Math. Phys}. \textbf{39}:3 (1998),  
1285--1315


\bibitem[CCLL]{CCLL}
Chapovalov D., Chapovalov M., Lebedev A.,
Leites D., The classification of almost affine (hyperbolic) Lie
superalgebras.  \textit{J. Nonlinear Math. Phys.}, \textbf{17} (2010), Special issue
1, 103--161; \url{https://arxiv.org/pdf/0906.1860}

\bibitem[CPT]{CPT}
Chapovskyi Y., Petravchuk O., Tyshchenko O., Maximal subalgebras of the Lie algebra $W_n(\Kee)$; \url{https://arxiv.org/pdf/2605.22970}



\bibitem[Dr]{Dr}
Drozd Yu.~A., Tame and wild matrix problems. In: \textit{Representation theory II}
(Ottawa, 1979), Lecture Notes in Math. \textbf{832}, Springer, Berlin (1980), 242--258


\bibitem[Dy]{Dy}
Dynkin E.B. Maximal subgroups of classical groups.  \textit{Trudy
Moskovskogo Matem. obshchestva}, (1952) \textbf{1} (Russian) = English
transl. in:  \textit{Moscow Math. Soc. Translations Ser. $2$}, \textbf{6} (1957) 245--378

\bibitem[Eld]{Eld}
Elduque A., Models of some simple modular Lie superalgebras. \textit{Pacific
J. Math.} \textbf{240}:1 (2009) 49--83;
\url{https://arxiv.org/pdf/0805.1304}

\bibitem[Eldi]{Eldi}
Elduque, A.,
Tits construction of the exceptional simple Lie algebras.  \textit{Pure Appl. Math. Q}. \textbf{7}:3 (2011),  Special Issue in honor of Jacques Tits, 559--586

\bibitem[GP]{GP}
Gelfand I.~M., Ponomarev V.~A., Problems of linear algebra and classification of quadruples of
subspaces in a~finite-dimensional vector space. In: \textit{Hilbert space operators and operator
algebras} (Tihany, 1970), Colloq. Math. Soc. J\'anos Bolyai \textbf{5}, North-Holland, Amsterdam
(1972), 163--237

\bibitem[Gel]{Gel}
Gelfand S., The Weil representation of a~Lie algebra~of type $G_2$ and related representations of $SL_3$.  \textit{Functional Anal. Appl}. \textbf{14}:1 (1980),  40--41

\bibitem[GOV]{GOV}
Gorbatsevich V. V., Onishchik A. L., Vinberg \'E. B. \textit{Foundations of Lie
theory and Lie transformation groups}. Translated from the Russian by A.
Kozlowski. Reprint of the 1993 translation [\textit{Lie groups and Lie algebras}. I,
Encyclopaedia Math. Sci., \textbf{20}, Springer, Berlin, 1993; MR 95f:22001].
Springer-Verlag, Berlin, (1997). vi+235 pp.

\bibitem[GOV1]{GOV1}
Gorbatsevich V. V., Onishchik A. L., Vinberg \'E. B. \textit{Lie groups and Lie
algebras, III. Structure of Lie groups and Lie algebras}. A~translation of
Current problems in mathematics. Fundamental directions. \textbf{41} (Russian),
Akad. Nauk SSSR, Vsesoyuz. Inst. Nauchn. i Tekhn. Inform., Moscow, 1990 [MR
91b:22001]. Translation by V. Minachin [V. V. Minakhin]. Translation edited by
A. L. Onishchik and \'E. B. Vinberg. Encyclopaedia of Mathematical Sciences, \textbf{41}.
Springer-Verlag, Berlin, (1994). iv+248 pp.

\bibitem[GL1]{GL1}
Grozman P., Leites D., Defining relations for Lie superalgebras
with Cartan matrix,  \textit{Czech. J. Phys.}, \textbf{51}:1 (2001),
1--22; \url{https://arxiv.org/pdf/hep-th/9702073}


\bibitem[He]{He}
Helgason, S., \textit{Differential Geometry, Lie Groups, and Symmetric Spaces}.
Corrected reprint of the 1978 original. Graduate Studies in Mathematics, \textbf{34}.
American Mathematical Society, Providence, RI, 2001. xxvi+641 pp.



\bibitem[Ko1]{Ko1}
Kotchetkoff [Kochetkov] Yu., D\'eformations des superalgèbres de Buttin et quantification. (French) [Deformations of Buttin superalgebras and quantization]  \textit{C. R. Acad. Sci. Paris S\'er. $I$ Math}. \textbf{299}:14 (1984),  643--645

\bibitem[KST]{KST}
Kruglikov B., Santi A., The D., G(3)-supergeometry and a~supersymmetric extension of the Hilbert-Cartan equation.  \textit{Adv. Math}. \textbf{376} (2021), 107420, 98 pp.; \url{https://arxiv.org/pdf/1908.09106}


\bibitem[L1]{L1}
Leites D., New Lie superalgebras, and mechanics.  \textit{Soviet Math. Dokl.},
\textbf{18}:5 (1977), 1277--1280

\bibitem[Lsos]{Lsos}
Leites D. (ed.), \textit{Seminar on supersymmetries $1$. Algebra
and Calculus on supermanifolds} (by J.~Bernstein, D.~Leites, V.~Molotkov, V.~Shander), MCCME,
Moscow, (2011), 410 pp. (In Russian; a~draft in English is available)


\bibitem[Ld]{Ld}
Leites D., On odd parameters in geometry. \textit{J. Lie Theory} \textbf{33}:4 (2023) 965--1004; \url{http://arxiv.org/pdf/2210.17096}


\bibitem[LeLo]{LOdd}
Leites D., Lozhechnyk O., Maximal subalgebras in three exceptional simple finite-dimensional Lie
superalgebras. (In preparation, to be put in arXiv ASAP)

\bibitem[LPS]{LPS}
Leites D., Poletaeva E., Serganova~V., On Einstein equations on manifolds and
supermanifolds,  \textit{J. Nonlinear Math. Physics}, \textbf{9}:4
(2002) 394--425; \url{https://arxiv.org/pdf/math.DG/0306209}


\bibitem[LSe]{LSe} 
Leites D., Serganova~V., Metasymmetry and Volichenko
algebras.  \textit{Phys. Lett. B} \textbf{252}:1 (1990)  91--96;

\bibitem[LSe2]{LSe2}
Leites D., Serganova~V., Symmetries wider than supersymmetries and
Volichenko algebras. In: Duplij S., Wess J. (eds.) \textit{Noncommutative Structures in Mathematics and Physics} Proceedings of
the NATO Advanced Research Workshop Kiev, Ukraine September 24--28,
2000, Kluwer, (2001), 13--30

\bibitem[LS1]{LS1}
Leites D., Shchepochkina I., The classification of simple complex Lie superalgebras of polynomial
 vector fields and their deformations; \url{https://arxiv.org/pdf/2406.16760}


\bibitem[LS3]{LS3}
Leites D., Shchepochkina I., How to quantize the antibracket.  \textit{Theor.
and Math. Physics}, \textbf{126}:3 (2001)  339--369; (In Russian. 
English translation: How should an[the] antibracket be
quantized?  \textit{Theoret. and Math. Phys.}, \textbf{126}:3 (2001),
281--306);
\url{https://arxiv.org/pdf/math-ph/0510048}


\bibitem[LieT]{LieT}
Sophus Lie's 1880 transformation group paper. In part a~translation
of ``Theorie der Trans\-for\-ma\-tions-grup\-pen" by S. Lie [Math. Ann.
{16} (1880), 441--528]. Translated by Michael Ackerman. Comments
by Robert Hermann. Lie Groups: History, Frontiers and Applications,
\textbf{I}. Math Sci Press, Brookline, Mass., (1975) xxiv+563 pp

\bibitem[LW]{LW}
Liu Wende, Wang Qi, Maximal subalgebras for modular graded Lie
superalgebras of odd Cartan type.  \textit{Transform. Groups} \textbf{20}:4 (2015),  1075--1106; \url{https://arxiv.org/pdf/1309.0171}

\bibitem[Me]{Me}
Melikyan H., Maximal subalgebras of simple modular Lie algebras.
 \textit{J.~of Algebra}, \textbf{284} (2005) 824--856

\bibitem[MSS]{MSS}
Mikhailov A. V.; Shabat A. B.; Sokolov V. V., 
The symmetry approach to classification of integrable equations. What is
integrability?  \textit{Springer Ser. Nonlinear Dynamics}, Springer,
Berlin, (1991), 115--184

\bibitem[Mo]{Mo}
Mostow G.~D., Fully reducible subgroups of algebraic groups. \textit{Amer. J. Math.} \textbf{78}
(1956), 200--221

\bibitem[Mu]{Mu}
Mulyar O.A. Maximal subalgebras of Frank algebras.  \textit{Vestnik NNGU,
ser. Math.}, \textbf{1}(3), (2005) 109--113

\bibitem[Og]{Og}
Ogievetsky V. I., Infinite-dimensional algebra~of general covariance
group as the closure of finite-dimensional algebras of conformal and
linear groups.  \textit{Lett. Nuovo Cimento}  \textbf{8}(2)
(1973), 988--990. Reprinted in: Wess J., Ivanov E. (eds.) \textit{Supersymmetries
and quantum symmetries}. Proceedings of the International Seminar
dedicated to the memory of V. I. Ogievetsky, held in Dubna, July
22--26, 1997. Lecture Notes in Physics, \textbf{524}. Springer-Verlag,
Berlin, (1999), 397--399

\bibitem[On]{On}
Onishchik A. L. \textit{Lectures on real semi-simple Lie algebras and their
representations}. Erwin Schr\"odinger Lectures in Mathematics and Physics. European Mathematical
Society (EMS), Z\"urich, (2004) x+86 pp.

\bibitem[OV]{OV}
Onishchik A. L., Vinberg \'E. B. \textit{Lie groups and algebraic groups}.
Translated from the Russian and with a~preface by D. A. Leites. Springer Series
in Soviet Mathematics. Springer-Verlag, Berlin, (1990) xx+328 pp.

\bibitem[Po1]{Po1}
Poletaeva E., Embedding of the Lie superalgebra~$D(2, 1 ; \alpha)$
into the Lie superalgebra~of pseudodifferential symbols on
$S^{1|2}$.  \textit{J. of Mathematical Physics}, \textbf{10/2007};
48(10): 103504-17; \url{https://arxiv.org/pdf/0709.0083}


\bibitem[P]{P}
Post G., A~class of graded Lie algebras of vector fields and first
order differential operators.  \textit{J. Math. Phys}. \textbf{35}:12 (1994),  
6838--6856.


\bibitem[PS]{PS}
Premet A., Stewart D., Classification of the maximal subalgebras of exceptional Lie algebras over fields of good characteristic.
 \textit{J. of the American Math. Soc}. \textbf{32}:4 (2019/20)  965--1008



\bibitem[Se]{Se}
Serganova, V. V. Classification of simple real Lie superalgebras and
symmetric superspaces.  \textit{Functional Anal.
Appl}. \textbf{17}:3 (1983),   200--207


\bibitem[Sh5]{Sh5}
Shchepochkina I., Five simple exceptional Lie superalgebras of vector
fields.  \textit{Funct. Anal. Appl.} \textbf{33}:3 (1999), 
208--219 (2000);

\bibitem[Sh14]{Sh14}
Shchepochkina I., The five exceptional simple Lie superalgebras of
vector fields and their fourteen regradings,  \textit{Represent. Theory,} \textbf{3}
(1999), 373--415 (for a~short version, see \url{https://arxiv.org/pdf/hep-th/9702121})

\bibitem[ShM]{ShM}
Shchepochkina I., Maximal subalgebras of matrix Lie superalgebras.
In: Duval C., Guieu L., Ovsienko V. (eds.) \textit{The orbit method in
geometry and physics.} In honor of A. A. Kirillov. Papers from the
International Conference held in Marseille, December 4--8,
2000. Progress in Mathematics, \textbf{213}. Birkh\"auser Boston, Inc.,
Boston, MA, (2003), 445--472. (for a~version, see \url{https://arxiv.org/pdf/1311.4131})

\bibitem[ShVf]{ShVf}
Shchepochkina I., How to realize Lie algebras by vector fields.
 \textit{Theoret. and Math. Phys.}, \textbf{147}:3 (2006)   821--838;
\url{https://arxiv.org/pdf/math.RT/0509472}

\bibitem[Shch]{Shch}
Shchepochkina I., 
Maximal Graded Solvable Subalgebras of Simple Finite-Dimensional Vectorial Lie Superalgebras. \textit{SIGMA} \textbf{21} (2025), 047; \url{https://arxiv.org/pdf/2503.03339}



\bibitem[S]{S}
Sokolov V. V., Symmetries of evolution equations.
 \textit{Russian Math. Surveys} \textbf{43}:5 (1988),   165--204;

\bibitem[SS]{SS}
Svinolupov S. I., Sokolov V. V., Representations of
contragredient Lie algebras in contact vector fields. \textit{Functional
Anal. Appl.} \textbf{25}:2 (1991),  146--147

\bibitem[The]{The}
The D., Exceptionally simple PDE.  \textit{Differential Geom. Appl.} \textbf{56} (2018), 13--41; \url{https://arxiv.org/pdf/1603.08251}

\bibitem[T1]{T1}
Turbiner A., Quasi-exactly solvable problem and $\fsl(2)$-algebra,
 \textit{Commun. Math. Phys.}, \textbf{118}, (1988) 467--474

\bibitem[W]{W}
Weisfeiler B., Infinite-dimensional filtered Lie algebras and their
connection with the graded Lie algebras. \textit{Funct. Anal. Appl.},
\textbf{2}:1 (1968), 88--89

\bibitem[Y]{Y}
Yamaguchi K., Differential systems associated with simple graded Lie
algebras.  \textit{Progress in differential geometry, Adv. Stud. Pure Math.},
\textbf{22}, Math. Soc. Japan, Tokyo, (1993) 413--494

\bibitem[ZZO]{ZZO}
Zharov O., Zolotarev B., Onishchik A., Reducible maximal subalgebra
of classical Lie superalgebras. In: Onishchik A. (ed.), \textit{Questions of
groups theory and homological algebra}, Yaroslavl Univ. press, Yaroslavl,
(1987), 161--163 (Russian)

\bibitem[Zi]{Zi}
Zirnbauer M. Riemannian symmetric superspaces and their origin in random-matrix
theory.  \textit{J. Math. Phys}. \textbf{37}:10 (1996), 4986--5018


\end{thebibliography}
\end{document}